\documentclass[12pt]{article}
\usepackage{amssymb, amsmath, amsthm, amscd}
\usepackage[dvips]{graphics}
\usepackage[utf8]{inputenc}
\usepackage[all,cmtip]{xy}
\usepackage{bbm}
\usepackage{enumitem}
\usepackage{setspace}
\usepackage[colorlinks=true]{hyperref}
\hypersetup{colorlinks   = true}
\hypersetup{linkcolor=blue}
\begin{document}

\newtheorem{The}{Theorem}[section]
\newtheorem{Lem}[The]{Lemma}
\newtheorem{Prop}[The]{Proposition}
\newtheorem{Cor}[The]{Corollary}
\newtheorem{Rem}[The]{Remark}
\newtheorem{Obs}[The]{Observation}
\newtheorem{SConj}[The]{Standard Conjecture}
\newtheorem{Titre}[The]{\!\!\!\! }
\newtheorem{Conj}[The]{Conjecture}
\newtheorem{Question}[The]{Question}
\newtheorem{Prob}[The]{Problem}
\newtheorem{Def}[The]{Definition}
\newtheorem{Not}[The]{Notation}
\newtheorem{Claim}[The]{Claim}
\newtheorem{Conc}[The]{Conclusion}
\newtheorem{Ex}[The]{Example}
\newtheorem{Fact}[The]{Fact}
\newtheorem{Formula}[The]{Formula}
\newtheorem{Formulae}[The]{Formulae}
\newtheorem{Prop-Def}[The]{Proposition and Definition}
\newtheorem{Cor-Def}[The]{Corollary and Definition}
\newcommand{\C}{\mathbb{C}}
\newcommand{\R}{\mathbb{R}}
\newcommand{\N}{\mathbb{N}}
\newcommand{\Z}{\mathbb{Z}}
\newcommand{\Q}{\mathbb{Q}}
\newcommand{\Proj}{\mathbb{P}}
\newcommand{\Rc}{\mathcal{R}}
\newcommand{\Oc}{\mathcal{O}}
\newcommand\dan[1]{{\textcolor{blue}{#1}}}
\begin{center}

{\Large\bf A Non-Integrable Ohsawa-Takegoshi-Type $L^2$ Extension Theorem}

\end{center}

\begin{center}

{\large Dan Popovici}

\end{center}

\vspace{1ex}

\noindent{\small{\bf Abstract.} Given a complete K\"ahler manifold $(X,\,\omega)$ with finite second Betti number, a smooth complex hypersurface $Y\subset X$ and a smooth real $d$-closed $(1,\,1)$-form $\alpha$ on $X$ with arbitrary, possibly non-rational, De Rham cohomology class $\{\alpha\}$ satisfying a certain assumption, we obtain extensions to $X$, with control of their $L^2$-norms, of smooth sections of the canonical bundle of $Y$ twisted by the restriction to $Y$ of any $C^\infty$ complex line bundle $L_k$ in a sequence of asymptotically holomorphic line bundles whose first Chern classes approximate the positive integer multiples $k\{\alpha\}$ of the original class. Besides a known non-integrable $(0,\,1)$-connection $\bar\partial_k$ on $L_k$, the proof uses two twisted Laplace-type elliptic differential operators that are introduced and investigated, leading to Bochner-Kodaira-Nakano-type (in-)equalities, a spectral gap result and an a priori $L^2$-estimate. The main difference from the classical Ohsawa-Takegoshi extension theorem is that the objects need not be holomorphic, but only asymptotically holomorphic as $k\to\infty$. The possibility that $\bar\partial_k$ does not square to $0$ accounts for its lack of commutation with the Laplacian $\Delta''_k$ it induces. We hope this study is a possible first step in a future attack on Siu's conjecture predicting the invariance of the plurigenera in K\"ahler families of compact complex manifolds.}

\vspace{2ex}

\section{Introduction}\label{section:Introduction} A huge literature has been generated by the seminal paper [OT87] of Ohsawa and Takegoshi asserting, for any given holomorphic function $f$ defined on the intersection $\Omega\cap H$ of a bounded pseudoconvex domain $\Omega\subset\C^n$ with a complex hyperplane $H\subset\C^n$ and satisfying the $L^2$-condition \begin{eqnarray*}\int\limits_{\Omega\cap H}|f|^2\,e^{-\psi}\,dV_{n-1}<\infty\end{eqnarray*} with respect to a plurisubharmonic (psh) weight function $\psi:\Omega\longrightarrow\R\cup\{-\infty\}$, the existence of a holomorphic function $F$ on the whole of $\Omega$ such that \begin{eqnarray*}F_{|\Omega\cap H} = f\end{eqnarray*} and \begin{eqnarray*}\int\limits_\Omega|F|^2\,e^{-\psi}\,dV_n \leq C\,\int\limits_{\Omega\cap H}|f|^2\,e^{-\psi}\,dV_{n-1},\end{eqnarray*} where, crucially, $C>0$ is a constant depending only on the diameter of $\Omega$, while $dV_{n-1}$ and $dV_n$ are the $(2n-2)$-dimensional and the $(2n)$-dimensional Lebesgue measures.

    It would be a daunting task to even enumerate the generalisations and applications of this result. Suffice to say, Ohsawa subsequently wrote a series of papers on the subject (e.g. [Ohs88], [Ohs94], [Ohs95]), while Manivel extended the result by replacing $\Omega$ with a more general K\"ahler manifold $X$, $H$ with a submanifold $Y\subset X$ of arbitrary codimension and the functions $f,F$ with holomorphic sections of a line bundle satisfying appropriate positivity properties (see [Man93]). A semiclassical version of the Ohsawa-Takegoshi extension theorem was recently obtained by Finski in [Fin21]. 

    \vspace{2ex}

    One of the numerous striking geometric applications of the Ohsawa-Takegoshi $L^2$ extension theorem was Siu's proof of the invariance of the plurigenera in projective families of compact complex manifolds, first in the case where the fibres are of general type ([Siu98]) and then in the general case ([Siu02]). Siu went on to conjecture the analogous result for K\"ahler fibres:

\begin{Conj}([Siu02, Conjecture 0.4.])\label{Conj:Siu_Kaehler_plurigenera} Let $\pi:{\cal X}\longrightarrow\Delta$ be a holomorphic family of compact {\bf K\"ahler} manifolds $X_t:=\pi^{-1}(t)$ over the open unit disc $\Delta\subset\C$. Then, for every positive integer $m$, the $m$-th plurigenus $p_m:=\mbox{dim}_\C H^0(X_t,\,mK_{X_t})$ of $X_t$ is independent of $t\in\Delta$.

\end{Conj}

One denotes by $K_Y$ the canonical bundle of a complex manifold $Y$. The case that was solved in the affirmative in [Siu98] and [Siu02] is the one where the family $\pi:{\cal X}\longrightarrow\Delta$ is supposed to be {\it projective} in the sense that there exists a positive holomorphic line bundle on its total space ${\cal X}$. (That it suffices for the fibres $X_t$ to be supposed {\it projective} was shown in [RT20].) Siu fixes a sufficiently positive holomorphic line bundle $A$ on ${\cal X}$ and first extends certain holomorphic sections of $m_0K_{X_0}+A$ from $X_0$ to ${\cal X}$ for a given $m_0$. At a later stage in his proof, he does away with the unwanted effects of the use of $A$ by twisting the $l$-th power of the section of $m_0K_{X_0}$ to be extended to ${\cal X}$ by a section of $A$ on ${\cal X}$, extends the twisted section to ${\cal X}$, takes the $l$-th root of this extension and eventually lets $l\to\infty$. These operations are based on repeated applications of a bespoke version of the Ohsawa-Takegoshi $L^2$ extension theorem that Siu derives as Theorem 3.1. in [Siu02]. 

It seems reasonable to expect that a generalisation of this approach to the K\"ahler context would aim at solving the case of Conjecture \ref{Conj:Siu_Kaehler_plurigenera} where the total space ${\cal X}$ of the family $\pi:{\cal X}\longrightarrow\Delta$ is supposed to carry a K\"ahler metric $\omega$. This $\omega$ would be a transcendental analogue of Siu's sufficiently ample line bundle $A$.

As a first step in such an undertaking, it seems legitimate to seek to prove a suitable transcendental version of the Ohsawa-Takegoshi $L^2$ extension theorem in which the holomorphic line bundle $L$ whose sections are to be extended from a submanifold $Y$ to the ambient K\"ahler manifold $X$ is replaced by a possibly non-rational De Rham cohomology class $\{\alpha\}\in H^2_{DR}(X,\,\R)$. This is what this paper sets out to achieve.

\subsection{Statement of the main result}\label{subsection:introd_statement_main} Let $X$ be a (not necessarily compact) complex manifold with $\mbox{dim}_\C X=n$ such that its {\it second Betti number} $b_2 = b_2(X)=\mbox{dim}_{\R}H^2(X,\,\R)$ is {\it finite}. We suppose there exists a {\it complete} K\"ahler metric $\omega$ on $X$. Let $\alpha$ be a ${\cal C}^\infty$ real $(1,\,1)$-form on $X$ such that $d\alpha=0$. 

The De Rham cohomology class $\{\alpha\}_{DR}$ can be viewed as an element of $\R^{b_2}$ which contains $\Q^{b_2}$ as a dense subset. Approximating $\{\alpha\}_{DR}$ by rational classes and using what we call {\it Assumption A} (see Definition \ref{Def:Assumption_A} -- this assumption is always satisfied when $X$ is either compact or the total space of a holomorphic family of compact complex manifolds), we deduce, via a well-known argument that was used at least in [Lae02] and in [Pop13], the existence of a sequence, indexed over some infinite subset $\Sigma\subset\N^\star$, of ${\cal C}^\infty$ complex line bundles $L_k\longrightarrow X$ endowed with ${\cal C}^\infty$ Hermitian metrics $h_k$ and Hermitian connections $D_k$ such that their curvature forms \begin{equation}\label{eqn:approximation_1_curvature}\frac{i}{2\pi}\Theta_{h_k}(D_k) = \alpha_k\end{equation} satisfy the condition: \begin{equation}\label{eqn:approximation_1}||\alpha_k-k\alpha||\leq_{{\cal C}^\infty}\frac{C}{k^{1/b_2}}, \hspace{3ex} k\in\Sigma\subset\N^{\star},\end{equation} where the $\alpha_k$'s are $d$-closed, ${\cal C}^\infty$ real $2$-forms on $X$ not necessarily of type $(1,\, 1)$. 

Condition (\ref{eqn:approximation_1}) means that, for every $l\in\N$ and every $\Omega\Subset X$, there exists a constant $C_l(\Omega)>0$ such that $||\alpha_k-k\alpha||_{C^l(\overline\Omega)}\leq\frac{C_l(\Omega)}{k^{1/b_2}}$ for all $k\in\Sigma\subset\N^{\star}$, where $||\,\,||_{{\cal C}^l(\overline\Omega)}$ stands for the ${\cal C}^l$-norm over the compact $\overline\Omega$. This further means that $\alpha_k$ gets closer and closer to $k\alpha$ in the ${\cal C}^\infty$ topology at the stated rate as $k\rightarrow +\infty$.

Now, the complex structure of $X$ induces a splitting $D_k=\partial_k + \bar\partial_k$ into components of types $(1,\, 0)$ and $(0,\, 1)$. Since $D_k^2 = \Theta_{h_k}(D_k)\wedge\cdot$, for every $k\in\Sigma\subset\N^{\star}$ we get: \begin{equation}\label{eqn:approximation_2}\partial_k^2 = -2\pi i\alpha_k^{2,\, 0}\wedge\cdot, \hspace{2ex} \partial_k\bar\partial_k + \bar\partial_k\partial_k = -2\pi i\alpha_k^{1,\, 1}\wedge\cdot, \hspace{2ex} \bar\partial_k^2 = -2\pi i\alpha_k^{0,\, 2}\wedge\cdot.\end{equation} Thus, if $L_k$ is not holomorphic, $\bar\partial_k^2\neq 0$ (i.e. $\bar\partial_k$ is a {\bf non-integrable} $(0,\, 1)$-connection on $L_k$) and $\alpha_k^{0,\, 2}\neq 0$.

Intuitively, the line bundles $L_k$ are {\it approximately holomorphic}, or {\it asymptotically holomorphic} as $k\to\infty$, in the sense that their holomorphicity defect is small and tends to zero. 

\begin{Def}\label{Def:approx-seq} A sequence $(L_k,\,h_k,\,D_k)_{k\in\Sigma}$ with the above properties is called an {\bf approximating sequence} for $\alpha$. 

\end{Def} 

As in the integrable case, we can define Laplacians $\Delta_k', \Delta_k'': {\cal C}^\infty_{p,\,q}(X,\,L_k)\longrightarrow {\cal C}^\infty_{p,\,q}(X,\,L_k)$ acting on smooth $L_k$-valued $(p,\,q)$-forms by $$\Delta_k'=\partial_k\partial_k^{\star} + \partial_k^{\star}\partial_k  \hspace{2ex} \mbox{and} \hspace{2ex} \Delta_k''=\bar\partial_k\bar\partial_k^{\star} + \bar\partial_k^{\star}\bar\partial_k.$$

In $\S$\ref{subsection:definitions_approx-hol}, we apply the classical theory of {\it spectral projections} to the unbounded, closed and densely defined operator $$\Delta''_k\,:\,\mbox{Dom}\,\Delta''_k\longrightarrow L^2_{p,\,q}(X,\,L_k),$$ whose domain consists in the $L_k$-valued $L^2$ $(p,\,q)$-forms $u$ on $X$ such that $\Delta''_k u$ computed in the sense of distributions is still $L^2$, to define the space ${\cal H}^{p,\,q}_{[0,\,\delta_k],\,\Delta''_k}(X,\,L_k)$ of {\bf approximately holomorphic} $L_k$-valued $(p,\,q)$-forms on $X$ (see Definition \ref{Def:Hk-def}) for a certain small constant $\delta_k>0$. In the special case when $X$ is compact, this space is the direct sum of the eigenspaces of $\Delta''_k$ corresponding to the eigenvalues lying in the interval $[0,\,\delta_k]$ and was used in [Lae02] and [Pop13].

\vspace{1ex}

The following theorem is the main result of this paper. Conclusion (b) improves the estimate obtained in conclusion (a) for the $L^2$-norm of the extension $F_k$ of $f_k\wedge dw$ to $X$ when the original $L_k$-valued $(n-1,\,0)$-form $f_k$ on $Y$ is supposed to lie in the space of {\it approximately holomorphic} sections introduced in (a) of Definition \ref{Def:Hk-def} for $(p,\,q)=(n-1,\,0)$ and $X$ replaced by $Y$.

\begin{The}\label{nonintOT} Let $(X,\,\omega)$ be a {\bf complete K\"ahler} manifold with $\mbox{dim}_{\C}X=n$ and $b_2=b_2(X)<\infty$. Let $\alpha$ be a ${\cal C}^\infty$ real $(1,\,1)$-form on $X$ such that $d\alpha=0$ and $\alpha$ satisfies {\bf Assumption A}. Choose any approximating sequence $(L_k,\,h_k,\,D_k=\partial_k+\bar\partial_k)_{k\in\Sigma}$ for $\alpha$.

 Let $Y:=\{x\in X\,;\, w(x)=0\}\subset X$ be a {\bf nonsingular complex hypersurface} of $X$ defined by a holomorphic function $w\,:\,X\longrightarrow\C$ such that $\sup_X|w|\leq 1$ and $dw(x)\neq 0$ for all $x\in Y$. 

Suppose $\alpha$ satisfies the following {\bf positivity assumption}: \begin{equation}\label{eqn:curvature-assumption}2\pi\alpha>\omega + idw\wedge d\overline{w}  \hspace{3ex} \mbox{on} \hspace{1ex} X.\end{equation}

Then, there exist constants $c_1, c_2, c_3>0$ and $C>0$ such that, for every $k\in\Sigma$, $k\gg 2$, and any constants $\varepsilon_0\in(0,\,1)$ and $0<\delta<2/b_2$, the following statements (a) and (b) hold.

\vspace{2ex}

 \noindent (a)\, For every $f=f_k\in C^\infty_{n-1,\,0}(Y,\,L_{k|Y})$ such that $\int\limits_Y|f_k|^2_{h_k,\,\omega}\,dV_{Y,\,\omega}<\infty$, there exists $F=F_k\in C^\infty_{n,\,0}(X,\,L_k)$ having the properties (i), (ii) and (iii) below:

\vspace{1ex}

(i)\, the holomorphicity defect of $F_k$ is estimated in $L^2$-norm as \begin{eqnarray*}||\bar\partial_k F_k||^2_{h_k,\,\omega}\leq \frac{12c}{k^{1+\delta}}\,\bigg(\frac{4}{c}\,k^{1+\delta} + \frac{4}{(1-\varepsilon_0)\,k}\bigg)\,\bigg(c_1\,\int\limits_Y|f|^2_{h_k,\,\omega}\,dV_{Y,\,\omega} & + & c_2\,\int\limits_Y|f|_{h_k,\,\omega}\,|\bar\partial_kf|_{h_k,\,\omega}\,dV_{Y,\,\omega} \\
& + & \frac{c_3}{4}\,\int\limits_Y|\bar\partial_kf|^2_{h_k,\,\omega}\,dV_{Y,\,\omega}\bigg);\end{eqnarray*}

\vspace{1ex}

(ii)\, $F_{k|Y} = f_k\wedge dw$; 

\vspace{1ex}

(iii)\, the squared $L^2$-norm of $F_k$ satisfies the estimate: \begin{eqnarray}\label{eqn:nonintOT_estimate}\nonumber\int\limits_X|F_k|^2_{h_k,\,\omega}\,dV_{X,\,\omega} \leq 3\cdot 2^6\bigg(\frac{k^{1+\delta}}{c} + \frac{1}{(1-\varepsilon_0)\,k}\bigg)\,\,\bigg(c_1\,\int\limits_Y|f_k|^2_{h_k,\,\omega}\,dV_{Y,\,\omega} & + & c_2\,\int\limits_Y|f_k|_{h_k,\,\omega}\,|\bar\partial_kf_k|_{h_k,\,\omega}\,dV_{Y,\,\omega} \\
    & + & \frac{c_3}{4}\,\int\limits_Y|\bar\partial_kf_k|^2_{h_k,\,\omega}\,dV_{Y,\,\omega}\bigg),\end{eqnarray} where $c:=C\,\max\bigg\{\sup\limits_X|dw|^2_\omega,\,1\bigg\}$.

\vspace{2ex}

\noindent (b)\, If $f=f_k\in {\cal H}^{n-1,\,0}_{[0,\,\delta_k],\,\Delta''_k}(Y,\,L_{k|Y})$, where $\delta_k:=\frac{a_0}{k^{1+\delta}}$ for an arbitrary constant $a_0>0$, its extension $F_k$ obtained in (a) satisfies the estimate: \begin{eqnarray}\label{eqn:nonintOT_estimate_bis}\int\limits_X|F_k|^2_{h_k,\,\omega}\,dV_{X,\,\omega} \leq 3\cdot 2^6\bigg(\frac{k^{1+\delta}}{c} + \frac{1}{(1-\varepsilon_0)\,k}\bigg)\,\,\bigg(c_1 + c_2\,\sqrt{\delta_k} + \frac{c_3}{4}\,\delta_k\bigg)\,\int\limits_Y|f_k|^2_{h_k,\,\omega}\,dV_{Y,\,\omega}.\end{eqnarray}

\end{The}

Note that the assumption $dw(x)\neq 0$ for all $x\in Y$ implies that the normal bundle $NY$ of $Y$ in $X$ is trivial. Hence, the adjunction formula reduces to $K_Y = (K_X)_{|Y}$.

Theorem \ref{nonintOT} is complementary to our earlier Ohsawa-Takegoshi-type $L^2$ extension theorem for {\it jets} of holomorphic line bundle sections obtained in [Pop05]. This time, we deal only with the sections themselves, rather than their jets, but the context is the one of a {\it non-integrable} $\bar\partial$-operator $\bar\partial_k$ as the line bundles $L_k$ and the sections involved need not be holomorphic.

Finally, we point out that in the special case where the cohomology class of $\alpha$ is {\it integral} (the classical, integrable, case of [OT87] and the other works on the subject that we are aware of), part (b) of Theorem \ref{nonintOT} follows from Theorem 1.1 of [Fin21] with a better estimate. However, the key new situation that is investigated in the present paper is the {\it non-integrable} one.

\subsection{Outline of our approach}\label{subsection:introd_outline_approach}

(I)\, Recall that the main point in the original proof of the classical {\it integrable} Ohsawa-Takegoshi theorem [OT87] was the beautiful idea to twist the $\bar\partial$-operator by an appropriately chosen function placed either on the left or the right of the operator, depending on the bidegree. (See [Siu96] for this take on the argument.) We reuse this idea by considering the sequence of vector spaces of $C^\infty$ $L_k$-valued forms on $X$ of bidegrees $(n,\,q)$ with $q=0,1,2$: \begin{equation}\label{eqn:eta-twisted-dbar-seq_n1}{\cal C}^\infty_{n,\,0}(X,\,L_k)\stackrel{\bar\partial_k^{\eta}}{\longrightarrow}{\cal C}^\infty_{n,\,1}(X,\,L_k)\stackrel{^{\eta}\bar\partial_k}{\longrightarrow}{\cal C}^\infty_{n,\,2}(X,\,L_k),\end{equation} where $\eta$ is a positive $C^\infty$ function on $X$, while the right and left twistings of $\bar\partial_k$ are the operators: \begin{eqnarray*}\bar\partial_k^{\eta}:=\bar\partial_k(\sqrt{\eta}\,\cdot) \hspace{3ex} \mbox{and} \hspace{3ex} ^{\eta}\bar\partial_k:=\sqrt{\eta}\,\bar\partial_k.\end{eqnarray*}

While the composition $^{\eta}\bar\partial_k\circ\bar\partial_k^{\eta} = \sqrt{\eta}\,\bar\partial_k^2(\sqrt{\eta}\,\cdot)$ vanishes in the integrable context of [OT87], making (\ref{eqn:eta-twisted-dbar-seq_n1}) into a complex, it need not vanish in our more general, possibly {\it non-integrable}, setup where it may happen that $\bar\partial_k^2\neq 0$.

Our starting idea is to consider the {\bf $\eta$-twisted Laplacian} \begin{equation*}\Delta''_{k,\,(\eta)}:= (\bar\partial_k^{\eta})(\bar\partial_k^{\eta})^{\star} + (^{\eta}\bar\partial_k)^{\star}(^{\eta}\bar\partial_k):{\cal C}^\infty_{n,\,1}(X,\,L_k)\longrightarrow {\cal C}^\infty_{n,\,1}(X,\,L_k)\end{equation*} in bidegree $(n,\,1)$ and the {\bf $\eta$-twisted Laplacian} \begin{equation*}\label{eqn:twistedDelta''_tilde_n0}\widetilde\Delta''_{k,\,(\eta)}:= (\bar\partial_k^{\eta})^{\star}(\bar\partial_k^{\eta}):{\cal C}^\infty_{n,\,0}(X,\,L_k)\longrightarrow {\cal C}^\infty_{n,\,0}(X,\,L_k),\end{equation*} in bidegree $(n,\,0)$, as dictated by the sequence (\ref{eqn:eta-twisted-dbar-seq_n1}).

These are both elliptic operators for which Bochner-Kodaira-Nakano-type inequalities are derived in $\S$\ref{section:twisted_laplacians} and whose spectra and the corresponding spectral projections are studied in $\S$\ref{section:spaces-approx-def}. As with $\Delta''_k$, we also consider their unbounded extensions as closed and densely defined operators on domains contained in the respective $L^2$-spaces.

\vspace{1ex}

(II)\, In $\S$\ref{subsection:OTapplication}, for every $\varepsilon>0$ small enough, we choose the same twisting function $\eta_{\varepsilon} + \lambda_{\epsilon}:X\longrightarrow(1,\,+\infty)$ as in the integrable context of [OT87]. It is defined, for an arbitrary constant $A>e$, by \begin{eqnarray*}\eta_{\varepsilon} + \lambda_{\epsilon}:= \log\frac{A}{|w|^2 + \varepsilon^2} + \frac{1}{|w|^2 + \varepsilon^2}\end{eqnarray*} and is smooth and bounded on $X$, with a {\it bump} that increases, as $\varepsilon\downarrow 0$, along the submanifold $Y\subset X$ from which the extension is to be performed.

As a consequence of the Bochner-Kodaira-Nakano-type inequalities, we obtain in (i) of Proposition \ref{Prop:BKN-consequence} the following two key inequalities under a positivity assumption that is satisfied when (\ref{eqn:curvature-assumption}) holds: \begin{eqnarray}\label{eqn:introd_l-bound_twisted-Laplacian}\Delta''_{k,\,(\eta_{\varepsilon}+\lambda_{\varepsilon})} \geq B_\varepsilon \geq  kq\,\mbox{Id} +\frac{\varepsilon^2}{(|w|^2 + \varepsilon^2)^2}\,(d\overline{w}\wedge\cdot)(d\overline{w}\wedge\cdot)^{\star}\end{eqnarray} on the space ${\cal D}^{n,\,q}(X,\,L_k)$ of compactly supported $L_k$-valued smooth $(n,\,q)$-forms. They involve the {\it twisted curvature operator} \begin{eqnarray*}B_\varepsilon:=\bigg[\eta_\varepsilon\,i\Theta(D_k)^{1,\,1} - i\partial\bar\partial\eta_\varepsilon - \frac{1}{\lambda_\varepsilon}\,i\partial\eta_\varepsilon\wedge\bar\partial\eta_\varepsilon,\,\Lambda\bigg],\end{eqnarray*} where $\Lambda=\Lambda_\omega$ is the adjoint of the multiplication operator $\omega\wedge\cdot$ w.r.t. the pointwise inner product induced by $\omega$.

\vspace{1ex}

(III)\, In $\S$\ref{section:spaces-approx-def}, we introduce our main techniques. 

$\bullet$ Firstly, using the theory of spectral projections [RS80], we propose in Definition \ref{Def:H-Q_spaces_arbitrary} two spaces ${\cal H}_{I,\,P}^{p,\,q}(X,\,L_k)$ and $Q_{I,\,P}^{p,\,q}(X,\,L_k)$ of $L_k$-valued $(p,\,q)$-forms associated with an arbitrary interval $I\subset\R$ and an arbitrary elliptic self-adjoint linear differential operator $P:\mbox{Dom}\,P\longrightarrow L^2_{p,\,q}(X,\,L_k)$ such that $P\geq 0$.


To define {\bf approximately holomorphic} $L_k$-valued $(p,\,q)$-forms, the spaces ${\cal H}_{I,\,P}^{p,\,q}(X,\,L_k)$ and $Q_{I,\,P}^{p,\,q}(X,\,L_k)$ are then specialised in Definition \ref{Def:Hk-def} to the case where $I=[0,\,\delta_k]$, for some small constant $\delta_k>0$, and $P$ is either the untwisted $\bar\partial_k$-Laplacian $\Delta''_k$ or one of its twisted counterparts $\Delta''_{k,\,(\varepsilon)}:= \Delta''_{k,\,(\eta_{\varepsilon}+\lambda_{\varepsilon})}$ and $\widetilde\Delta''_{k,\,(\varepsilon)}:= \widetilde\Delta''_{k,\,(\eta_{\varepsilon}+\lambda_{\varepsilon})}$.

\vspace{1ex}

$\bullet$ Secondly, we estimate in $\S$\ref{subsection:commutation-defect}, w.r.t. the $L^2$-norm, the {\bf commutation defect} between $\bar\partial_k$ and $\Delta''_k$ (Lemma \ref{Lem:commutation-defect_untwisted}) and then the one between the twisted counterparts of these operators (Proposition \ref{Prop:comm-defect-twisted} for an arbitrary twisting function $\eta$ and Proposition \ref{Prop:comm-defect-twisted_epsilon} for the twisting function $\eta_\varepsilon + \lambda_\varepsilon$). In the twisted case, it is key to use the twisting $\Delta''_{k,\,(\eta_\varepsilon + \lambda_\varepsilon)}$ in bidegree $(n,\,1)$ and the other twisting, $\widetilde\Delta''_{k,\,(\eta_\varepsilon + \lambda_\varepsilon)}$, in bidegree $(n,\,0)$. The actual estimate we get in the last case, for $s\in{\cal D}^{n,\,0}(X,\,L_k)$, is \begin{eqnarray*}\label{eqn:introd_comm-defect-twisted_epsilon}\nonumber\bigg|\bigg|\bigg(\Delta''_{k,\,(\eta_\varepsilon + \lambda_\varepsilon)}\bar\partial_k^{\eta_\varepsilon + \lambda_\varepsilon} & - & \bar\partial_k^{\eta_\varepsilon + \lambda_\varepsilon}\widetilde\Delta''_{k,\,(\eta_\varepsilon + \lambda_\varepsilon)}\bigg)\,s\bigg|\bigg|^2 \\
& \leq & \max\bigg\{\sup\limits_X|dw|^2_\omega,\,1\bigg\}\,\frac{1}{\varepsilon^6}\,\frac{C}{k^{\frac{2}{b_2}}}\,\bigg(\frac{k}{\varepsilon^3}\,||s||^2 + ||\bar\partial_k^{\eta_\varepsilon + \lambda_\varepsilon}s)||^2\bigg),\end{eqnarray*} for some constant $C>0$ independent of $\varepsilon$, $k$ and $s$.

\vspace{1ex}

$\bullet$ Thirdly, we get in Proposition \ref{Prop:spectra_n0-n1}, under our positivity assumption (\ref{eqn:curvature-assumption}), the following {\bf spectral gap} result in bidegree $(n,\,0)$: \begin{eqnarray*}\label{eqn:introd_spectrum_n0}\mbox{Spec}^{n,\,0}(\widetilde\Delta''_{k,\,(\varepsilon)})\bigcap\bigg(\frac{c}{k^{1+\delta}},\,(1-\varepsilon_0)k\bigg] = \emptyset\end{eqnarray*} under the crucial extra hypothesis that $\varepsilon$ is {\it constrained by $k$}. This essentially means that, for a fixed constant $0<\delta<\frac{2}{b_2}$, we have \begin{eqnarray*}\varepsilon>\frac{Const}{k^{\frac{1}{9}\,(\frac{2}{b_2}-\delta)}},\end{eqnarray*} in which case we cannot let $\varepsilon\downarrow 0$ while keeping $k$ fixed. We hope to be able, in future work, to first let $k\to\infty$  and then let $\varepsilon\downarrow 0$, or to simultaneously perform these operations by choosing $\varepsilon=k^{-a}$ for a suitable constant $a$ chosen as in Proposition \ref{Prop:spectra_n0-n1}.

\vspace{1ex}

$\bullet$ Fourthly, we derive in Theorem \ref{The:main-a-priori-estimate} our main {\bf a priori twisted $L^2$-estimate} as a non-integrable analogue of Ohsawa's a priori $L^2$ estimate of [Ohs95]. Put simply, given any $L_k$-valued $L^2$ $(n,\,0)$-form $s$ on $X$, whether $\varepsilon$ is {\it constrained by $k$}, in which case $s$ splits orthogonally as $$s = s_{h,\,\varepsilon} + s_{nh,\,\varepsilon}$$ with $s_{h,\,\varepsilon}\in{\cal H}^{n,\,0}_{[0,\,c/k^{1+\delta}],\,\widetilde\Delta''_{k,\,(\varepsilon)}}(X,\,L_k)$ and $s_{nh,\,\varepsilon}\in{\cal H}^{n,\,0}_{((1-\varepsilon_0)\,k,\,+\infty),\,\widetilde\Delta''_{k,\,(\varepsilon)}}(X,\,L_k)$, or whether $\varepsilon$ and $k$ are {\it mutually independent}, in which case $s$ splits orthogonally as $$s = s_{h,\,\varepsilon} + t_\varepsilon + s_{nh,\,\varepsilon}$$ with $t_\varepsilon\in{\cal H}^{n,\,0}_{(c/k^{1+\delta},\,(1-\varepsilon_0)\,k],\,\widetilde\Delta''_{k,\,(\varepsilon)}}(X,\,L_k)$ and $s_{h,\,\varepsilon}$, $s_{nh,\,\varepsilon}$ as above, the components of $s$ that are {\it far from being holomorphic} satisfy the $L^2$-estimates: \begin{equation*}||t_\varepsilon||^2\leq \frac{k^{1+\delta}}{c}\,||\bar\partial_k^{\eta_\varepsilon + \lambda_\varepsilon}s||^2 \hspace{3ex} \mbox{and} \hspace{3ex} ||s_{nh,\,\varepsilon}||^2\leq \frac{1}{(1-\varepsilon_0)\,k}\,||\bar\partial_k^{\eta_\varepsilon + \lambda_\varepsilon}s||^2.\end{equation*} 

\vspace{1ex}

(IV)\, Finally, we give in $\S$\ref{section:main-theorem} the proof of the main Theorem \ref{nonintOT}. Building on the ingredients introduced in $\S$\ref{section:preliminaries}--\ref{section:spaces-approx-def}, it follows the pattern of the proof of the integrable Ohsawa-Takegoshi $L^2$ extension theorem of [OT87] as it was reworked in [Siu96] and [Dem01], but with the key difference that the correction of the initial rough extension is now multiplied by $w-2\varepsilon$ rather than the standard $w$ in order to ensure the convergence, as $\varepsilon\downarrow 0$, of one of the two extra integrals that appear in our non-integrable context.

\vspace{2ex}

\noindent {\it Acknowledgements}. The author is very grateful to G. Marinescu for many useful discussions that took place during the author's several visits to Cologne at his invitation. Thanks are also due to S. Finski and S. Dinew who made very useful comments on early versions of this paper.

\section{Preliminaries: Rational approximations of transcendental classes}\label{section:preliminaries} We fix a $\Z$-basis of $H^2(X,\,\Z)$ that is also an $\R$-basis of $H^2(X,\,\R)$. This gives isomorphisms $\Z^{b_2}\simeq H^2(X,\,\Z)\subset H^2(X,\,\R)\simeq\R^{b_2}$, so we can make sense of the Euclidean norm $||\,\cdot\,||$ on $H^2(X,\,\R)\simeq\R^{b_2}$. 

 The De Rham cohomology class $\{\alpha\}\in H^2(X,\,\R)\simeq\R^{b_2}$ can be approximated by rational classes in $H^2(X,\,\Q)$ (cf. [Lae02, Th\'eor\`eme 1.3, p. 57] for the case of a compact $X$, but the argument is still valid for any $X$ with $b_2$ finite), so after clearing denominators we get a family of integral classes \begin{eqnarray}\label{eqn:alpha-approx}\gamma_k\in H^2(X,\,\Z) \hspace{1ex}\mbox{for}\hspace{1ex} k\in\Sigma\subset\N^{\star} \hspace{2ex}\mbox{such that}\hspace{1ex} ||\gamma_k-k\,\{\alpha\}||\leq\frac{C}{k^{1/b_2}},\end{eqnarray} \noindent where $\Sigma$ is an infinite subset of $\N^{\star}$ and $C = C(X)>0$ is a constant independent of $k$, depending only on $b_2$, hence on $X$. Note that the family $(\gamma_k)_{k\in\Sigma}$ is neither unique, nor canonically defined; we simply make an arbitrary choice.

 Recall that the exponential exact sequence of sheaves $0\rightarrow\underline{\Z}\rightarrow {\cal E}\rightarrow {\cal E}^{\star}\rightarrow 0$ induces a long exact sequence in cohomology containing the sequence

$$0=H^1(X,\,{\cal E})\longrightarrow H^1(X,\,{\cal E}^{\star})\stackrel{\stackrel{c_1}{\simeq}}{\longrightarrow}H^2(X,\,\Z)\longrightarrow H^2(X,\,{\cal E})=0,$$

\noindent where the map $c_1$ associates with any isomorphism class of ${\cal C}^\infty$ line bundles on $X$ its first Chern class, while $H^q(X,\,{\cal E})$ vanishes for all $q\geq 1$ since the sheaf ${\cal E}$ of germs of ${\cal C}^\infty$ functions on $X$ is acyclic (on any paracompact, possibly non-compact, differentiable manifold $X$). The map $c_1: H^1(X,\,{\cal E}^{\star})\rightarrow H^2(X,\,\Z)$ is thus an isomorphism, hence

\begin{equation}\label{eqn:L_k-def}\forall\, k\in\Sigma, \hspace{1ex} \exists\,!\hspace{1ex} L_k\longrightarrow X \hspace{1ex} {\cal C}^\infty\hspace{1ex}\mbox{complex line bundle such that} \hspace{1ex} c_1(L_k)=\gamma_k.\end{equation}

 Fix now a Hermitian metric $\omega$ on $X$. We will suppose that representatives of the approximating integral classes $\gamma_k$ can be found to approximate multiples of the original real $d$-closed $(1,\,1)$-form $\alpha$.

\vspace{2ex}

\begin{Def}\label{Def:Assumption_A} We say that $\alpha$ satisfies {\bf Assumption A} if there exist ${\cal C}^\infty$ real $2$-forms $\alpha_k$ on $X$ such that: 

\vspace{1ex}

$(i)$\, $d\alpha_k=0$ for all $k\in\Sigma\subset\N^{\star}$;

$(ii$\, $\gamma_k=\{\alpha_k\}\in H^2(X,\,\Z)$ for all $k\in\Sigma\subset\N^{\star}$;

$(iii)$\, there is a Hermitian metric $\omega$ on $X$ such that \begin{equation}\label{not:Cinfty_norm}||\alpha_k-k\alpha||\leq_{{\cal C}^\infty}\frac{C}{k^{1/b_2}} \hspace{5ex} \mbox{for all} \hspace{2ex} k\in\Sigma\subset\N^{\star}.\end{equation}

\end{Def} 

\vspace{2ex}

For the meaning of the symbol $\leq_{{\cal C}^\infty}$, see explanation after (\ref{eqn:approximation_1}) in the introduction. Note that if (\ref{not:Cinfty_norm}) ($=$ (\ref{eqn:approximation_1})) is satisfied for some metric $\omega$, it is satisfied for any Hermitian metric on $X$ since any two such metrics are comparable on any compact $\overline\Omega\subset X$.

\vspace{2ex}

\noindent {\bf Consequence of Assumption A.} If $\alpha_k=\alpha_k^{2,\,0} + \alpha_k^{1,\,1} + \alpha_k^{0,\,2}$ is the splitting of $\alpha_k$ into forms of pure type, Assumption $A$ implies the following estimates: \begin{equation}\label{eqn:approximation_3} (i)\, ||\alpha_k^{1,\,1}-k\alpha||\leq_{{\cal C}^\infty}\frac{C}{k^{1/b_2}} \hspace{3ex} \mbox{and} \hspace{3ex} (ii)\, ||\alpha_k^{0,\,2}||\leq_{{\cal C}^\infty}\frac{C}{k^{1/b_2}} \end{equation}

\noindent for all $k\in\Sigma\subset\N^{\star}$, where the constants $C>0$ depend on the open subset $\Omega\Subset X$ on which $\alpha$ and its pure-type components are evaluated.

\vspace{2ex}

 We note two important cases where this assumption is satisfied.

 \begin{Lem}\label{Lem:Assumption_A_satisfied} $(1)$\, Assumption $(A)$ is satisfied if $X$ is compact (as was shown in [Lae02, p.57]).

\vspace{1ex}

$(2)$\, On the total space ${\cal X}$ of any holomorphic family of compact complex manifolds $\pi:{\cal X}\longrightarrow S$ over a contractible base $S$, Assumption $(A)$ is satisfied with the same constants on the whole of $S$.

\vspace{1ex}

(Such a $\pi$ is, by definition, a proper holomorphic submersion from a complex manifold ${\cal X}$ to a ball $S\subset\C^m$ containing the origin.)

\end{Lem}

\noindent {\it Proof.} $(1)$\, If $X$ is compact, ${\cal C}^\infty_2(X,\,\C)$ is the space of smooth complex-valued $2$-forms on $X$ and $\omega$ is a Hermitian metric on $X$, we have orthogonal decompositions

\vspace{1ex}

\hspace{3ex} ${\cal C}^\infty_2(X,\,\C)=\ker\Delta\oplus\mbox{Im}\,d\oplus\mbox{Im}\,d^{\star} \hspace{2ex} \mbox{and} \hspace{2ex} \ker d = \ker\Delta\oplus\mbox{Im}\,d$

\vspace{1ex}

\noindent inducing a splitting $\ker d\ni\alpha = h + dv$ with $h\in\ker\Delta$ and $v$ a $1$-form. If $h_k$ is the $\Delta$-harmonic representative of the class $\gamma_k$, we set (cf. [Lae02, p. 57])

\vspace{1ex}

\hspace{16ex} $\alpha_k:=h_k + kdv, \hspace{2ex} k\in\Sigma\subset\N^{\star}.$

\vspace{1ex}

$(2)$\, Since $S$ is contractible and the family is ${\cal C}^\infty$ trivial (i.e. the fibres $X_t:=\pi^{-1}(t)$ are all ${\cal C}^\infty$ diffeomorphic to a fixed ${\cal C}^\infty$ manifold $X$ thanks to Ehresmann's Theorem), we have canonical isomorphisms

\vspace{1ex}

\hspace{10ex}   $H^2({\cal X},\,\R)\simeq H^2(X_t,\,\R)\simeq H^2(X,\,\R), \hspace{2ex} t\in S.$

\vspace{1ex}

\noindent Thus $\mbox{dim}_{\R}H^2({\cal X},\,\R)<+\infty$ since the fibres $X_t$ are compact. 

We fix a Riemannian metric $g$ on $X$ and we denote by $\Delta:=dd^\star + d^\star d$ the $d$-Laplacian induced by $g$ on the differential forms of $X$, where $d^\star$ is the formal adjoint of $d$ w.r.t. the $L^2$ inner product defined by $g$. If ${\cal C}^\infty_2(X,\,\C)$ stands for the space of smooth complex-valued $2$-forms on $X$, standard Hodge Theory yields orthogonal decompositions: $${\cal C}^\infty_2(X,\,\C)=\ker\Delta\oplus\mbox{Im}\,d\oplus\mbox{Im}\,d^{\star} \hspace{5ex} \mbox{and} \hspace{5ex} \ker d = \ker\Delta\oplus\mbox{Im}\,d$$ inducing the Hodge isomorphism $\ker\Delta\simeq H^2(X,\,\C)$.

If we are given a ${\cal C}^\infty$ real $(1,\,1)$-form $\alpha$ on ${\cal X}$ such that $d\alpha=0$, we let $\gamma_k\in H^2({\cal X},\,\Z)$ (with $k\in\Sigma\subset\N^\star$) be a choice of integral classes such that $||\gamma_k - k\,\{\alpha\}||\leq C/k^{1/b_2}$ for all $k\in\Sigma$ (cf. (\ref{eqn:alpha-approx})). The classes $\{\alpha\}$ and $\gamma_k$ can be viewed as elements in $H^2(X,\,\C)$. We will choose representatives $\alpha_k\in\gamma_k$ in the following way.

Let $h$, resp. $h_k$, be the $\Delta$-harmonic representative of the class $\{\alpha\}\in H^2(X,\,\C)$, resp. $\gamma_k\in H^2(X,\,\C)$. The restricted form $\alpha_{|X_t}\in\ker d$ splits as $\alpha_{|X_t} = h + dv_t$, with $v_t$ a smooth $1$-form on $X$ depending on $t$. We can choose the forms $v_t$ to depend in a ${\cal C}^\infty$ way on $t$. Put $\alpha_{k|X_t}:=h_k + dv_t$ for every $k$ and every $t$. Thus, $\alpha_{k|X_t}$ is a smooth form on $X$ lying in the class $\gamma_k\in H^2(X,\,\C)$. The forms $\alpha_{k|X_t}$ depend in a ${\cal C}^\infty$ way on $t$, so they define a smooth $2$-form on ${\cal X}$ that we denote by $\alpha_k$. 

 For any given $l\in\N$, if we equip ${\cal H}^2_\Delta(X,\,\C)$ with the $C^l$-norm induced by the $C^l$-norm on the space ${\cal C}^\infty_2(X,\,\C)$ of smooth $2$-forms on $X$, the norm induced under the Hodge isomorphism $T:{\cal H}^2_\Delta(X,\,\C)\longrightarrow H^2(X,\,\C)$ on the finite-dimensional space $H^2(X,\,\C)$ is equivalent to the Euclidean norm defined by a fixed basis of $H^2(X,\,\C)$. The continuity of the inverse map $T^{-1}$ w.r.t. these norms implies the existence of a constant $C_l>0$ such that $||v||_{C^l}\leq C_l\,||T v||$ for all $v\in{\cal H}^2_\Delta(X,\,\C)$. In particular, for every $l\in\N$ we get the first inequality below: \begin{eqnarray*}||(\alpha_k - k\alpha)_{|X_t}||_{C^l(X)} =||h_k - kh||_{C^l(X)} \leq C_l\,||\gamma_k - k\,\{\alpha\}|| \leq \frac{C\,C_l}{k^{\frac{1}{b_2}}},\end{eqnarray*} for all $t\in S$ and all $k\in\Sigma$.

 Now, the forms $h_k$ and $h$ are independent of $t$, by construction, so they can be viewed as $2$-forms on ${\cal X}$ that are constant w.r.t. $t$. Hence, the derivatives in the $t$-directions of $h_k - kh$ vanish and we get

\vspace{1ex}

\hspace{6ex} $||\alpha_k - k\alpha||_{C^l({\cal X})} =||h_k - kh||_{C^l({\cal X})} \leq C_l\,||\gamma_k - k\,\{\alpha\}|| \leq \frac{C\,C_l}{k^{\frac{1}{b_2}}}.$

\vspace{1ex}

\noindent for every $l\in\N$.     

 Thus, the forms $(\alpha_k)_{k\in\Sigma}$ satisfy properties $(i)-(iii)$ in Assumption (A).   \hfill $\Box$

\section{Twisted Laplacians for non-integrable connections}\label{section:twisted_laplacians}

Let $X$ be a complex manifold ($\mbox{dim}_{\C}X=n$) endowed with a {\it complete} K\"ahler metric $\omega$. Suppose that $b_2$ is finite. Let $\alpha$ be a ${\cal C}^\infty$ real $(1,\, 1)$-form on $X$ such that $d\alpha=0$ and $\alpha$ satisfies Assumption A. Let $(L_k,\,h_k,\,D_k)_{k\in\Sigma}$ be a choice of {\bf approximating sequence} for $\alpha$ (see Definition \ref{Def:approx-seq}).

Fix $k\in\Sigma\subset\N^{\star}$ and a bidegree $(p,\,q)$. It is standard that the operator $\Delta''_k\,:\,{\cal C}^\infty_{p,\,q}(X,\,L_k)\longrightarrow {\cal C}^\infty_{p,\,q}(X,\,L_k)$ extends to an unbounded, closed and densely defined operator $$\Delta''_k\,:\,\mbox{Dom}\,\Delta''_k\longrightarrow L^2_{p,\,q}(X,\,L_k)$$ where the domain $\mbox{Dom}\,\Delta''_k\subset L^2_{p,\,q}(X,\,L_k)$ of $\Delta''_k$ consists of the forms $u\in L^2_{p,\,q}(X,\,L_k)$ for which $\Delta''_ku$ computed in the sense of distributions still belongs to $L^2_{p,\,q}(X,\,L_k)$. Similar extensions are defined for $\bar\partial_k$ and $\bar\partial_k^{\star}$. The {\it completeness} assumption on $\omega$ ensures (see e.g. [Dem97, chapter VIII]) that the Von Neumann adjoint of $\bar\partial_k$ coincides with the unbounded extension of the formal adjoint $\bar\partial_k^{\star}$ and extends the following identity from ${\cal D}_{p,\,q}(X,\,L_k)$ to $\mbox{Dom}\,\Delta''_k$: \begin{equation}\label{eqn:Delta''-scalarprod}\langle\langle\Delta''_ku,\,u\rangle\rangle = ||\bar\partial_ku||^2 + ||\bar\partial_k^{\star}u||^2  \hspace{2ex}\mbox{for all}\hspace{1ex} u\in\mbox{Dom}\,\Delta''_k\subset L^2_{p,\,q}(X,\,L_k)\end{equation}

\noindent So, in particular, $\mbox{Dom}\,\Delta''_k\subset\mbox{Dom}\,\bar\partial_k\cap\mbox{Dom}\,\bar\partial_k^{\star}$ and the unbounded operator $\Delta''_k$ is self-adjoint. 

These facts are still valid if $\Delta''_k$ is replaced by any elliptic self-adjoint non-negative operator defined in a similar way as $AA^\star + B^\star B$ for some operators $A, B$ that replace $\bar\partial_k$. In particular, they remain valid for all the Laplacians we use in this work, including $\Delta'_k$, $\Delta''_k$ and their twisted analogues introduced below. We will often write (in)equalities for {\it compactly supported} $L_k$-valued smooth forms $u\in{\cal D}^{p,\,q}(X,\,L_k)$ and we will tacitly use their validity throughout the domain of the Laplacian involved. Indeed, the completeness of $\omega$ ensures the density of ${\cal D}^{p,\,q}(X,\,L_k)$ in this domain w.r.t. the graph norm (see e.g. [Dem97, chapter VIII]).

\subsection{Definition of twisted Laplacians}\label{subsection:twisted-lap-def}

 We will need to twist the operators $\partial_k$ and $\bar\partial_k$ by a well-chosen positive smooth function as in the proof of the classical Ohsawa-Takegoshi theorem (cf. point of view taken in [Siu96]). For the sake of generality, we let $\eta>0$ be an arbitrary ${\cal C}^\infty$ positive function on $X$ at this point. A special choice of $\eta$ will be made later on. 

Fix a bidegree $(p,\,q)$ and consider for every $k$ the sequence of vector spaces \begin{eqnarray}\label{eqn:eta-twisted-dbar-seq}{\cal C}^\infty_{p,\,q-1}(X,\,L_k)\stackrel{\bar\partial_k^{\eta}}{\longrightarrow}{\cal C}^\infty_{p,\,q}(X,\,L_k)\stackrel{^{\eta}\bar\partial_k}{\longrightarrow}{\cal C}^\infty_{p,\,q+1}(X,\,L_k)\end{eqnarray} where we define: \begin{eqnarray}\label{eqn:twisted_d-bar_def}\bar\partial_k^{\eta}:=\bar\partial_k(\sqrt{\eta}\,\cdot) \hspace{3ex} \mbox{and} \hspace{3ex} ^{\eta}\bar\partial_k:=\sqrt{\eta}\,\bar\partial_k.\end{eqnarray} This sequence is not a complex if $\bar\partial_k^2\neq 0$. The metrics $\omega$ of $X$ and $h_k$ of $L_k$ define $L^2$ inner products on the spaces ${\cal C}^\infty_{p,\,q}(X,\,L_k)$ with respect to which we can consider the formal adjoints of the above twisted operators

$${\cal C}^\infty_{p,\,q+1}(X,\,L_k)\stackrel{(^{\eta}\bar\partial_k)^{\star}}{\longrightarrow}{\cal C}^\infty_{p,\,q}(X,\,L_k)\stackrel{(\bar\partial_k^{\eta})^{\star}}{\longrightarrow}{\cal C}^\infty_{p,\,q-1}(X,\,L_k)$$

\noindent given explicitly by $(\bar\partial_k^{\eta})^{\star}=\sqrt{\eta}\,\bar\partial_k^{\star}$ and $(^{\eta}\bar\partial_k)^{\star}=\bar\partial_k^{\star}(\sqrt{\eta}\,\cdot)$. 

\vspace{2ex}

The main new tool that we introduce and study in this work is the operator induced by (\ref{eqn:eta-twisted-dbar-seq}) in two incarnations. 

\subsubsection{The first form of the twisted Laplacian}\label{subsub:1st-form_twistedL} Having fixed a bidegree $(p,\,q)$, the twisted $\bar\partial$-operators of the sequence (\ref{eqn:eta-twisted-dbar-seq}) naturally induce the twisted Laplacian described below acting on the middle space of the sequence.

\begin{Def}\label{Def:twisted-d-bar-laplacian} Let $X$ be a complex manifold equipped with a complete K\"ahler metric $\omega$ such that $b_2$ is finite. Let $(L_k,\,h_k,\,D_k = \partial_k + \bar\partial_k)\longrightarrow X$ be a choice of approximately holomorphic ${\cal C}^\infty$ Hermitian complex line bundles associated with a given closed ${\cal C}^\infty$ real $(1,\,1)$-form $\alpha$ and let $\eta>0$ be a ${\cal C}^\infty$ positive function on $X$. 

 For every bidegree $(p,\,q)$, we define the $\eta$-{\bf twisted $\bar\partial$-Laplacian} by \begin{equation}\label{eqn:twistedDelta''def}\Delta''_{k,\,(\eta)}:= (\bar\partial_k^{\eta})(\bar\partial_k^{\eta})^{\star} + (^{\eta}\bar\partial_k)^{\star}(^{\eta}\bar\partial_k):{\cal C}^\infty_{p,\,q}(X,\,L_k)\longrightarrow {\cal C}^\infty_{p,\,q}(X,\,L_k),\end{equation}

\noindent where the formal adjoints are computed w.r.t. $\omega$ and $h_k$.

\end{Def}

  Thus, $\Delta''_{k,\,(\eta)}\geq 0$ and $(\Delta''_{k,\,(\eta)})^{\star} = \Delta''_{k,\,(\eta)}$, while a trivial computation shows that \begin{equation}\label{eqn:twisted-lap''-formula}\Delta''_{k,\,(\eta)} = \bar\partial_k\eta\bar\partial_k^{\star} + \bar\partial_k^{\star}\eta\bar\partial_k,\end{equation} where the symbol $\eta$ on the right denotes the operator of multiplication by the function $\eta$.

 Hence, for all compactly supported $v\in {\cal C}^\infty_{p,\,q}(X,\,L_k)$ we have

\begin{equation}\label{eqn:twisted-lap''-scal}\langle\langle\Delta''_{k,\,(\eta)}v,\, v\rangle\rangle = \langle\langle\eta\,\bar\partial_kv,\, \bar\partial_kv\rangle\rangle +\langle\langle\eta\,\bar\partial_k^{\star}v,\, \bar\partial_k^{\star}v\rangle\rangle = ||\sqrt{\eta}\,\bar\partial_kv||^2 + ||\sqrt{\eta}\,\bar\partial_k^{\star}v||^2.\end{equation}


 Running the same construction with $\partial_k$ in place of $\bar\partial_k$, we get sequences

\begin{equation}\label{eqn:eta-twisted-del-seq}{\cal C}^\infty_{p-1,\,q}(X,\,L_k)\stackrel{\partial_k^{\eta}}{\longrightarrow}{\cal C}^\infty_{p,\,q}(X,\,L_k)\stackrel{^{\eta}\partial_k}{\longrightarrow}{\cal C}^\infty_{p+1,\,q}(X,\,L_k)\end{equation}

\noindent in which we define: \begin{eqnarray*}\label{eqn:twisted_del_def}\partial_k^{\eta}:=\partial_k(\sqrt{\eta}\,\cdot), \hspace{3ex} ^{\eta}\partial_k:=\sqrt{\eta}\,\partial_k\end{eqnarray*} and \begin{eqnarray*}{\cal C}^\infty_{p+1,\,q}(X,\,L_k)\stackrel{(^{\eta}\partial_k)^{\star}}{\longrightarrow}{\cal C}^\infty_{p,\,q}(X,\,L_k)\stackrel{(\partial_k^{\eta})^{\star}}{\longrightarrow}{\cal C}^\infty_{p-1,\,q}(X,\,L_k)\end{eqnarray*}

\noindent where $(\partial_k^{\eta})^{\star}=\sqrt{\eta}\,\partial_k^{\star}$ and $(^{\eta}\partial_k)^{\star}=\partial_k^{\star}(\sqrt{\eta}\,\cdot)$. We also get the following

\begin{Def}\label{Def:twisted-d-laplacian} Under the assumptions of Definition \ref{Def:twisted-d-bar-laplacian}, for every bidegree $(p,\,q)$ we define the $\eta$-{\bf twisted $\partial$-Laplacian} by \begin{equation}\label{eqn:twistedDelta'def}\Delta'_{k,\,(\eta)}:= (\partial_k^{\eta})(\partial_k^{\eta})^{\star} + (^{\eta}\partial_k)^{\star}(^{\eta}\partial_k):{\cal C}^\infty_{p,\,q}(X,\,L_k)\longrightarrow {\cal C}^\infty_{p,\,q}(X,\,L_k).\end{equation}
\end{Def}

\vspace{2ex}

 We have $\Delta'_{k,\,(\eta)}\geq 0$, $(\Delta'_{k,\,(\eta)})^{\star} = \Delta'_{k,\,(\eta)}$ and 

\begin{equation}\label{eqn:twisted-lap'-formula}\Delta'_{k,\,(\eta)} = \partial_k\eta\partial_k^{\star} + \partial_k^{\star}\eta\partial_k.\end{equation}

\noindent Consequently,

\begin{equation}\label{eqn:twisted-lap'-scal}\langle\langle\Delta'_{k,\,(\eta)}v,\, v\rangle\rangle = ||\sqrt{\eta}\,\partial_kv||^2 + ||\sqrt{\eta}\,\partial_k^{\star}v||^2, \hspace{3ex} v\in {\cal D}^{p,\,q}(X,\,L_k).\end{equation}


Note that the operators $\Delta''_k, \Delta'_k, \Delta''_{k,\,(\eta)}$ and $\Delta'_{k,\,(\eta)}$ are {\it self-adjoint} since they are formally so and the metric $\omega$ is complete.

\subsubsection{The second form of the twisted Laplacian}\label{subsub:2nd-form_twistedL} For a reason that will become apparent later, we will also use another twisting $\widetilde\Delta''_{k,\,(\eta)}$ of the Laplacian $\Delta''_k$. It differs from the original twisting $\Delta''_{k,\,(\eta)}$ in that the position of the twisting function $\eta$ with respect to the operator $\bar\partial_k$ is reversed. 


\begin{Def}\label{Def:twisted-d-bar-laplacian_tilde} In the context of Definition \ref{Def:twisted-d-bar-laplacian}, for every bidegree $(p,\,q)$, we define the following operator: \begin{equation}\label{eqn:twistedDelta''_tilde_def}\widetilde\Delta''_{k,\,(\eta)}:= (^{\eta}\bar\partial_k)(^{\eta}\bar\partial_k)^{\star} + (\bar\partial_k^{\eta})^{\star}(\bar\partial_k^{\eta}):{\cal C}^\infty_{p,\,q}(X,\,L_k)\longrightarrow {\cal C}^\infty_{p,\,q}(X,\,L_k),\end{equation}

\noindent where the formal adjoints are computed w.r.t. $\omega$ and $h_k$.

\end{Def}

The motivation behind Definition \ref{Def:twisted-d-bar-laplacian_tilde} was explained in the introduction. We shall now observe a link between the kernels and the spectra of the untwisted
Laplacian $\Delta''_k:=\bar\partial_k\bar\partial_k^{\star} + \bar\partial_k^{\star}\bar\partial_k:\,C^\infty_{p,\,q}(X,\,L_k)\rightarrow C^\infty_{p,\,q}(X,\,L_k)$ and those of the twisted Laplacian $\widetilde\Delta''_{k,\,(\eta)}:C^\infty_{p,\,q}(X,\,L_k)\longrightarrow C^\infty_{p,\,q}(X,\,L_k)$.

\begin{Lem}\label{Lem:link_Laplacians_ker-spectra} In the context of Definitions \ref{Def:twisted-d-bar-laplacian} and \ref{Def:twisted-d-bar-laplacian_tilde}, for every bidegree $(p,\,q)$ and for every constant $\lambda\geq 0$, the following equivalence holds for every $u\in{\cal D}^{p,\,q}(X,\,L_k)$: \begin{equation}\label{eqn:equiv_Laplacians_ker-spectra}\widetilde\Delta''_{k,\,(\eta)}u = \lambda\,u \iff \Delta''_k(\sqrt{\eta}\,u) = \frac{\lambda}{\eta}\,(\sqrt{\eta}\,u).\end{equation}

In particular, the linear map: \begin{equation}\label{eqn:map_Laplacians_ker}\ker\widetilde\Delta''_{k,\,(\eta)}\ni u\longmapsto \sqrt{\eta}\,u\in\ker\Delta''_k =\ker\Delta''_{k,\,(\eta)}\end{equation} is an isomorphism.

\end{Lem}

\noindent {\it Proof.} Equivalence (\ref{eqn:equiv_Laplacians_ker-spectra}) follows from the following sequence of equivalences: \begin{eqnarray*}\widetilde\Delta''_{k,\,(\eta)}u = \lambda\,u & \iff & \langle\langle\widetilde\Delta''_{k,\,(\eta)}u,\,v\rangle\rangle = \langle\langle\lambda u,\,v\rangle\rangle \hspace{2ex} \forall\,v\in L^2_{p,\,q}(X,\,L_k) \\
 & \iff & \langle\langle\bar\partial_k^{\eta}u,\,\bar\partial_k^{\eta}v\rangle\rangle + \langle\langle(^{\eta}\bar\partial_k)^{\star}u,\, (^{\eta}\bar\partial_k)^{\star}v\rangle\rangle = \langle\langle\lambda u,\,v\rangle\rangle \hspace{2ex} \forall\,v\in {\cal D}^{p,\,q}(X,\,L_k) \\  
 & \iff & \langle\langle\bar\partial_k(\sqrt{\eta}\,u),\,\bar\partial_k(\sqrt{\eta}\,v)\rangle\rangle + \langle\langle(\bar\partial_k^\star(\sqrt{\eta}\,u),\, \bar\partial_k^\star(\sqrt{\eta}\, v)\rangle\rangle = \langle\langle\lambda u,\,v\rangle\rangle \hspace{2ex} \forall\,v\in {\cal D}^{p,\,q}(X,\,L_k) \\
 & \iff & \langle\langle\Delta''_k(\sqrt{\eta}\,u),\,\sqrt{\eta}\,v\rangle\rangle = \langle\langle(\lambda/\eta)\,\sqrt{\eta}\,u,\,\sqrt{\eta}\,v\rangle\rangle \hspace{2ex} \forall\,v\in {\cal D}^{p,\,q}(X,\,L_k).\end{eqnarray*}

The equality between the kernels of $\Delta''_k$ and $\Delta''_{k,\,(\eta)}$ follows by noticing that \begin{eqnarray*}\ker\Delta''_{k,\,(\eta)} = \ker(^{\eta}\bar\partial_k)\cap\ker(\bar\partial_k^{\eta})^\star = \ker(\sqrt{\eta}\,\bar\partial_k)\cap\ker(\sqrt{\eta}\,\bar\partial_k^\star) = \ker\bar\partial_k\cap\ker\bar\partial_k^\star = \ker\Delta''_k.\end{eqnarray*}

The proof is complete. \hfill $\Box$

\vspace{3ex}

Note that the above proof even gives the following equality: \begin{equation}\label{eqn:map_Laplacians_ker_proof-extra}\langle\langle\widetilde\Delta''_{k,\,(\eta)}u,\,u\rangle\rangle = \langle\langle\Delta''_k(\sqrt{\eta}\,u),\,\sqrt{\eta}\,u\rangle\rangle\end{equation} for every $u\in {\cal D}^{p,\,q}(X,\,L_k)$.

\subsection{The special integrable case} 

 Although we shall be working on a possibly non-compact complex manifold $X$ with a possibly non-integrable $(0,\,1)$-connection $\bar\partial_k$ on $L_k$, we pause to point out that the standard $L^2$ theory still applies to the twisted $\bar\partial_k$-operators in the special integrable case on a compact manifold. This special framework will serve later on as a model that will be modified in our more general context.

\begin{Obs}\label{Obs:compact-integrable-special} Suppose $X$ is compact and $\bar\partial_k^2=0$. Fix any $(p,\,q)$.

\vspace{1ex}

\noindent $(i)$\, There is an orthogonal three-space decomposition

\begin{equation}\label{eqn:3-space-decomp_integrable}{\cal C}^\infty_{p,\,q}(X,\,L_k) = \ker\Delta''_{k,\,(\eta)}\oplus\mbox{Im}\,\bar\partial_k^{\eta}\oplus\mbox{Im}\,(^{\eta}\bar\partial_k)^{\star}\end{equation}

\noindent in which $\ker\,(^{\eta}\bar\partial_k) = \ker\Delta''_{k,\,(\eta)}\oplus\mbox{Im}\,\bar\partial_k^{\eta}$ and $\ker\,(\bar\partial_k^{\eta})^{\star}=\ker\Delta''_{k,\,(\eta)}\oplus\mbox{Im}\,(^{\eta}\bar\partial_k)^{\star}$.

\vspace{1ex}

\noindent $(ii)$\, For every $v\in\mbox{Im}\,\bar\partial_k^{\eta}\subset {\cal C}^\infty_{p,\,q}(X,\,L_k)$, the minimal $L^2$-norm solution $u\in {\cal C}^\infty_{p,\,q-1}(X,\,L_k)$ of the equation \begin{equation}\label{eqn:dbar-eq-integrable}\bar\partial_k^{\eta}u=v\end{equation}

\noindent is given explicitly by the Neumann-type formula: \begin{equation}\label{eqn:dbar-eq-integrable-sol}u=(\bar\partial_k^{\eta})^{\star}\Delta^{''-1}_{k,\,(\eta)}v,\end{equation}

\noindent where $\Delta^{''-1}_{k,\,(\eta)}$ is the Green operator of $\Delta''_{k,\,(\eta)}$. 

 Note that equation (\ref{eqn:dbar-eq-integrable}) is equivalent to the equation $\bar\partial_k(\sqrt{\eta}\,u) =v$, so its solutions are obtained from those of the untwisted equation $\bar\partial_k u =v$ by a division by $\sqrt{\eta}$. However, what changes in the twisted equation is the $L^2$ estimate of the minimal $L^2$-norm solution (see $(iii)$ below). Indeed, if $u$ is the minimal $L^2$-norm solution of the untwisted equation $\bar\partial_k u =v$, then $(1/\sqrt{\eta})\,u$, though a solution, need not be the minimal $L^2$-norm solution of the twisted equation $\bar\partial_k^{\eta}u=v$.

\vspace{1ex}

\noindent $(iii)$\, In practice, the Bochner-Kodaira-Nakano inequality (see next subsection) yields an operator of order zero $B_{k,\,(\eta)}:{\cal C}^\infty_{p,\,q-1}(X,\,L_k)\longrightarrow {\cal C}^\infty_{p,\,q-1}(X,\,L_k)$ depending on the curvature of $(L_k,\, D_k,\, h_k)$ such that $\Delta''_{k,\,(\eta)}\geq B_{k,\,(\eta)}$. If, under suitable positivity assumptions on $L_k$, the operator $B_{k,\,(\eta)}$ has a positive lower bound (i.e. if $B_{k,\,(\eta)}\geq c\,\mbox{Id}$ everywhere on $\Lambda^{p,\,q-1}T^\star X\otimes L_k$ for a constant $c>0$), then the minimal $L^2$-norm solution $u\in {\cal C}^\infty_{p,\,q-1}(X,\,L_k)$ of equation (\ref{eqn:dbar-eq-integrable}) satisfies the $L^2$ estimate

\begin{equation}\label{eqn:dbar-eq-estimate_bis}||u||^2=\langle\langle\Delta^{''-1}_{k,\,(\eta)}v,\, v\rangle\rangle\leq\langle\langle B^{-1}_{k,\,(\eta)}v,\, v\rangle\rangle\leq\frac{1}{c}\,||v||^2.\end{equation}

\end{Obs}

\noindent {\it Proof.} $(i)$\, Since $X$ is compact and $\Delta''_{k,\,(\eta)}$ is elliptic (see (\ref{eqn:twisted-lap''-formula2}) below), G$\mathring{a}$rding's inequality implies the two-space orthogonal decomposition ${\cal C}^\infty_{p,\,q}(X,\,L_k) = \ker\Delta''_{k,\,(\eta)}\oplus\mbox{Im}\,\Delta''_{k,\,(\eta)}$, while the assumption $\bar\partial_k^2=0$ implies $^{\eta}\bar\partial_k\circ\bar\partial_k^{\eta}=0$ which, in turn, gives the further orthogonal decomposition $\mbox{Im}\,\Delta''_{k,\,(\eta)} = \mbox{Im}\,\bar\partial_k^{\eta}\oplus\mbox{Im}\,(^{\eta}\bar\partial_k)^{\star}.$  

$(ii)$\, The solution of equation (\ref{eqn:dbar-eq-integrable}) is unique modulo $\ker\bar\partial_k^{\eta}$, so the minimal $L^2$ norm solution is the unique solution lying in $(\ker\bar\partial_k^{\eta})^{\perp}=\mbox{Im}\,(\bar\partial_k^{\eta})^{\star}$. Since the form $u$ defined in (\ref{eqn:dbar-eq-integrable-sol}) is obviously $(\bar\partial_k^{\eta})^{\star}$-exact, it remains to verify that $\bar\partial_k^{\eta}(\bar\partial_k^{\eta})^{\star}\Delta^{''-1}_{k,\,(\eta)}v = v$, which is straightforward.    \hfill $\Box$

\subsection{Bochner-Kodaira-Nakano formula for the twisted Laplacians}\label{subsection:twistedBKN}

 We now calculate the difference $\Delta''_{k,\,(\eta)}-\Delta'_{k,\,(\eta)}$.

Formula (\ref{eqn:twisted-lap''-formula}) gives: \begin{equation}\label{eqn:twisted-lap''-formula1}\Delta''_{k,\,(\eta)} = \eta\,\bar\partial_k\bar\partial_k^{\star} + \bar\partial\eta\wedge\bar\partial_k^{\star} + \bar\partial_k^{\star}\eta\bar\partial_k.\end{equation}

\noindent Since the metric $\omega$ of $X$ is supposed to be K\"ahler, $\bar\partial_k^{\star} = -i\,[\Lambda,\,\partial_k]$ (see Proposition \ref{Prop:Hermitian_commutation-rel}), so we get: \begin{eqnarray}\nonumber\bar\partial_k^{\star}\eta\bar\partial_k & = & -i\Lambda\partial_k\eta\bar\partial_k + i\partial_k\Lambda\eta\bar\partial_k\\
 \nonumber & = & -i\Lambda(\eta\partial_k\bar\partial_k + \partial\eta\wedge\bar\partial_k) + \eta\,i\partial_k\Lambda\bar\partial_k + \partial\eta\wedge\,i\Lambda\bar\partial_k\\
\nonumber & = & -\eta\,i\,[\Lambda,\,\partial_k]\,\bar\partial_k -i\,[\Lambda,\,\partial\eta\wedge\cdot]\,\bar\partial_k.\end{eqnarray}

Since $i\,[\Lambda,\,\partial_k] = -\bar\partial_k^{\star}$ (see above) and $i\,[\Lambda,\,\partial\eta\wedge\cdot] = (\bar\partial\eta\wedge\cdot)^{\star}$ (well known, see e.g. [Dem01, p.97] or (a) of Lemma \ref{Lem:com}), we get: \begin{equation}\label{eqn:twisted-lap''-formula2}\Delta''_{k,\,(\eta)} = \eta\,\Delta''_k + \bar\partial\eta\wedge\bar\partial_k^{\star} - (\bar\partial\eta\wedge\cdot)^{\star}\,\bar\partial_k.\end{equation} This shows, in particular, that the twisted Laplacian $\Delta''_{k,\,(\eta)}$ is elliptic since $\Delta''_k$ is, hence so is $\eta\Delta''_k$, while $\Delta''_{k,\,(\eta)}$ and $\eta\Delta''_k$ have the same principal part by (\ref{eqn:twisted-lap''-formula2}).

 In a similar way, we get from formula (\ref{eqn:twisted-lap'-formula}) the identity

\begin{equation}\label{eqn:twisted-lap'-formula2}\Delta'_{k,\,(\eta)} = \eta\,\Delta'_k + \partial\eta\wedge\partial_k^{\star} - (\partial\eta\wedge\cdot)^{\star}\,\partial_k\end{equation}

\noindent which shows, in particular, that $\Delta'_{k,\,(\eta)}$ is elliptic. Using the non-integrable version of the standard Bochner-Kodaira-Nakano identity in the K\"ahler case (see [Lae02]): $$\Delta''_k = \Delta'_k + [i\Theta(D_k)^{1,\,1},\,\Lambda] = \Delta'_k + 2\pi\,[\alpha_k^{1,\,1},\,\Lambda],$$

\noindent a combination of (\ref{eqn:twisted-lap''-formula2}) and (\ref{eqn:twisted-lap'-formula2}) yields: \begin{eqnarray}\label{eqn:pre-twistedBKN}\nonumber\Delta''_{k,\,(\eta)} - \Delta'_{k,\,(\eta)} = \eta\,[i\Theta(D_k)^{1,\,1},\,\Lambda] & + & \bar\partial\eta\wedge\bar\partial_k^{\star} - (\bar\partial\eta\wedge\cdot)^{\star}\,\bar\partial_k\\
 & - & \partial\eta\wedge\partial_k^{\star} + (\partial\eta\wedge\cdot)^{\star}\,\partial_k.\end{eqnarray}

\noindent To transform the term $(\bar\partial\eta\wedge\cdot)^{\star}\,\bar\partial_k = i\,[\Lambda,\,\partial\eta\wedge\cdot]\,\bar\partial_k$, we use the Jacobi identity

$$[[\Lambda,\,\partial\eta\wedge\cdot],\, \bar\partial_k] + [[\partial\eta\wedge\cdot,\,\bar\partial_k],\,\Lambda] - [[\bar\partial_k,\,\Lambda],\,\partial\eta\wedge\cdot]=0.$$

\noindent Since $[\partial\eta\wedge\cdot,\,\bar\partial_k] = -\partial\bar\partial\eta\wedge\cdot$ (obvious) and $[\bar\partial_k,\,\Lambda] = i\partial_k^{\star}$ (by K\"ahler commutation), the Jacobi identity translates to

\begin{equation}\label{eqn:jacobi-transf}(\bar\partial\eta\wedge\cdot)^{\star}\bar\partial_k + \bar\partial_k(\bar\partial\eta\wedge\cdot)^{\star}-[i\partial\bar\partial\eta\wedge\cdot,\,\Lambda] + \partial_k^{\star}\partial\eta\wedge\cdot + \partial\eta\wedge\partial_k^{\star}=0.\end{equation}

\noindent Combining (\ref{eqn:pre-twistedBKN}) and (\ref{eqn:jacobi-transf}), we get: \begin{eqnarray}\nonumber\Delta''_{k,\,(\eta)} - \Delta'_{k,\,(\eta)} = \eta\,[i\Theta(D_k)^{1,\,1},\,\Lambda] & + & \bar\partial\eta\wedge\bar\partial_k^{\star} + \bar\partial_k(\bar\partial\eta\wedge\cdot)^{\star}\\
\nonumber  & + & (\partial\eta\wedge\cdot)^{\star}\partial_k + \partial_k^{\star}(\partial\eta\wedge\cdot)\\
\nonumber & - & [i\partial\bar\partial\eta\wedge\cdot,\,\Lambda]  \hspace{6ex} \mbox{on}\hspace{1ex} C^\infty_{p,\,q}(X,\,L_k).\end{eqnarray}

\noindent When $p=n$, the operators $\partial_k$ and $\partial\eta\wedge\cdot$ vanish for bidegree reasons, so we get the following.

\begin{Prop}\label{Prop:twistedBKN} For any ${\cal C}^\infty$ function $\eta:X\longrightarrow (0,\, +\infty)$, the following {\bf twisted Bochner-Kodaira-Nakano identity} holds

\begin{equation}\label{eqn:twistedBKN}\Delta''_{k,\,(\eta)} = \Delta'_{k,\,(\eta)} + [\eta\,i\Theta(D_k)^{1,\,1} - i\partial\bar\partial\eta\wedge\cdot,\,\Lambda] + \bar\partial\eta\wedge\bar\partial_k^{\star} + \bar\partial_k(\bar\partial\eta\wedge\cdot)^{\star}\end{equation}

\noindent on ${\cal C}^\infty_{n,_,q}(X,\,L_k)$.

\end{Prop}

 A second positive function $\lambda >0$ is now introduced to estimate the loose terms above. It follows from (\ref{eqn:twistedBKN}), (\ref{eqn:twisted-lap''-scal}) and (\ref{eqn:twisted-lap'-scal}) that for every $v\in {\cal D}^{n,_,q}(X,\,L_k)$ and for every ${\cal C}^\infty$ functions $\eta, \lambda:X\longrightarrow (0,\, +\infty)$, we have: \begin{eqnarray}\label{eqn:BKNineq1}\nonumber  & & ||\sqrt{\eta}\,\bar\partial_kv||^2 + ||\sqrt{\eta}\,\bar\partial_k^{\star}v||^2 - ||\sqrt{\eta}\,\partial_kv||^2 - ||\sqrt{\eta}\,\partial_k^{\star}v||^2\\
\nonumber & = & \langle\langle[\eta\,i\Theta(D_k)^{1,\,1} - i\partial\bar\partial\eta\wedge\cdot,\,\Lambda]v,\,v\rangle\rangle + 2\mbox{Re}\,\langle\langle\bar\partial_k^{\star}v,\,(\bar\partial\eta\wedge\cdot)^{\star}v\rangle\rangle\\
\nonumber & \geq & \langle\langle[\eta\,i\Theta(D_k)^{1,\,1} - i\partial\bar\partial\eta\wedge\cdot,\,\Lambda]v,\,v\rangle\rangle - \langle\langle\lambda\,\bar\partial_k^{\star}v,\,\bar\partial_k^{\star}v\rangle\rangle - \langle\langle\frac{1}{\lambda}\,(\bar\partial\eta\wedge\cdot)(\bar\partial\eta\wedge\cdot)^{\star}v,\,v\rangle\rangle,\end{eqnarray}

\noindent where the Cauchy-Schwarz inequality has been used. Neglecting the non-positive terms $-||\sqrt{\eta}\,\partial_kv||^2$ and $-||\sqrt{\eta}\,\partial_k^{\star}v||^2$, we get for all $v\in {\cal D}^{n,\,q}(X,\, L_k)$: \begin{eqnarray}\label{eqn:BKNineq2}\nonumber\langle\langle(\eta + \lambda)\,\bar\partial_k^{\star}v,\,\bar\partial_k^{\star}v\rangle\rangle + \langle\langle\eta\,\bar\partial_kv,\,\bar\partial_kv\rangle\rangle & \geq & \langle\langle[\eta\,i\Theta(D_k)^{1,\,1} - i\partial\bar\partial\eta\wedge\cdot,\,\Lambda]v,\,v\rangle\rangle\\
  & - & \langle\langle\frac{1}{\lambda}\,(\bar\partial\eta\wedge\cdot)(\bar\partial\eta\wedge\cdot)^{\star}v,\,v\rangle\rangle.\end{eqnarray}

\begin{Lem} For any real-valued smooth function $\eta$, the following identity holds on $(n,\,q)$-forms at every point of $X$ \begin{equation}\label{eqn:dbar-dbar-star-eta}(\bar\partial\eta\wedge\cdot)(\bar\partial\eta\wedge\cdot)^{\star} = [i\partial\eta\wedge\bar\partial\eta\wedge\cdot,\,\Lambda].\end{equation}

\end{Lem}

\noindent {\it Proof.} Part $(b)$ of Lemma \ref{Lem:com} in the Appendix yields in every bidegree: $$(\partial\eta\wedge\cdot)^{\star}(\partial\eta\wedge\cdot) - (\bar\partial\eta\wedge\cdot)(\bar\partial\eta\wedge\cdot)^{\star} = [i\bar\partial\eta\wedge\partial\eta\wedge\cdot,\,\Lambda].$$

\noindent This gives (\ref{eqn:dbar-dbar-star-eta}) when applied to $v\in {\cal C}^\infty_{n,\,q}(X,\,L_k)$ since $\partial\eta\wedge v=0$ for bidegree reasons. \hfill $\Box$

\vspace{3ex}

Combining (\ref{eqn:BKNineq2}) and (\ref{eqn:dbar-dbar-star-eta}) (and dropping $\wedge\cdot$), we get for any $v\in {\cal D}^{n,\,q}(X,\,L_k)$: $$\langle\langle(\eta + \lambda)\,\bar\partial_k^{\star}v,\,\bar\partial_k^{\star}v\rangle\rangle + \langle\langle\eta\,\bar\partial_kv,\,\bar\partial_kv\rangle\rangle \geq \langle\langle[\eta\,i\Theta(D_k)^{1,\,1} - i\partial\bar\partial\eta - \frac{1}{\lambda}\,i\partial\eta\wedge\bar\partial\eta,\,\Lambda]v,\,v\rangle\rangle.$$

\noindent Since $\lambda>0$, the left-hand side above is bounded above by 

$$\langle\langle(\eta + \lambda)\,\bar\partial_k^{\star}v,\,\bar\partial_k^{\star}v\rangle\rangle + \langle\langle(\eta + \lambda)\,\bar\partial_kv,\,\bar\partial_kv\rangle\rangle = \langle\langle\Delta''_{k,\,(\eta+\lambda)}v,\, v\rangle\rangle \hspace{2ex} \mbox{(cf.}\hspace{1ex} (\ref{eqn:twisted-lap''-scal})),$$

\noindent so we have obtained the following non-integrable analogue of the {\it a priori estimate} of [Ohs95].

\begin{Prop}\label{Prop:twistedBKNineq} For any ${\cal C}^\infty$ functions $\eta, \lambda:X\longrightarrow (0,\, +\infty)$, the following {\bf twisted Bochner-Kodaira-Nakano inequality} holds: \begin{equation}\label{eqn:twistedBKNineq}\langle\langle\Delta''_{k,\,(\eta+\lambda)}v,\, v\rangle\rangle\geq\langle\langle[\eta\,i\Theta(D_k)^{1,\,1} - i\partial\bar\partial\eta - \frac{1}{\lambda}\,i\partial\eta\wedge\bar\partial\eta,\,\Lambda]v,\,v\rangle\rangle\end{equation} \noindent for all $v\in {\cal D}^{n,\,q}(X,\,L_k)$. (Recall that $i\Theta(D_k)^{1,\,1} = 2\pi\alpha_k^{1,\,1 }$ is close to $2\pi k\alpha$ in the ${\cal C}^\infty$ topology.)

\end{Prop}

 In other words, the following inequality between linear operators holds on ${\cal D}^{n,\,q}(X,\,L_k)$, hence also on $\mbox{Dom}\,\Delta''_{k,\,(\eta+\lambda)}\subset L^2_{n,\,q}(X,\,L_k)$ by density of ${\cal D}^{n,\,q}(X,\,L_k)$ in $\mbox{Dom}\,\Delta''_{k,\,(\eta+\lambda)}$ w.r.t. the graph norm (which is a consequence of the metric $\omega$ being {\bf complete}): \begin{equation}\label{eqn:twistedBKNineq_condensed}\Delta''_{k,\,(\eta+\lambda)}\geq B=B^{(k)}_{\eta,\,\lambda}:=[\eta\,i\Theta(D_k)^{1,\,1} - i\partial\bar\partial\eta - \frac{1}{\lambda}\,i\partial\eta\wedge\bar\partial\eta,\,\Lambda].\end{equation}

\noindent The zeroth-order operator $B=B^{(k)}_{\eta,\,\lambda}:L^2_{n,\,q}(X,\,L_k)\longrightarrow L^2_{n,\,q}(X,\,L_k)$ will be called the {\it twisted curvature operator} in bidegree $(n,\,q)$.

Note that the lower bound (\ref{eqn:twistedBKNineq_condensed}) on the twisted Laplacian $\Delta''_{k,\,(\eta+\lambda)}$ in terms of the twisted curvature operator $B$ does not require any positivity assumption on the curvature form $i\Theta(D_k)^{1,\,1}$. The second stage in our quest for a positive lower bound on $\Delta''_{k,\,(\eta+\lambda)}$, performed in the next subsection, will be to get a positive lower bound on $B$ under a suitable positivity assumption on the curvature.

\subsection{Application to the $L^2$ extension context}\label{subsection:OTapplication}

As before, $X$ is supposed to be a complex manifold of dimension $n$ endowed with a {\it complete} K\"ahler metric $\omega$ such that $b_2$ is finite. Suppose that $w:X\longrightarrow\C$ is a bounded holomorphic function such that

$$Y:=\{x\in X\,;\,w(x)=0\}\subset X$$

\noindent is nonsingular. In other words, we assume that $dw(x)\neq 0$ for all $x\in Y$. Assume for convenience that $\sup_X|w|\leq 1$. We wish to extend a certain class of approximately holomorphic sections from $Y$ to $X$ in such a way that the $L^2$-norm on $X$ of the extension is controlled in terms of the $L^2$-norm on $Y$ of the original section.

 The functions $\eta, \lambda$ of the previous subsections will be chosen to suit $Y$. They are the same as in [OT87] and [Siu02, $\S.3$]. Specifically, fix any real number $A>e$ and, for every $0<\varepsilon<\delta_0:=\sqrt{\frac{A}{e}-1}$, consider the bounded ${\cal C}^\infty$ functions on $X$: \begin{eqnarray}\label{eqn:eta-lambda-eps-def}\nonumber 1<\eta_{\varepsilon} &:= &\log\frac{A}{|w|^2 + \varepsilon^2}\leq\log\frac{A}{\varepsilon^2},\\
  0< \lambda_{\epsilon} & := & \frac{1}{|w|^2 + \varepsilon^2}\leq\frac{1}{\varepsilon^2}.\end{eqnarray}

\noindent Thus, for every $x\in Y$, $\eta_{\varepsilon}(x)=\log\frac{A}{\varepsilon^2}\nearrow +\infty$ as $\varepsilon\downarrow 0$, creating a ``bump'' along $Y$. We can choose $w$ as one of the $n$ coordinates along $Y$, i.e. for every $x\in Y$, we can find local holomorphic coordinates $(z_1,\dots , z_{n-1}, w)$ on $X$ centred at $x$.  

 We now calculate how (\ref{eqn:twistedBKNineq}) transforms with this choice of functions $\eta_{\varepsilon}, \lambda_{\varepsilon}$. We have: \begin{eqnarray}\label{eqn:iddbar-eta-epsilon} i\partial\bar\partial\eta_{\varepsilon} = -\frac{\partial^2\log(|w|^2 + \varepsilon^2)}{\partial w\partial\overline{w}}\,idw\wedge d\overline{w} = -\frac{\varepsilon^2}{(|w|^2 + \varepsilon^2)^2}\,idw\wedge d\overline{w}.\end{eqnarray} On the other hand, \begin{eqnarray}\label{eqn:idel-eta-epsilon}\nonumber i\partial\eta_{\varepsilon} & = & -i\partial_w\log(|w|^2 + \varepsilon^2) = -\frac{\overline{w}}{|w|^2 + \varepsilon^2}\,idw,\\
\nonumber i\bar\partial\eta_{\varepsilon}  & = & -i\bar\partial_w\log(|w|^2 + \varepsilon^2) = -\frac{w}{|w|^2 + \varepsilon^2}\,id\overline{w},\end{eqnarray}

\noindent hence $\frac{1}{\lambda_{\varepsilon}}\,i\partial\eta_{\varepsilon}\wedge\bar\partial\eta_{\varepsilon} = \frac{|w|^2}{|w|^2 + \varepsilon^2}\,idw\wedge d\overline{w}$. Combining this with (\ref{eqn:iddbar-eta-epsilon}), we get

$$i\partial\bar\partial\eta_{\varepsilon} = -\frac{\varepsilon^2}{|w|^2\,(|w|^2 + \varepsilon^2)}\,\frac{1}{\lambda_{\varepsilon}}\,i\partial\eta_{\varepsilon}\wedge\bar\partial\eta_{\varepsilon}.$$

\noindent Hence, we finally get \begin{eqnarray}\nonumber -i\partial\bar\partial\eta_{\varepsilon} - \frac{1}{\lambda_{\varepsilon}}\,i\partial\eta_{\varepsilon}\wedge\bar\partial\eta_{\varepsilon} & = & \bigg(\frac{\varepsilon^2}{|w|^2\,(|w|^2 + \varepsilon^2)} - 1\bigg)\,\frac{|w|^2}{|w|^2 + \varepsilon^2}\,idw\wedge d\overline{w} \\
\nonumber & = & \frac{\varepsilon^2 - |w|^2\,(|w|^2 + \varepsilon^2)}{(|w|^2 + \varepsilon^2)^2}\,idw\wedge d\overline{w}\\
\nonumber  & \geq & \bigg(\frac{\varepsilon^2}{(|w|^2 + \varepsilon^2)^2}-1\bigg)\,idw\wedge d\overline{w} = \frac{\varepsilon^2}{(|w|^2 + \varepsilon^2)^2}\,idw\wedge d\overline{w} - idw\wedge d\overline{w},\end{eqnarray}

\noindent where the last inequality follows from $\frac{|w|^2\,(|w|^2 + \varepsilon^2)}{(|w|^2 + \varepsilon^2)^2}\leq 1$. (Note also that $idw\wedge d\overline{w} = i\partial\bar\partial|w|^2$.) Since $\eta_{\varepsilon}>1$ on $X$, we infer that \begin{eqnarray}\label{eqn:idel-eta-epsilon2}\eta_{\varepsilon}\,idw\wedge d\overline{w}-i\partial\bar\partial\eta_{\varepsilon} - \frac{1}{\lambda_{\varepsilon}}\,i\partial\eta_{\varepsilon}\wedge\bar\partial\eta_{\varepsilon}\geq \frac{\varepsilon^2}{(|w|^2 + \varepsilon^2)^2}\,idw\wedge d\overline{w}.\end{eqnarray}

After these preliminaries, we get round to deriving a positive lower bound for the twisted curvature operator $B_\varepsilon=B^{(k)}_{\eta_{\varepsilon},\,\lambda_{\varepsilon}}$ (see (\ref{eqn:twistedBKNineq_condensed})) induced by the auxiliary functions $\eta=\eta_{\varepsilon}$ and $\lambda=\lambda_{\varepsilon}$ under a positivity assumption on the curvature form $i\Theta(D_k)^{1,\,1}$.

\begin{Prop}\label{Prop:curvature-lbound} Suppose $\alpha$ is a ${\cal C}^\infty$ real $(1,\,1)$-form on $X$ such that $d\alpha=0$, $\alpha$ satisfies Assumption A and a choice $(L_k,\,h_k,\,D_k)_{k\in\Sigma}$ of approximating sequence (cf. Definition \ref{Def:approx-seq}) has been fixed for $\alpha$.

If $i\Theta(D_k)^{1,\,1}\geq k\omega + idw\wedge d\overline{w}$ on $X$ for some $k\geq 1$, then for every $0<\varepsilon<\delta_0=\sqrt{(A/e)-1}$ we have: \begin{equation}\label{eqn:curvature-lbound}\nonumber\eta_{\varepsilon}\,i\Theta(D_k)^{1,\,1} - i\partial\bar\partial\eta_{\varepsilon} - \frac{1}{\lambda_{\varepsilon}}\,i\partial\eta_{\varepsilon}\wedge\bar\partial\eta_{\varepsilon} \geq k\omega + \frac{\varepsilon^2}{(|w|^2 + \varepsilon^2)^2}\,idw\wedge d\overline{w}.\end{equation}

\end{Prop}

\noindent {\it Proof.} Since $\eta_{\varepsilon}>1>0$ on $X$, the stated inequality follows from (\ref{eqn:idel-eta-epsilon2}).     \hfill $\Box$

\vspace{2ex}

The following fact is well known. We recall its proof for the reader's convenience.

\begin{Lem}\label{Lem:alpha-Lambda} If $\alpha$ is a ${\cal C}^\infty$ $(1,\,1)$-form on $X$, the following implication holds: $$\alpha\geq 0 \hspace{2ex} \Longrightarrow \hspace{2ex}[\alpha,\,\Lambda]\geq 0 \hspace{2ex} \mbox{pointwise, as an operator on}\hspace{1ex} {\cal C}^\infty_{n,\,q}(X,\,\C)\hspace{1ex}\mbox{or}\hspace{1ex}{\cal C}^\infty_{n,\,q}(X,\,L_k).$$

\end{Lem}

\noindent {\it Proof.} Fix any point $x\in X$ and choose local holomorphic coordinates $z_1, \dots , z_n$ about $x$ such that $$\omega(x) = i\sum\limits_{j=1}^ndz_j\wedge d\bar{z_j}  \hspace{2ex} \mbox{and} \hspace{2ex} \alpha(x) = i\sum\limits_{j=1}^n\alpha_j\,dz_j\wedge d\bar{z_j},$$ \noindent where $0\leq\alpha_1\leq\dots\leq\alpha_n$ are the eigenvalues of $\alpha$ w.r.t. $\omega$. Then, for every $v=\sum\limits_{|K|=q}v_K\,dz_1\wedge\dots\wedge dz_n\wedge d\bar{z}_K\in {\cal C}^\infty_{n,\,q}(X,\,\C)$, the following well-known formula (cf. e.g. [Dem01]) holds for the pointwise inner product induced by $\omega$:

$$\langle[\alpha,\,\Lambda]v,\,v\rangle = \sum\limits_{|K|=q}\bigg((\alpha_1 + \dots + \alpha_n) + \sum\limits_{j\in K}\alpha_j - (\alpha_1 + \dots + \alpha_n)\bigg)\,|v_K|^2\geq 0.$$  \hfill $\Box$

\vspace{3ex}

 We are now in a position to derive the desired positive lower bound for the twisted Laplacian and the twisted curvature operator. The discussion of this section is summed up in the following

 \begin{Prop}\label{Prop:BKN-consequence} Let $X$ be a complex manifold of dimension $n$ endowed with a {\bf complete} K\"ahler metric $\omega$ such that $b_2$ is {\bf finite}. Suppose $\alpha$ is a ${\cal C}^\infty$ real $(1,\,1)$-form on $X$ such that $d\alpha=0$ and $\alpha$ satisfies Assumption A for a choice $(L_k,\,h_k,\,D_k)_{k\in\Sigma}$ of approximating sequence (cf. Definition \ref{Def:approx-seq}). 

   Suppose, moreover, that $i\Theta(D_k)^{1,\,1}\geq k\omega + idw\wedge d\overline{w}$ on $X$ for some $k\in\Sigma$ with $k\geq 2$.

   \vspace{1ex}

   Then, for all $0<\varepsilon<\delta_0=\sqrt{(A/e)-1}$ and all $q$,  

\vspace{1ex}

$(i)$\, the following inequalities hold: \begin{eqnarray}\label{eqn:l-bound_twisted-Laplacian}\nonumber\Delta''_{k,\,(\eta_{\varepsilon}+\lambda_{\varepsilon})} & \geq & B_\varepsilon:=[\eta_\varepsilon\,i\Theta(D_k)^{1,\,1} - i\partial\bar\partial\eta_\varepsilon - \frac{1}{\lambda_\varepsilon}\,i\partial\eta_\varepsilon\wedge\bar\partial\eta_\varepsilon,\,\Lambda] \\
\nonumber & \geq & kq\,\mbox{Id} +\frac{\varepsilon^2}{(|w|^2 + \varepsilon^2)^2}\,(d\overline{w}\wedge\cdot)(d\overline{w}\wedge\cdot)^{\star}\end{eqnarray}

\noindent on ${\cal D}^{n,\,q}(X,\,L_k)$. The first inequality holds in the sense of linear operators on $\mbox{Dom}\,\Delta''_{k,\,(\eta_{\varepsilon}+\lambda_{\varepsilon})}\subset L^2_{n,\,q}(X,\,L_k)$. The latter inequality holds pointwise on $\Lambda^{n,\,q}T^\star X\otimes L_k$.

\vspace{1ex}

$(ii)$\, the following inequality holds {\bf pointwise} on $X$ for every $u\in {\cal C}^\infty_{n,\,q}(X,\,L_k)$:

\begin{equation}\label{eqn:B-inverse_upper-bound}\nonumber\langle B_\varepsilon^{-1}(d\overline{w}\wedge u),\,d\overline{w}\wedge u\rangle\leq\frac{|d\overline{w}|_\omega^2}{kq +\frac{\varepsilon^2}{(|w|^2 + \varepsilon^2)^2}\,|d\overline{w}|_\omega^2}\,|u|_\omega^2.\end{equation} 

\noindent Hence, for every $u\in{\cal D}^{n,\,q}(X,\,L_k)$, \begin{eqnarray}\label{eqn:twisted-Laplace-inverse_upper-bound}\nonumber\langle\langle\Delta_{k,\,(\eta_{\varepsilon}+\lambda_{\varepsilon})}^{''-1}(d\overline{w}\wedge u),\,d\overline{w}\wedge u\rangle\rangle  & \leq & \langle\langle B_\varepsilon^{-1}(d\overline{w}\wedge u),\,d\overline{w}\wedge u\rangle\rangle \leq \int\limits_X\frac{|d\overline{w}|_\omega^2}{kq +\frac{\varepsilon^2}{(|w|^2 + \varepsilon^2)^2}\,|d\overline{w}|_\omega^2}\,|u|_\omega^2\,dV_\omega.\end{eqnarray}

\vspace{1ex}

$(iii)$\, the following inequality holds {\bf pointwise} on $X$ for every $v\in {\cal C}^\infty_{n,\,q}(X,\,L_k)$:

$$\langle B_\varepsilon^{-1}v,\,v\rangle\leq\frac{1}{kq}\,|v|^2_\omega.$$

\noindent Hence, for every $v\in{\cal D}^{n,\,q}(X,\,L_k)$ (not necessarily divisible by $d\overline{w}$), we get the $L^2$ estimates: \begin{eqnarray}\label{eqn:twisted-Laplace-inverse_upper-bound}\nonumber\langle\langle\Delta_{k,\,(\eta_{\varepsilon}+\lambda_{\varepsilon})}^{''-1}v,\,v\rangle\rangle  & \leq & \langle\langle B_\varepsilon^{-1}v,\,v\rangle\rangle\leq\frac{1}{kq}\,||v||^2_\omega.\end{eqnarray}

\end{Prop}

\noindent {\it Proof.} $(i)$\, The stated inequalities follow by putting together Propositions \ref{Prop:twistedBKNineq}, \ref{Prop:curvature-lbound} and Lemma \ref{Lem:alpha-Lambda}. For the latter inequality, we also use the identity $|w|^2\,idw\wedge d\overline{w} = i\partial|w|^2\wedge\bar\partial|w|^2$ and (\ref{eqn:dbar-dbar-star-eta}) that between them yield $[idw\wedge d\overline{w}\wedge\cdot,\,\Lambda] = (d\overline{w}\wedge\cdot)(d\overline{w}\wedge\cdot)^{\star}$ pointwise on $(n,\,q)$-forms. The standard identity $[L_\omega,\,\Lambda_\omega] = (p+q-n)\,\mbox{Id}$ on $(p,\,q)$-forms (hence yielding $q\,\mbox{Id}$ on $(n,\,q)$-forms) was also used, where $L_\omega=\omega\wedge\cdot$.

\vspace{1ex}

$(ii)$\, By the latter inequality in $(i)$, for every $(n,\,q)$-form $u$ we have

\vspace{1ex}

$\langle B_\varepsilon(d\overline{w}\wedge u),\,d\overline{w}\wedge u\rangle \geq
kq\,|d\overline{w}\wedge u|^2 + \frac{\varepsilon^2}{(|w|^2 + \varepsilon^2)^2}\,\langle(d\overline{w}\wedge\cdot)(d\overline{w}\wedge\cdot)^{\star}(d\overline{w}\wedge u),\,d\overline{w}\wedge u\rangle.$

\vspace{1ex}

\noindent Now, from Lemma \ref{Lem:d-bar_eta_star_no-star} in the Appendix and the fact that $d\overline{w} = \bar\partial\overline{w}$ (since $w$ is holomorphic), we get

\vspace{1ex}

\hspace{6ex} $(d\overline{w}\wedge\cdot)^{\star}(d\overline{w}\wedge u) = |d\overline{w}|^2_\omega\,u - d\overline{w}\wedge(d\overline{w}\wedge\cdot)^{\star}u.$

\vspace{1ex}

\noindent However, the term $d\overline{w}\wedge(d\overline{w}\wedge\cdot)^{\star}u$ vanishes after it gets multiplied by $d\overline{w}$ since $d\overline{w}\wedge d\overline{w}=0$. Hence

$$\langle B_\varepsilon(d\overline{w}\wedge u),\,d\overline{w}\wedge u\rangle \geq \bigg[kq +\frac{\varepsilon^2}{(|w|^2 + \varepsilon^2)^2}\,|d\overline{w}|^2_\omega\bigg]\,|d\overline{w}\wedge u|^2_\omega$$

\noindent for all $(n,\,q)$-forms $u$. This means that

\vspace{1ex}

$\displaystyle B_\varepsilon\geq\bigg[kq +\frac{\varepsilon^2}{(|w|^2 + \varepsilon^2)^2}\,|d\overline{w}|^2_\omega\bigg]\,\mbox{Id}, \hspace{3ex} \mbox{hence} \hspace{3ex} B_\varepsilon^{-1}\leq \frac{1}{kq +\frac{\varepsilon^2}{(|w|^2 + \varepsilon^2)^2}\,|d\overline{w}|^2_\omega}\,\mbox{Id}$

\vspace{1ex}

\noindent pointwise on $(n,\,q+1)$-forms that are divisible by $d\overline{w}$, where $\mbox{Id}$ is the identity operator. Consequently,

$$\langle B_\varepsilon^{-1}(d\overline{w}\wedge u),\,d\overline{w}\wedge u\rangle \leq \frac{1}{kq +\frac{\varepsilon^2}{(|w|^2 + \varepsilon^2)^2}\,|d\overline{w}|^2_\omega}\,|d\overline{w}\wedge u|^2_\omega.$$

\noindent Since $|d\overline{w}\wedge u|^2_\omega\leq |d\overline{w}|^2_\omega \,|u|^2_\omega$, the stated pointwise inequality follows.

 The $L^2$ inequality (\ref{eqn:twisted-Laplace-inverse_upper-bound}) follows by integrating the pointwise inequality.  

\vspace{1ex}

$(iii)$\, Both inequalities are immediate consequences of $(i)$ after ignoring the non-negative term $[\varepsilon^2/(|w|^2 + \varepsilon^2)^2]\,(d\overline{w}\wedge\cdot)(d\overline{w}\wedge\cdot)^{\star}$.     \hfill $\Box$

\section{Approximately holomorphic sections}\label{section:spaces-approx-def}

 The setting and the notation are the same as in the previous sections.

\subsection{Definitions and basic properties}\label{subsection:definitions_approx-hol} Fix $k\in\Sigma\subset\N^{\star}$ and a bidegree $(p,\,q)$. Besides $\Delta''_k:\mbox{Dom}\,\Delta''_k\longrightarrow L^2_{p,\,q}(X,\,L_k)$, the following discussion also applies to $\Delta'_k$ and to various twisted versions of these operators.

Self-adjointness ensures that the spectral theorem holds for $\Delta''_k$. In particular, for every Borel set $\Omega\subset\R$, we can define (cf. [RS80, chapter VIII]) the {\it spectral projection} of $\Delta''_k$ (a bounded operator): $$\chi_{\Omega}(\Delta''_k)\,:\,L^2_{p,\,q}(X,\,L_k)\longrightarrow L^2_{p,\,q}(X,\,L_k),$$

\noindent where $\chi_{\Omega}$ is the {\it characteristic function} of $\Omega$. If $X$ is compact, $\chi_{\Omega}(\Delta''_k)$ is the orthogonal projection onto the direct sum of the eigenspaces of $\Delta''_k$ corresponding to the eigenvalues $\lambda\in\Omega$.

The family of spectral projections $\{\chi_{\Omega}(\Delta''_k)\}$ enjoys the following properties (cf. [RS80]):

\vspace{1ex}

$(i)$\, $\chi_{\Omega}(\Delta''_k) = \chi_{\Omega}(\Delta''_k)^2 = \chi_{\Omega}(\Delta''_k)^{\star}$ for any Borel set $\Omega\subset\R$;

\vspace{1ex}

$(ii)$\, $\chi_{\emptyset}(\Delta''_k) = 0$ and $\chi_{(-\infty,\,+\infty)}(\Delta''_k) = \mbox{Id}_{L^2_{p,\,q}(X,\,L_k)}$;

\vspace{1ex}

$(iii)$\, if $\Omega=\bigcup\limits_{l=1}^{+\infty}\Omega_l$ with $\Omega_l\cap\Omega_r=\emptyset$ for all $l\neq r$, then

$$\chi_{\Omega}(\Delta''_k) = \lim\limits_{N\rightarrow +\infty}\sum\limits_{l=1}^{N}\chi_{\Omega_l}(\Delta''_k),$$

\noindent in the strong operator topology limit;

\vspace{1ex}

$(iv)$\, $\chi_{\Omega_1}(\Delta''_k)\,\chi_{\Omega_2}(\Delta''_k) = \chi_{\Omega_1\cap\Omega_2}(\Delta''_k)$ for all Borel sets $\Omega_1,\Omega_2\subset\R$.

\vspace{2ex}

Whether $P$ is $\Delta''_k$ or any of the other Laplacians defined above, we introduce the following 

\begin{Def}\label{Def:H-Q_spaces_arbitrary} Fix a bidegree $(p,\,q)$, an integer $k\in\Sigma\subset\N^{\star}$ and an interval  $I\subset\R$. Let $P:\mbox{Dom}\,P\longrightarrow L^2_{p,\,q}(X,\,L_k)$ be an elliptic self-adjoint linear differential operator such that $P\geq 0$, where $\mbox{Dom}\,P\subset L^2_{p,\,q}(X,\,L_k)$ consists of the forms $u\in L^2_{p,\,q}(X,\,L_k)$ for which $Pu$ computed in the sense of distributions lies in $L^2_{p,\,q}(X,\,L_k)$.

\vspace{1ex}  

 We set: \begin{equation*}\label{eqn:H_spaces_arbitrary}{\cal H}_{I,\,P}^{p,\,q}(X,\,L_k):= \mbox{Im}\,\bigg(\chi_I(P):L^2_{p,\,q}(X,\,L_k)\rightarrow L^2_{p,\,q}(X,\,L_k)\bigg)\subset L^2_{p,\,q}(X,\,L_k)\end{equation*} and, if $I=[a,\,b]$ with $-\infty\leq a\leq b\leq +\infty$, \begin{equation*}\label{eqn:H_spaces_arbitrary}Q_{I,\,P}^{p,\,q}(X,\,L_k):=\bigg\{u\in\mbox{Dom}\,P\,\bigg|\, a\,||u||^2\leq\langle\langle Pu,\,u\rangle\rangle\leq b\,||u||^2\bigg\}.\end{equation*}

If $I$ is one of the intervals $(a,\,b]$, $[a,\,b)$ or $(a,\,b)$, the inequalities in the definition of $Q_{I,\,P}^{p,\,q}$ are adapted accordingly in the obvious manner. 

\end{Def}

\vspace{2ex}

 From the discussion preceding Definition \ref{Def:H-Q_spaces_arbitrary}, we infer the following

\begin{Cor}\label{Cor:H-I-pq} Let $P:\mbox{Dom}\,P\longrightarrow L^2_{p,\,q}(X,\,L_k)$ be an operator as in Definition \ref{Def:H-Q_spaces_arbitrary}.

\vspace{1ex}

$(a)$\, For all intervals $I,J\subset\R$ such that $I\cap J=\emptyset$ and $I\cup J=\R$, the following orthogonal decomposition holds: \begin{equation*}L^2_{p,\,q}(X,\,L_k)={\cal H}_{I,\,P}^{p,\,q}(X,\,L_k)\oplus{\cal H}_{J,\,P}^{p,\,q}(X,\,L_k).\end{equation*}

\vspace{1ex}

\noindent $(b)$\, For all $-\infty\leq a\leq b\leq +\infty$, the following inclusion holds: \begin{equation*}{\cal H}^{p,\,q}_{[a,\,b],\,P}(X,\,L_k)\subset Q^{p,\,q}_{[a,\,b],\,P}(X,\,L_k).\end{equation*}

\vspace{1ex}

\noindent $(c)$\, For any interval $I\subset (-\infty,\,0)$, we have ${\cal H}_I^{p,\,q}(X,\,L_k)=\{0\}.$

\end{Cor}

\noindent {\it Proof.} We spell out the arguments in the case $P=\Delta''_k$. They are identical for other choices of $P$.

\vspace{1ex}

$\bullet$ By properties $(ii)$ and $(iii)$ of the spectral projections, we have $$\mbox{Id}_{L^2_{p,\,q}(X,\,L_k)} = \chi_{(-\infty,\,+\infty)}(\Delta''_k) = \chi_I(\Delta''_k) + \chi_J(\Delta''_k),$$

\noindent so taking images we get $L^2_{p,\,q}(X,\,L_k)={\cal H}_{I}^{p,\,q}(X,\,L_k) + {\cal H}_{I}^{p,\,q}(X,\,L_k)$. 

To prove ${\cal H}_{I}^{p,\,q}(X,\,L_k)\perp{\cal H}_{J}^{p,\,q}(X,\,L_k)$, we write for any $u,v\in L^2_{p,\,q}(X,\,L_k)$: \begin{eqnarray}\nonumber\langle\langle\chi_I(\Delta''_k)u,\,\chi_J(\Delta''_k)v\rangle\rangle & \stackrel{(1)}{=} & \langle\langle\chi_J(\Delta''_k)\chi_I(\Delta''_k)u,\,v\rangle\rangle \stackrel{(2)}{=}\langle\langle\chi_{I\cap J}(\Delta''_k)u,\,v\rangle\rangle\\
\nonumber   & = & \langle\langle\chi_{\emptyset}(\Delta''_k)u,\,v\rangle\rangle \stackrel{(3)}{=}0,\end{eqnarray}

\noindent where $(1)$ followed from $\chi_J(\Delta''_k)^{\star} = \chi_J(\Delta''_k)$, while $(2)$ (resp. $(3)$) followed from property $(iv)$ (resp. $(ii)$) of the spectral projections. This proves $(a)$.

\vspace{1ex}

$\bullet$ To prove $(b)$, recall that for every bounded Borel function $h:\R\rightarrow\R$, the bounded operator $h(\Delta''_k):L^2_{p,\,q}(X,\,L_k)\rightarrow L^2_{p,\,q}(X,\,L_k)$ is defined (cf. bounded measurable functional calculus in [RS80, chapter VIII]) by the formula: \begin{equation}\label{eqn:h-delta''}\langle\langle h(\Delta''_k)u,\,u\rangle\rangle = \int\limits_{-\infty}^{+\infty}h\,d\nu_{\Delta''_k,\,u} \hspace{5ex}\mbox{for all}\hspace{1ex}u\in L^2_{p,\,q}(X,\,L_k),\end{equation}

\noindent where $\nu_{\Delta''_k,\,u}$ is the {\it spectral measure} associated with $u$ and defined by

$$\nu_{\Delta''_k,\,u}(\Omega):=\langle\langle\chi_{\Omega}(\Delta''_k)u,\,u\rangle\rangle\in\R_{+}$$

\noindent for every Borel set $\Omega\subset\R$. 

In this case, $\nu_{\Delta''_k,\,u}$ is a {\it positive measure} since $\langle\langle\Delta''_ku,\,u\rangle\rangle\geq 0$ for all $u\in\mbox{Dom}\,\Delta''_k$ by (\ref{eqn:Delta''-scalarprod}) (itself a consequence of the completeness assumption on $\omega$). Taking $h\equiv 1$ in (\ref{eqn:h-delta''}), we get

\begin{equation}\label{eqn:u-spectral-measure}||u||^2 = \nu_{\Delta''_k,\,u}(\R) \hspace{2ex}\mbox{for all}\hspace{1ex} u\in L^2_{p,\,q}(X,\,L_k),\end{equation}

\noindent while taking $h=\chi_{\Omega}$ for an arbitrary Borel set $\Omega\subset\R$, we get

\begin{equation}\label{eqn:chi-spectral-measure}\langle\langle\chi_{\Omega}(\Delta''_k)u,\,u\rangle\rangle = \nu_{\Delta''_k,\,u}(\Omega) \hspace{2ex}\mbox{for all}\hspace{1ex} u\in L^2_{p,\,q}(X,\,L_k).\end{equation}
 
Meanwhile, formula (\ref{eqn:h-delta''}) extends to unbounded Borel functions $h$ (cf. [RS80, chapter VIII]), so taking $h=\mbox{Id}_{\R}$, we get \begin{equation}\label{eqn:id-spectral-measure}\langle\langle\Delta''_ku,\,u\rangle\rangle = \int\limits_{-\infty}^{+\infty}\mbox{Id}_{\R}\, d\nu_{\Delta''_k,\,u} = \int\limits_{-\infty}^{+\infty}\lambda\, d\nu_{\Delta''_k,\,u}(\lambda)\end{equation}

\noindent for all $u\in\mbox{Dom}\,\Delta''_k\subset L^2_{p,\,q}(X,\,L_k)$.

 Finally, we have the following easy observation.

\begin{Lem}\label{Lem:spectral-measures-intersection} If $u=\chi_{\Omega_2}(\Delta''_k)v$ for some $v\in L^2_{p,\,q}(X,\,L_k)$ and some Borel set $\Omega_2\subset\R$, then \begin{equation}\label{eqn:spectral-measures-intersection}\nu_{\Delta''_k,\,u}(\Omega_1) = \nu_{\Delta''_k,\,v}(\Omega_1\cap\Omega_2) \hspace{2ex}\mbox{for any Borel set}\hspace{1ex}\Omega_1\subset\R.\end{equation}

\end{Lem}

\noindent {\it Proof.} For any $u=\chi_{\Omega_2}(\Delta''_k)v$, we get: \begin{eqnarray*}\nu_{\Delta''_k,\,u}(\Omega_1) & \stackrel{(1)}{=} & \langle\langle\chi_{\Omega_1}(\Delta''_k)\chi_{\Omega_2}(\Delta''_k)v,\,\chi_{\Omega_2}(\Delta''_k)v\rangle\rangle\\
 & \stackrel{(2)}{=} & \langle\langle\chi_{\Omega_1\cap\Omega_2}(\Delta''_k)v,\,\chi_{\Omega_1\cap\Omega_2}(\Delta''_k)v\rangle\rangle + \langle\langle\chi_{\Omega_1\cap\Omega_2}(\Delta''_k)v,\,\chi_{\Omega_2\setminus\Omega_1}(\Delta''_k)v\rangle\rangle\\
  & \stackrel{(3)}{=} & \langle\langle\chi_{\Omega_1\cap\Omega_2}^2(\Delta''_k)v,\,v\rangle\rangle + \langle\langle\chi_{\Omega_2\setminus\Omega_1}(\Delta''_k)\chi_{\Omega_1\cap\Omega_2}(\Delta''_k)v,\,v\rangle\rangle\\
 & \stackrel{(4)}{=} & \langle\langle\chi_{\Omega_1\cap\Omega_2}(\Delta''_k)v,\,v\rangle\rangle + \langle\langle\chi_{\emptyset}(\Delta''_k)v,\,v\rangle\rangle\\ 
 & \stackrel{(5)}{=} & \langle\langle\chi_{\Omega_1\cap\Omega_2}(\Delta''_k)v,\,v\rangle\rangle = \nu_{\Delta''_k,\,v}(\Omega_1\cap\Omega_2),\end{eqnarray*}

\noindent where $(1)$ follows from (\ref{eqn:chi-spectral-measure}), $(2)$ follows from properties $(iv)$ and $(iii)$ of the spectral projections and from $\chi_{\Omega_2}(\Delta''_k)v = \chi_{\Omega_1\cap\Omega_2}(\Delta''_k)v + \chi_{\Omega_2\setminus\Omega_1}(\Delta''_k)v$, $(3)$ follows from the self-adjointness of the spectral projections, $(4)$ follows from properties $(i)$ and $(iv)$ of the spectral projections, while $(5)$ follows from property $(ii)$ and from (\ref{eqn:chi-spectral-measure}).  \hfill $\Box$

\vspace{2ex}

\noindent {\it End of proof of Corollary \ref{Cor:H-I-pq}.} Going back to the proof of $(b)$ of Corollary \ref{Cor:H-I-pq}, let us fix an arbitrary $u=\chi_{[a,\,b]}(\Delta''_k)v\in{\cal H}_{[a,\,b]}^{p,\,q}(X,\,L_k)$ with some $v\in L^2_{p,\,q}(X,\,L_k)$. Using (\ref{eqn:id-spectral-measure}) first and then (\ref{eqn:spectral-measures-intersection}) with $\Omega_2=[a,\,b]$, we get: $$\langle\langle\Delta''_ku,\,u\rangle\rangle = \int\limits_{-\infty}^{+\infty}\mbox{Id}_{\R}\,d\nu_{\Delta''_k,\,u} = \int\limits_a^b\mbox{Id}_{\R}\,d\nu_{\Delta''_k,\,v}  = \int\limits_a^b\lambda\,d\nu_{\Delta''_k,\,v}(\lambda).$$

 Since $\nu_{\Delta''_k,\,v}$ is a positive measure, we infer: \begin{equation}\label{eqn:Delta_kuu1}a\,\nu_{\Delta''_k,\,v}([a,\,b])\leq\langle\langle\Delta''_ku,\,u\rangle\rangle\leq b\,\nu_{\Delta''_k,\,v}([a,\,b]).\end{equation}

 Now, using (\ref{eqn:u-spectral-measure}) first and then (\ref{eqn:spectral-measures-intersection}) with $\Omega_1=\R$ and $\Omega_2=[a,\,b]$, we get: \begin{equation}\label{eqn:norm-u-squared}||u||^2 = \nu_{\Delta''_k,\,u}(\R) = \nu_{\Delta''_k,\,v}([a,\,b]).\end{equation}

 The combined (\ref{eqn:Delta_kuu1}) and (\ref{eqn:norm-u-squared}) give $$a\,||u||^2\leq\langle\langle\Delta''_ku,\,u\rangle\rangle\leq b\,||u||^2,$$

\noindent proving that $u\in Q^{p,\,q}_{[a,\,b]}(X,\,L_k)$. This proves $(b)$.

\vspace{1ex}

 Part $(c)$ is an immediate consequence of $(b)$ since $\langle\langle\Delta''_ku,\,u\rangle\rangle\geq 0$ for all $u\in\mbox{Dom}\,\Delta''_k$ (cf. (\ref{eqn:Delta''-scalarprod})). The proof of Corollary \ref{Cor:H-I-pq} is complete.  \hfill $\Box$

\vspace{3ex}

 The above construction and arguments for $P=\Delta''_k$ can be repeated for the twisted Laplacians. For example, for any ${\cal C}^\infty$ positive function $\eta>0$ on $X$, the $\eta$-twisted Laplacian extends to an unbounded, closed and densely defined operator

$$\Delta''_{k,\,(\eta)}\,:\,\mbox{Dom}\,\Delta''_{k,\,(\eta)}\longrightarrow L^2_{p,\,q}(X,\,L_k)$$

\noindent where $\mbox{Dom}\,\Delta''_{k,\,(\eta)}\subset L^2_{p,\,q}(X,\,L_k)$ consists of forms $u\in L^2_{p,\,q}(X,\,L_k)$ for which $\Delta''_{k,\,(\eta)}u\in L^2_{p,\,q}(X,\,L_k)$.

\begin{Obs}\label{Obs:equal-domains} If $0<C_1<\eta<C_2<\infty$ and $|\partial\eta|<C_2$ at every point of $X$ for some constants $C_1,C_2$, then

$$\mbox{Dom}\,\Delta''_{k,\,(\eta)} = \mbox{Dom}\,\Delta''_k.$$

\end{Obs}

\noindent {\it Proof.} Suppose $u\in\mbox{Dom}\,\Delta''_k\subset\mbox{Dom}\,\bar\partial_k\cap\mbox{Dom}\,\bar\partial_k^{\star}$. Then $\Delta''_ku$, $\bar\partial_ku$ and $\bar\partial_k^{\star}u$ are all $L^2$. The assumptions on $\eta$ ensure that $\eta\,\Delta''_ku$, $\bar\partial\eta\wedge\bar\partial_k^{\star}u$ and $(\bar\partial\eta\wedge\cdot)^{\star}\bar\partial_ku$ are still in $L^2$. It follows from formula (\ref{eqn:twisted-lap''-formula2}) that $\Delta''_{k,\,(\eta)}u$ is again in $L^2$, hence $u\in\mbox{Dom}\,\Delta''_{k,\,(\eta)}$. This proves the inclusion $\supset$.

 The reverse inclusion is proved similarly after noticing that $\Delta''_k$ and $\Delta''_{k,\,(\eta)}$ can be permuted in formula (\ref{eqn:twisted-lap''-formula2}) if $\eta$ is replaced with $1/\eta$.   \hfill $\Box$

\vspace{2ex}


In the $L^2$ extension context of $\S.$\ref{subsection:OTapplication}, to lighten the notation we set: $$\Delta''_{k,\,(\varepsilon)}:= \Delta''_{k,\,(\eta_{\varepsilon}+\lambda_{\varepsilon})}\hspace{5ex} \mbox{and} \hspace{5ex} \widetilde\Delta''_{k,\,(\varepsilon)}:= \widetilde\Delta''_{k,\,(\eta_{\varepsilon}+\lambda_{\varepsilon})}$$ for all $k\in\Sigma\subset\N^{\star}, \varepsilon>0$. Note that with the choice of ${\cal C}^\infty$ functions $\eta_{\varepsilon}$ and $\lambda_{\varepsilon}$ of $\S.$\ref{subsection:OTapplication}, the function $\eta_{\varepsilon} + \lambda_{\varepsilon}$ satisfies the assumptions of Observation \ref{Obs:equal-domains} for every $\varepsilon>0$ (with constants $C_1$ and $C_2$ depending on $\varepsilon$). Thus \begin{equation}\label{epsilon-equal-domains}\mbox{Dom}\,\Delta''_{k,\,(\varepsilon)} = \mbox{Dom}\,\Delta''_k \hspace{2ex}\mbox{for all}\hspace{1ex} \varepsilon>0, k\in\Sigma\subset\N^{\star}.\end{equation}

Definition \ref{Def:H-Q_spaces_arbitrary} specialises to the following

\begin{Def}\label{Def:Hk-def} (a)\, For every bidegree $(p,\,q)$ and every $k\in\Sigma\subset\N^{\star}$, we set: \begin{equation*}\label{eqn:Hk-def1}{\cal H}^{p,\,q}(X,\,L_k)={\cal H}^{p,\,q}_{[0,\,\delta_k],\,\Delta''_k}(X,\,L_k): = \mbox{Im}\,\bigg(\chi_{[0,\,\delta_k]}(\Delta''_k)\bigg),\end{equation*}  \begin{equation*}\label{eqn:Hk-def2}{\cal H}^{p,\,q}_{(\varepsilon)}(X,\,L_k)={\cal H}^{p,\,q}_{[0,\,\delta_k],\,\Delta''_{k,\,(\varepsilon)}}(X,\,L_k): = \mbox{Im}\,\bigg(\chi_{[0,\,\delta_k]}(\Delta''_{k,\,(\varepsilon)})\bigg),\end{equation*} and \begin{equation*}\label{eqn:Hk-def3}{\cal H}^{p,\,q}_{k,\,\varepsilon}(X,\,L_k)= {\cal H}^{p,\,q}_{[0,\,\delta_k],\,\widetilde\Delta''_{k,\,(\varepsilon)}}(X,\,L_k):= \mbox{Im}\,\bigg(\chi_{[0,\,\delta_k]}(\widetilde\Delta''_{k,\,(\varepsilon)})\bigg),\end{equation*} where $\delta_k>0$ is an arbitrary constant. 

 The elements of the vector space ${\cal H}^{p,\,q}(X,\,L_k)$ are called {\bf approximately holomorphic} $L_k$-valued $(p,\,q)$-forms on $X$. 

\vspace{2ex}

(b)\, In the set-up of (a), we also set \begin{equation}\label{eqn:Qk-def1}\nonumber Q^{p,\,q}(X,\,L_k):=Q^{p,\,q}_{[0,\,\delta_k],\,\Delta''_k}(X,\,L_k)=\bigg\{u\in\mbox{Dom}\,\Delta''_k\,;\, 0\leq\langle\langle\Delta''_ku,\,u\rangle\rangle\leq\delta_k\,||u||^2\bigg\}\end{equation} and \begin{equation}\label{eqn:Qk-def2}\nonumber Q^{p,\,q}_{(\varepsilon)}(X,\,L_k):=Q^{p,\,q}_{[0,\,\delta_k],\,\Delta''_{k,\,(\varepsilon)}}(X,\,L_k)=\bigg\{u\in\mbox{Dom}\,\Delta''_{k,\,(\varepsilon)}\,;\, 0\leq\langle\langle\Delta''_{k,\,(\varepsilon)}u,\,u\rangle\rangle\leq\delta_k\,||u||^2\bigg\}.\end{equation}

\end{Def}

The first properties of these spaces are summed up in the following

\begin{Obs}\label{Obs:inclusion_Qs} For all $k\in\Sigma\subset\N^{\star}$ and $\varepsilon>0$, we have 

\vspace{1ex}

\noindent $(i)$\, ${\cal H}^{p,\,q}(X,\,L_k)\subset Q^{p,\,q}(X,\,L_k)$ and ${\cal H}^{p,\,q}_{(\varepsilon)}(X,\,L_k)\subset Q^{p,\,q}_{(\varepsilon)}(X,\,L_k)$; \\

\noindent $(ii)$\, $\Delta''_{k,\,(\varepsilon)}\geq\Delta''_k$ and \begin{equation}\label{eqn:ahspace-incl}Q^{p,\,q}_{(\varepsilon)}(X,\,L_k)\subset Q^{p,\,q}(X,\,L_k);\end{equation}



\end{Obs}

\noindent {\it Proof.} Part $(i)$ follows from (b) of Corollary \ref{Cor:H-I-pq}, while using (\ref{eqn:twisted-lap''-scal}) and the fact that $\eta_{\varepsilon} + \lambda_{\varepsilon}>1$ on $X$, we see that \begin{eqnarray}\nonumber\langle\langle\Delta''_{k,\,(\varepsilon)}v,\,v\rangle\rangle & = & \int\limits_X(\eta_{\varepsilon} + \lambda_{\varepsilon})\,|\bar\partial_kv|^2_{h_k}\,dV_{\omega} + \int\limits_X(\eta_{\varepsilon} + \lambda_{\varepsilon})\,|\bar\partial_k^{\star}v|^2_{h_k}\,dV_{\omega}\\
\nonumber & \geq & \int\limits_X|\bar\partial_kv|^2_{h_k}\,dV_{\omega} + \int\limits_X|\bar\partial_k^{\star}v|^2_{h_k}\,dV_{\omega}  = \langle\langle\Delta''_k v,\,v\rangle\rangle\end{eqnarray}

\noindent for $v\in{\cal D}^{p,\,q}(X,\,L_k)$ (hence for $v\in\mbox{Dom}\,\Delta''_k$) and $\varepsilon>0$. This implies part $(ii)$.    \hfill $\Box$

\begin{Obs}\label{Obs:C_infty_Hs} Fix an arbitrary bidegree $(p,\,q)$. For all $k\in\Sigma\subset\N^{\star}$ and $\varepsilon>0$, the following inclusions hold  

\begin{equation}\label{eqn:C_infty_Hs}{\cal H}^{p,\,q}(X,\,L_k)\subset {\cal C}^\infty_{p,\,q}(X,\,L_k) \hspace{3ex} \mbox{and} \hspace{3ex} {\cal H}^{p,\,q}_{(\varepsilon)}(X,\,L_k)\subset {\cal C}^\infty_{p,\,q}(X,\,L_k).\end{equation}

\noindent Moreover, ${\cal H}^{p,\,q}(X,\,L_k)$ and ${\cal H}^{p,\,q}_{(\varepsilon)}(X,\,L_k)$ are closed subspaces of $L^2_{p,\,q}(X,\,L_k)$.

 More generally, for every {\bf finite-length} interval $I\subset\R$, the space 
${\cal H}_{I,\,P}^{p,\,q}(X,\,L_k)$ introduced in Definition \ref{Def:H-Q_spaces_arbitrary} is contained in ${\cal C}^\infty_{p,\,q}(X,\,L_k)$.

\end{Obs}

\noindent {\it Proof.} Let $k$ and $\varepsilon$ be fixed as in the statement. We first note that, by taking $h=x^l\chi_I$ in (\ref{eqn:h-delta''}), we get the following inclusion for every finite-length interval $I\subset\R$: \begin{equation}\label{eqn:image_chiI_domain_inclusion}\mbox{Im}\,\chi_I(\Delta''_k)\subset\bigcap\limits_{l=1}^{+\infty}\mbox{Dom}\,((\Delta''_k)^l).\end{equation} The analogous inclusion with $\Delta''_{k,\,(\varepsilon)}$ in place of $\Delta''_k$ holds as well. 

 Then, the ellipticity of $\Delta''_k$ implies the inclusion \begin{equation}\label{eqn:domain_inclusion_C_infty}\bigcap\limits_{l=1}^{+\infty}\mbox{Dom}\,((\Delta''_k)^l)\subset {\cal C}^\infty_{p,\,q}(X,\,L_k),\end{equation}  

 \noindent while  the ellipticity of $\Delta''_{k,\,(\varepsilon)}$ implies the analogous inclusion with $\Delta''_{k,\,(\varepsilon)}$ in place of $\Delta''_k$.

 Thus, (\ref{eqn:image_chiI_domain_inclusion}) and (\ref{eqn:domain_inclusion_C_infty}) yield ${\cal H}_I^{p,\,q}(X,\,L_k)\subset {\cal C}^\infty_{p,\,q}(X,\,L_k)$. \hfill $\Box$

 \vspace{3ex}

In anticipation of further use, we now give the version that is relevant to our setting of a standard inequality. The Laplacian $\Delta''_k$ in the next statement can be replaced by any non-negative self-adjoint elliptic operator.

\begin{Lem}\label{Lem:iterated-laplacian} For any $k\in\Sigma\subset\N^\star$ and any constant $\delta_k>0$, the following inequality holds: \begin{eqnarray}\label{eqn:iterated-laplacian}\langle\langle(\Delta''_k)^lu,\,u\rangle\rangle\leq\delta_k^l\,||u||^2,\end{eqnarray} for every positive integer $l$ and every $u\in{\cal H}^{n,\,0}_{[0,\,\delta_k],\,\Delta''_k}(X,\,L_k)$.

\end{Lem}  

\noindent {\it Proof.} When $l=1$, (\ref{eqn:iterated-laplacian}) follows from the definition of ${\cal H}^{n,\,0}_{[0,\,\delta_k],\,\Delta''_k}(X,\,L_k)$.

We can continue by induction on $l\geq 1$. We have:  \begin{eqnarray*}\langle\langle(\Delta''_k)^lu,\,u\rangle\rangle & = & \bigg\langle\bigg\langle\Delta''_k\bigg((\Delta''_k)^{\frac{l-1}{2}}u\bigg),\,(\Delta''_k)^{\frac{l-1}{2}}u\bigg\rangle\bigg\rangle \stackrel{(a)}{\leq} \delta_k\,\bigg|\bigg|(\Delta''_k)^{\frac{l-1}{2}}u\bigg|\bigg|^2 = \delta_k\,\bigg\langle\bigg\langle(\Delta''_k)^{l-1}u,\,u\bigg\rangle\bigg\rangle \\
& \leq & \delta_k\,\delta_k^{l-1}||u||^2,\end{eqnarray*} where the last inequality follows from the induction hypothesis, while inequality (a) follows from $(\Delta''_k)^{\frac{l-1}{2}}u\in{\cal H}^{n,\,0}_{[0,\,\delta_k],\,\Delta''_k}(X,\,L_k)$. (Indeed, applying any power of $\Delta''_k$ to $u$ keeps $u$ in ${\cal H}^{n,\,0}_{[0,\,\delta_k],\,\Delta''_k}(X,\,L_k)$.)  \hfill $\Box$

\vspace{3ex}

All the definitions and results of this $\S.$\ref{subsection:definitions_approx-hol} still make sense if $X$ is replaced by $Y$. We will use the versions for $Y$ at the appropriate spots in the sequel.

\subsection{The commutation defect}\label{subsection:commutation-defect} As above, $(X,\,\omega)$ will be a complete K\"ahler manifold with a fixed ${\cal C}^\infty$ $2$-form $\alpha$ satisfying Assumption A. However, for the sake of generality, we will write the calculations in this subsection for an arbitrary Hermitian metric $\omega$. 

Since $\bar\partial^2_k\neq 0$, the twisted Laplacian $\Delta''_{k,\,(\eta)}$ of Definition \ref {Def:twisted-d-bar-laplacian} does not commute with $\bar\partial_k^{\eta}$. We will estimate the commutation defect in bidegree $(n,\, 0)$ in the $L^2$-norm in a way similar to Laeng's estimation of the commutation defect of $\Delta''_k$ with $\bar\partial_k$ in bidegree $(0,\, 0)$ (cf. [Lae02, $\S.4.2.1$]).

\vspace{1ex}

 $\bullet$ We first analyse the untwisted situation. The following estimate in bidegree $(n,\,0)$ is the same as Laeng's in bidegree $(0,\, 0)$.

\begin{Lem}\label{Lem:commutation-defect_untwisted} On any complete $n$-dimensional Hermitian manifold $(X,\,\omega)$, the following inequality holds in the $L^2$-norm: \begin{equation}\label{eqn:comm-defect-ntwisted}||(\Delta''_k\bar\partial_k - \bar\partial_k\Delta''_k)\,s||^2 = ||\bar\partial_k^{\star}\bar\partial_k^2s||^2\leq\frac{C}{k^{2/b_2}}\,(k||s||^2 + ||\bar\partial_ks||^2)\end{equation}

\noindent for all $s\in{\cal D}^{n,\,0}(X,\,L_k)$.

\end{Lem}

\noindent {\it Proof.} For $s\in {\cal C}^\infty_{n,\, 0}(X,\, L_k)$, we have

$$(\Delta''_k\bar\partial_k - \bar\partial_k\Delta''_k)\,s = \bar\partial_k^{\star}\bar\partial_k^2s.$$

\noindent Now, $\bar\partial_k^{\star}\bar\partial_k^2s = -2\pi i\,\bar\partial_k^{\star}(\alpha_k^{0,\, 2}\wedge s)$. Thus, if $s$ has compact support, we get: \begin{eqnarray}\label{eqn:ntwisted-est1}\nonumber||\bar\partial_k^{\star}\bar\partial_k^2s||^2 & = & (2\pi)^2\,||\bar\partial_k^{\star}(\alpha_k^{0,\, 2}\wedge s)||^2 \stackrel{(a)}{\leq} (2\pi)^2\,\langle\langle\Delta''_k(\alpha_k^{0,\, 2}\wedge s), \, \alpha_k^{0,\, 2}\wedge s\rangle\rangle\\
\nonumber  & \stackrel{(b)}{=} & (2\pi)^2\,||\partial_k(\alpha_k^{0,\, 2}\wedge s) + \tau (\alpha_k^{0,\, 2}\wedge s)||^2\\
 \nonumber   & + & (2\pi)^2\,||\partial_k^{\star}(\alpha_k^{0,\, 2}\wedge s) + \tau^{\star}(\alpha_k^{0,\, 2}\wedge s)||^2\\
\nonumber & + & (2\pi)^3 \, \langle\langle[\alpha_k^{1,\,1},\,\Lambda](\alpha_k^{0,\, 2}\wedge s),\, \alpha_k^{0,\, 2}\wedge s)\rangle\rangle\\
       & + & \langle\langle T_{\omega}(\alpha_k^{0,\, 2}\wedge s),\, \alpha_k^{0,\, 2}\wedge s)\rangle\rangle,\end{eqnarray}

\noindent where inequality $(a)$ follows from $\langle\langle\Delta''_k(\alpha_k^{0,\, 2}\wedge s), \, \alpha_k^{0,\, 2}\wedge s\rangle\rangle = ||\bar\partial_k(\alpha_k^{0,\, 2}\wedge s)||^2 + ||\bar\partial_k^{\star}(\alpha_k^{0,\, 2}\wedge s)||^2 $ which holds for any compactly supported $s$ thanks to the completeness of $\omega$, while identity $(b)$ follows from $s$ being compactly supported and from the Bochner-Kodaira-Nakano identity $\Delta''_k = \Delta'_{k,_,\tau} + [i\Theta(D_k)^{1,\, 1},\,\Lambda] + T_{\omega}$ (see [Lae02]) induced by an arbitrary (i.e. possibly non-K\"ahler) metric $\omega$. Here $\tau:=[\Lambda,\,\partial\omega]$ is the torsion $(1,\, 0)$-form of $\omega$ (which vanishes if $\omega$ is K\"ahler), $\Delta'_{k,_,\tau}:=[\partial_k+\tau,\, \partial_k^{\star} + \tau^{\star}]$ and $T_{\omega}:=[\Lambda,\,[\Lambda,\,\frac{i}{2}\,\partial\bar\partial\omega]] - [\partial\omega,\,(\partial\omega)^{\star}]$ (which also vanishes if $\omega$ is K\"ahler).

 Now $\partial_k(\alpha_k^{0,\,2}\wedge s) = (\partial\alpha_k^{0,\, 2})\wedge s$ since $\partial_ks=0$ for bidegree reasons, so $||\partial_k(\alpha_k^{0,\,2}\wedge s)||^2\leq\frac{C}{k^{2/b_2}}\,||s||^2$. Meanwhile, the Hermitian commutation relation $-i\,(\partial_k^{\star} + \tau^{\star}) = [\Lambda,\,\bar\partial_k]$ (see Proposition \ref{Prop:Hermitian_commutation-rel}) yields: \begin{eqnarray}\label{eqn:del-star-alpha-est}\partial_k^{\star}(\alpha_k^{0,\,2}\wedge s) & = & i[\Lambda,\,\bar\partial_k](\alpha_k^{0,\,2}\wedge s) - \tau^{\star}(\alpha_k^{0,\, 2}\wedge s)\\
 \nonumber & = & i\Lambda\bar\partial_k(\alpha_k^{0,\,2}\wedge s) - i\bar\partial_k\Lambda(\alpha_k^{0,\,2}\wedge s) - \tau^{\star}(\alpha_k^{0,\, 2}\wedge s)\\
\nonumber  & = & i\Lambda(\bar\partial\alpha_k^{0,\,2}\wedge s + \alpha_k^{0,\, 2}\wedge\bar\partial_ks) - i\bar\partial_k\Lambda(\alpha_k^{0,\,2}\wedge s) - \tau^{\star}(\alpha_k^{0,\, 2}\wedge s).\end{eqnarray}

 Since $\Lambda$, $\tau$ and $T_{\omega}$ are of order zero and independent of $k$, we get:

\vspace{2ex}

\hspace{5ex} $||i\bar\partial_k\Lambda(\alpha_k^{0,\,2}\wedge s)||^2,\hspace{2ex} ||i\Lambda(\bar\partial\alpha_k^{0,\,2}\wedge s + \alpha_k^{0,\, 2}\wedge\bar\partial_ks)||^2 \leq\frac{C}{k^{2/b_2}}\,(||s||^2 + ||\bar\partial_ks||^2),$

\vspace{1ex}

\hspace{5ex} $||\tau(\alpha_k^{0,\, 2}\wedge s)||^2, \hspace{2ex} ||\tau^{\star}(\alpha_k^{0,\, 2}\wedge s)||^2\leq\frac{C}{k^{2/b_2}}\,||s||^2,$ 

\vspace{1ex}

\hspace{5ex} $|\langle\langle T_{\omega}(\alpha_k^{0,\, 2}\wedge s),\, \alpha_k^{0,\, 2}\wedge s)\rangle\rangle|\leq\frac{C}{k^{2/b_2}}\,||s||^2.$

\vspace{2ex}

\noindent Hence, from this and from (\ref{eqn:del-star-alpha-est}) we get: $$||\partial_k^{\star}(\alpha_k^{0,\,2}\wedge s)||^2\leq\frac{C}{k^{2/b_2}}\,(||s||^2 + ||\bar\partial_ks||^2).$$

\vspace{1ex}

Finally, since $\alpha_k^{1,\, 1}$ is close in the ${\cal C}^\infty$ topology to $k\alpha$, we get: $$|\langle\langle[\alpha_k^{1,\,1},\,\Lambda](\alpha_k^{0,\, 2}\wedge s),\, \alpha_k^{0,\, 2}\wedge s)\rangle\rangle|\leq\frac{C}{k^{2/b_2}}\,k||s||^2$$ since $||\alpha_k^{0,\, 2}\wedge s||^2\leq\frac{C}{k^{2/b_2}}\,||s||^2.$

\vspace{2ex}

 The proof is complete.   \hfill $\Box$

\vspace{2ex}

$\bullet$ We now go to the twisted situation. We shall need the following brief computation.

\begin{Lem}\label{Lem:dbar-star-prod} If $\rho\,:\,X\rightarrow\R$ is any $C^1$ function, $q\geq 0$ is an integer and $\gamma\in {\cal C}^\infty_{n,\,q}(X,\,L_k)$, then $$\bar\partial_k^{\star}(\rho\gamma) = \rho\bar\partial_k^{\star}\gamma + i(\partial\rho)\wedge\Lambda\gamma,$$

\noindent where $\Lambda=\Lambda_\omega$ is the adjoint of the operator of multiplication by the fixed Hermitian metric $\omega$ w.r.t. the pointwise inner product induced by $\omega$.

More generally, in arbitrary bidegree $(p,\,q)$, we have

$$\bar\partial_k^{\star}(\rho\gamma) = \rho\bar\partial_k^{\star}\gamma -[\Lambda,\, i\partial\rho\wedge\cdot]\,(\gamma),  \hspace{3ex} \gamma\in{\cal C}^\infty_{p,\,q}(X,\,L_k).$$

\end{Lem}

\noindent {\it Proof.} We will only prove the formula in bidegree $(n,\,q)$, the only one that will be needed in this paper. By the Hermitian commutation relation $i\,(\bar\partial_k^{\star} + \bar\tau^{\star}) = [\Lambda,\,\partial_k]$ (see Proposition\ref{Prop:Hermitian_commutation-rel}), we get: \begin{eqnarray}\label{eqn:dbar-star-calc}\nonumber\bar\partial_k^{\star}(\rho\,\gamma) & = & -i\,[\Lambda,\,\partial_k](\rho\,\gamma) - \bar\tau^{\star}(\rho\,\gamma) = -i\,\bigg(\Lambda\partial_k(\rho\,\gamma) - \partial_k\Lambda(\rho\,\gamma)\bigg) - \rho\,\bar\tau^{\star}(\gamma)\\
\nonumber & \stackrel{(a)}{=} & i\,\partial_k(\rho\,\Lambda\gamma) - \rho\,\bar\tau^{\star}(\gamma) = i\,(\partial\rho)\wedge\Lambda\gamma + i\rho\,\partial_k\Lambda\gamma - \rho\,\bar\tau^{\star}(\gamma)\\
\nonumber & \stackrel{(b)}{=} &  i\,(\partial\rho)\wedge\Lambda\gamma + \rho\,\bigg(i\,[\partial_k,\,\Lambda]\gamma - \bar\tau^{\star}(\gamma)\bigg) = i\,(\partial\rho)\wedge\Lambda\gamma + \rho\,\bar\partial_k^{\star}\gamma,\end{eqnarray}

\noindent where we have used the fact that $\partial_k(\rho\gamma)=\partial_k\gamma=0$ for bidegree reasons (since $\rho\gamma$ and $\gamma$ are of bidegree $(n,\,q)$), so neglecting the term $\Lambda\partial_k(\rho\gamma)=0$ gave $(a)$ and adding the term $\Lambda\partial_k\gamma=0$ gave $(b)$.   \hfill $\Box$

\vspace{2ex}

We now fix a $C^\infty$ positive function $\eta$ on $X$ and turn to the sequence (\ref{eqn:eta-twisted-dbar-seq_n1}). As already argued, this sequence imposes the $\eta$-twisted Laplacians: \begin{equation*}\Delta''_{k,\,(\eta)}:= (\bar\partial_k^{\eta})(\bar\partial_k^{\eta})^{\star} + (^{\eta}\bar\partial_k)^{\star}(^{\eta}\bar\partial_k):{\cal C}^\infty_{n,\,1}(X,\,L_k)\longrightarrow {\cal C}^\infty_{n,\,1}(X,\,L_k)\end{equation*} defined in (\ref{eqn:twistedDelta''def}) in bidegree $(n,\,1)$ and $\widetilde\Delta''_{k,\,(\eta)}:= (\bar\partial_k^{\eta})^{\star}(\bar\partial_k^{\eta}):{\cal C}^\infty_{n,\,0}(X,\,L_k)\longrightarrow {\cal C}^\infty_{n,\,0}(X,\,L_k)$ introduced in Definition \ref{Def:twisted-d-bar-laplacian_tilde} in bidegree $(n,\,0)$.

\begin{Prop}\label{Prop:comm-defect-twisted} On any complete $n$-dimensional Hermitian manifold $(X,\,\omega)$, for every $k$ and every $s\in{\cal D}^{n,\,0}(X,\,L_k)$ we have the following $L^2$-estimate: $$\bigg|\bigg|\bigg(\Delta''_{k,\,(\eta)}\bar\partial_k^{\eta}-\bar\partial_k^{\eta}\widetilde\Delta''_{k,\,(\eta)}\bigg)\,s\bigg|\bigg|^2 \leq \max\bigg((\sup\eta)^2,\,\sup|\partial\eta|_\omega^2\bigg)\cdot\frac{C}{k^{\frac{2}{b_2}}}\,\bigg(k||\sqrt{\eta}\,s||^2 + ||\bar\partial_k(\sqrt{\eta}\,s)||^2\bigg),$$

\noindent where $C>0$ is a constant independent of $k$, $\eta$ and $s$, while $||\,\,||$ is the $L^2$-norm w.r.t. $\omega$ and $h_k$.

\end{Prop}

\noindent {\it Proof.} For any $s\in {\cal C}^\infty_{n,\,0}(X,\,L_k)$, we have: \begin{eqnarray}\label{eqn:comm-defect-twisted1}\nonumber(\Delta''_{k,\,(\eta)}\bar\partial_k^{\eta}-\bar\partial_k^{\eta}\widetilde\Delta''_{k,\,(\eta)})\,s & = & (^{\eta}\bar\partial_k)^{\star}\,(^{\eta}\bar\partial_k)(\bar\partial_k^{\eta})s = \bar\partial_k^{\star}\bigg(\eta\,\bar\partial_k^2(\sqrt{\eta}\,s)\bigg)\\
\nonumber & = & -2\pi i\,\bar\partial_k^{\star}\bigg(\eta\,\alpha_k^{0,\,2}\wedge(\sqrt{\eta}\,s)\bigg).\end{eqnarray}

\noindent Combining this with Lemma \ref{Lem:dbar-star-prod}, we get: $$(\Delta''_{k,\,(\eta)}\bar\partial_k^{\eta}-\bar\partial_k^{\eta}\widetilde\Delta''_{k,\,(\eta)})\,s = \eta\,\bar\partial_k^{\star}\bar\partial_k^2(\sqrt{\eta}\,s) + 2\pi\,(\partial\eta)\wedge\Lambda(\alpha_k^{0,\,2}\wedge(\sqrt{\eta}\,s))$$ for all $s\in {\cal C}^\infty_{n,\,0}(X,\,L_k)$. Now, the $L^2$-norm of $\bar\partial_k^{\star}\bar\partial_k^2(\sqrt{\eta}\,s)$ was estimated in (\ref{eqn:comm-defect-ntwisted}), while the squared $L^2$-norm of $\Lambda(\alpha_k^{0,\,2}\wedge(\sqrt{\eta}\, s))$ is dominated by $C/k^{\frac{2}{b_2}}\cdot||\sqrt{\eta}\,s||^2$  when $s$ is compactly supported.   \hfill $\Box$

\vspace{3ex}

 It remains to estimate the quantities depending on the auxiliary function $\eta$ when this is chosen to be $\eta=\eta_\varepsilon + \lambda_\varepsilon$ with $\eta_\varepsilon$ and $\lambda_\varepsilon$ introduced in (\ref{eqn:eta-lambda-eps-def}) to suit the $L^2$ extension procedure.

\begin{Prop}\label{Prop:comm-defect-twisted_epsilon} Let $(X,\,\omega)$ be a complete Hermitian manifold of dimension $n$ which has a smooth complex hypersurface $Y=\{w=0\}\subset X$ defined by a bounded holomorphic function $w:X\longrightarrow\C$ with $\sup_X|w|\leq 1$. Let $\eta_\varepsilon:=\log (A/(|w|^2+\varepsilon^2))$ and $\lambda_\varepsilon:=1/(|w|^2+\varepsilon^2) $ be ${\cal C}^\infty$ functions on $X$ defined for every small $\varepsilon>0$ by a constant $A>e$. 

For every $k$ and every $s\in{\cal D}^{n,\,0}(X,\,L_k)$, the following estimate holds in the $L^2$-norm: \begin{eqnarray}\nonumber & & \bigg|\bigg|\bigg(\Delta''_{k,\,(\eta_\varepsilon + \lambda_\varepsilon)}\bar\partial_k^{\eta_\varepsilon + \lambda_\varepsilon}-\bar\partial_k^{\eta_\varepsilon + \lambda_\varepsilon}\widetilde\Delta''_{k,\,(\eta_\varepsilon + \lambda_\varepsilon)}\bigg)\,s\bigg|\bigg|^2 \\
\nonumber & \leq & \max\bigg\{\sup\limits_X|dw|^2_\omega,\,1\bigg\}\,\frac{1}{\varepsilon^6}\,\frac{C}{k^{\frac{2}{b_2}}}\,\bigg(k\,||\sqrt{\eta_\varepsilon + \lambda_\varepsilon}\,s||^2 + ||\bar\partial_k(\sqrt{\eta_\varepsilon + \lambda_\varepsilon}\,s)||^2\bigg)\\
\nonumber & \leq & \max\bigg\{\sup\limits_X|dw|^2_\omega,\,1\bigg\}\,\frac{1}{\varepsilon^6}\,\frac{C}{k^{\frac{2}{b_2}}}\,\bigg(\frac{k}{\varepsilon^3}\,||s||^2 + ||\bar\partial_k(\sqrt{\eta_\varepsilon + \lambda_\varepsilon}\,s)||^2\bigg)\end{eqnarray} for every $\varepsilon>0$ small enough, where $C>0$ is a constant independent of $\varepsilon$, $k$ and $s$.

\end{Prop}

\noindent{\it Proof.} We will derive this estimate from the one given in Proposition \ref{Prop:comm-defect-twisted} after we have estimated $\sup\,(\eta_\varepsilon+ \lambda_\varepsilon)$ and $\sup|\partial\eta_\varepsilon + \partial\lambda_\varepsilon|_\omega$.

From the definitions of these functions, we get: \begin{eqnarray}\label{eqn:comm-defect-twisted_epsilon_proof_1}\bigg(\sup\,(\eta_\varepsilon + \lambda_\varepsilon)\bigg)^2\leq\bigg(\log\frac{A}{\varepsilon^2} + \frac{1}{\varepsilon^2}\bigg)^2\leq\frac{4}{\varepsilon^4}<\frac{1}{\varepsilon^6}\end{eqnarray} whenever $\varepsilon>0$ is small enough.

Meanwhile, a direct computation shows that $\partial\eta_\varepsilon = -\bigg(\overline{w}/(|w|^2 + \varepsilon^2)\bigg)\,dw$, hence \begin{eqnarray*}|\partial\eta_\varepsilon|_\omega = \frac{|w|}{|w|^2 + \varepsilon^2}\,|dw|_\omega\leq\frac{1}{2\varepsilon}\,|dw|_\omega.\end{eqnarray*} Similarly, $\partial\lambda_\varepsilon = -\bigg(\overline{w}/(|w|^2 + \varepsilon^2)^2\bigg)\,dw$, hence \begin{eqnarray*}|\partial\lambda_\varepsilon|_\omega = \frac{|w|}{(|w|^2 + \varepsilon^2)^2}\,|dw|_\omega\leq\frac{1}{2\varepsilon}\,\frac{1}{|w|^2 + \varepsilon^2}\,|dw|_\omega\leq\frac{1}{2\varepsilon^3}\,|dw|_\omega.\end{eqnarray*} Thus, we get: \begin{eqnarray}\label{eqn:comm-defect-twisted_epsilon_proof_2}\sup|\partial\eta_\varepsilon + \partial\lambda_\varepsilon|_\omega^2 \leq 2\,\bigg(\frac{1}{4\varepsilon^2} + \frac{1}{4\varepsilon^6}\bigg)\,\sup|dw|^2_\omega \leq\frac{1}{\varepsilon^6}\,\sup|dw|^2_\omega\end{eqnarray} whenever $\varepsilon>0$ is small enough.

Putting together estimates (\ref{eqn:comm-defect-twisted_epsilon_proof_1}) and (\ref{eqn:comm-defect-twisted_epsilon_proof_2}) and taking $\eta=\eta_\varepsilon + \lambda_\varepsilon$ in the estimate given by Proposition \ref{Prop:comm-defect-twisted}, we get the first stated inequality. The second stated inequality follows from the first after a further application of (\ref{eqn:comm-defect-twisted_epsilon_proof_1}).  \hfill $\Box$

\subsection{The spectra of the twisted Laplacians in bidegrees $(n,\,0)$ and $(n,\,1)$}\label{subsection:spectra_n0-n1} The lower bound in bidegree $(n,\,1)$ obtained for the twisted Laplacian $\Delta''_{k,\,(\eta_\varepsilon+\lambda_\varepsilon)}:=\Delta''_{k,\,(\varepsilon)}$ (defined in (\ref{eqn:twistedDelta''def})) by taking $q=1$ in (i) of Proposition \ref{Prop:BKN-consequence} implies that the spectrum of this operator is contained in the interval $[k,\,+\infty)$. This may prove especially useful if we let $k\to\infty$ in future work. 

Before proving our main a priori $L^2$-estimate, we now use this information on the spectrum of $\Delta''_{k,\,(\varepsilon)}$ in bidegree $(n,\,1)$ and the commutation defect estimated in Proposition \ref{Prop:comm-defect-twisted_epsilon} to get information about the spectrum of $\widetilde\Delta''_{k,\,(\eta_\varepsilon+\lambda_\varepsilon)}:=\widetilde\Delta''_{k,\,(\varepsilon)}$ (defined in (\ref{eqn:twistedDelta''_tilde_def})) in bidegree $(n,\,0)$. The approach is similar, albeit with extra technical difficulties, to the one we used to prove Corollary 2.2. in [Pop13] which had itself been inspired by [Lae02]. 

Let $0<\lambda_1<\lambda_2$ and $\varepsilon_0'>0$ be constants that will be specified later. We consider the intervals $I=(\lambda_1,\,\lambda_2],\,J=[0,\,\lambda_2+\varepsilon_0']\subset\R$ and the vector subspaces \begin{eqnarray*}\widetilde{E}_I^{n,\,0} = \mbox{Im}\,\bigg(\chi_I(\widetilde\Delta''_{k,\,(\varepsilon)})\bigg)\subset L^2_{n,\,0}(X,\,L_k) \hspace{2ex}\mbox{and}\hspace{2ex} E_J^{n,\,1} = \mbox{Im}\,\bigg(\chi_J(\Delta''_{k,\,(\varepsilon)})\bigg)\subset L^2_{n,\,1}(X,\,L_k).\end{eqnarray*} If $X$ is compact, these are the sums of the eigenspaces of $\widetilde\Delta''_{k,\,(\varepsilon)}$ acting in bidegree $(n,\,0)$, resp. $\Delta''_{k,\,(\varepsilon)}$ acting in bidegree $(n,\,1)$, corresponding to the eigenvalues lying in $I$, resp. $J$. 

We will now be looking to determine appropriate choices of $\lambda_1$, $\lambda_2$ and $\varepsilon_0'$ such that the map \begin{eqnarray}\label{eqn:map_spectra_injective}\chi_J(\Delta''_{k,\,(\varepsilon)})\circ\bar\partial_k^{\eta_\varepsilon + \lambda_\varepsilon}:\widetilde{E}_I^{n,\,0}\longrightarrow  E_J^{n,\,1}\end{eqnarray} is {\it injective} when either $k$ is large enough or $\varepsilon>0$ is not too small relative to $k$. (Recall that in the special case when $X$ is compact, the map $\chi_J(\Delta''_{k,\,(\varepsilon)})$ is the projection onto the sum of the eigenspaces of $\Delta''_{k,\,(\varepsilon)}$ corresponding to the eigenvalues lying in $J$ -- see $\S.$\ref{subsection:definitions_approx-hol}.) 

We will reason by contradiction. Suppose the map (\ref{eqn:map_spectra_injective}) is not injective. Then, there exists $s\in\widetilde{E}_I^{n,\,0}\setminus\{0\}$ such that $\chi_J(\Delta''_{k,\,(\varepsilon)})(\bar\partial_k^{\eta_\varepsilon + \lambda_\varepsilon}s) = 0$. The last equality implies that \begin{eqnarray}\label{eqn:map_spectra_injective_proof_1}\bigg|\bigg|\Delta''_{k,\,(\varepsilon)}(\bar\partial_k^{\eta_\varepsilon + \lambda_\varepsilon}s)\bigg|\bigg|\geq(\lambda_2+\varepsilon_0')\,||\bar\partial_k^{\eta_\varepsilon + \lambda_\varepsilon}s||.\end{eqnarray}

On the other hand, since $s$ is of bidegree $(n,\,0)$, we get the equality below: \begin{eqnarray}\label{eqn:map_spectra_injective_proof_2}||\bar\partial_k^{\eta_\varepsilon + \lambda_\varepsilon}s||^2 = \langle\langle\widetilde\Delta''_{k,\,(\varepsilon)}s,\,s\rangle\rangle \geq \lambda_1\,||s||^2,\end{eqnarray} where the inequality follows $s\in\widetilde{E}_I^{n,\,0}$ and from $I=(\lambda_1,\,\lambda_2]$. 

Unlike the lower bound obtained in (\ref{eqn:map_spectra_injective_proof_1}), we now get an upper bound for the $L^2$-norm of $\Delta''_{k,\,(\varepsilon)}(\bar\partial_k^{\eta_\varepsilon + \lambda_\varepsilon}s)$ by first using the triangle inequality and then the commutation defect given by Proposition \ref{Prop:comm-defect-twisted_epsilon}: \begin{eqnarray*}\bigg|\bigg|\Delta''_{k,\,(\varepsilon)}(\bar\partial_k^{\eta_\varepsilon + \lambda_\varepsilon}s)\bigg|\bigg|\leq\frac{CC_1}{\varepsilon^3\,k^{\frac{1}{b_2}}}\,\bigg(\frac{k}{\varepsilon^3}\,||s||^2 + ||\bar\partial_k^{\eta_\varepsilon + \lambda_\varepsilon}\,s)||^2\bigg)^{\frac{1}{2}} + \bigg|\bigg|\bar\partial_k^{\eta_\varepsilon + \lambda_\varepsilon}\widetilde\Delta''_{k,\,(\varepsilon)}s\bigg|\bigg|,\end{eqnarray*} where $C_1^2:=\max\bigg\{\sup\limits_X|dw|^2_\omega,\,1\bigg\}$. Now, using the upper bound on $||s||^2$ in terms of $||\bar\partial_k^{\eta_\varepsilon + \lambda_\varepsilon}s||^2$ given by (\ref{eqn:map_spectra_injective_proof_2}), from this we infer the first inequality below: \begin{eqnarray}\label{eqn:map_spectra_injective_proof_3}\nonumber\bigg|\bigg|\Delta''_{k,\,(\varepsilon)}(\bar\partial_k^{\eta_\varepsilon + \lambda_\varepsilon}s)\bigg|\bigg| & \leq & \frac{CC_1}{\varepsilon^3\,k^{\frac{1}{b_2}}}\,\bigg(\frac{k}{\lambda_1\,\varepsilon^3} + 1\bigg)^{\frac{1}{2}}\,||\bar\partial_k^{\eta_\varepsilon + \lambda_\varepsilon}\,s|| + \bigg|\bigg|\bar\partial_k^{\eta_\varepsilon + \lambda_\varepsilon}\widetilde\Delta''_{k,\,(\varepsilon)}s\bigg|\bigg|\\
 & \leq & \frac{CC_1}{\varepsilon^3\,k^{\frac{1}{b_2}}}\,\bigg(\frac{2k}{\lambda_1\,\varepsilon^3}\bigg)^{\frac{1}{2}}\,||\bar\partial_k^{\eta_\varepsilon + \lambda_\varepsilon}\,s|| + \bigg|\bigg|\bar\partial_k^{\eta_\varepsilon + \lambda_\varepsilon}\widetilde\Delta''_{k,\,(\varepsilon)}s\bigg|\bigg|,\end{eqnarray} where the second inequality follows trivially from the first if we choose \begin{eqnarray}\label{eqn:1st-choice_lambda_1_epsilon}\lambda_1<k \hspace{3ex}\mbox{and}\hspace{3ex} 0<\varepsilon<1.\end{eqnarray}

Now, we estimate the last term in (\ref{eqn:map_spectra_injective_proof_3}) in the following way: \begin{eqnarray}\label{eqn:map_spectra_injective_proof_4}\nonumber\bigg|\bigg|\bar\partial_k^{\eta_\varepsilon + \lambda_\varepsilon}\widetilde\Delta''_{k,\,(\varepsilon)}s\bigg|\bigg|^2 & = & \langle\langle\widetilde\Delta''_{k,\,(\varepsilon)}\widetilde\Delta''_{k,\,(\varepsilon)}s,\,\widetilde\Delta''_{k,\,(\varepsilon)}s\rangle\rangle = \bigg|\bigg|\widetilde\Delta''_{k,\,(\varepsilon)}(\widetilde\Delta^{''\frac{1}{2}}_{k,\,(\varepsilon)} s)\bigg|\bigg|^2 \\
& \stackrel{(a)}{\leq} & \lambda_2^2\,\bigg|\bigg|\widetilde\Delta^{''\frac{1}{2}}_{k,\,(\varepsilon)} s\bigg|\bigg|^2 = \lambda_2^2\,\langle\langle\widetilde\Delta''_{k,\,(\varepsilon)}s,\,s\rangle\rangle = \lambda_2^2\,\bigg|\bigg|\bar\partial_k^{\eta_\varepsilon + \lambda_\varepsilon}s\bigg|\bigg|^2,\end{eqnarray} where the definition of the twisted Laplacian $\widetilde\Delta''_{k,\,(\varepsilon)}$ in bidegree $(n,\,0)$ and its self-adjointness were used and (a) followed from the assumption $s\in\widetilde{E}_I^{n,\,0}$ (which implies $\widetilde\Delta^{''\frac{1}{2}}_{k,\,(\varepsilon)}s\in\widetilde{E}_I^{n,\,0}$) and the fact that $I=(\lambda_1,\,\lambda_2]$.

Putting (\ref{eqn:map_spectra_injective_proof_3}) and (\ref{eqn:map_spectra_injective_proof_4}) together, we get the following upper estimate for the $L^2$-norm of $\Delta''_{k,\,(\varepsilon)}(\bar\partial_k^{\eta_\varepsilon + \lambda_\varepsilon}s)$:
 
\begin{eqnarray}\label{eqn:map_spectra_injective_proof_5}\bigg|\bigg|\Delta''_{k,\,(\varepsilon)}(\bar\partial_k^{\eta_\varepsilon + \lambda_\varepsilon}s)\bigg|\bigg| \leq \bigg(\lambda_2 + \frac{CC_1}{\varepsilon^3\,k^{\frac{1}{b_2}}}\,\bigg(\frac{2k}{\lambda_1\,\varepsilon^3}\bigg)^{\frac{1}{2}}\bigg)\,||\bar\partial_k^{\eta_\varepsilon + \lambda_\varepsilon}\,s||.\end{eqnarray}

Comparing (\ref{eqn:map_spectra_injective_proof_5}) with the lower estimate of the same quantity obtained in (\ref{eqn:map_spectra_injective_proof_1}) and using the fact that $\bar\partial_k^{\eta_\varepsilon + \lambda_\varepsilon}\,s\neq 0$ (since otherwise $s$ would lie in the kernel of $\widetilde\Delta''_{k,\,(\varepsilon)}$ and this would contradict the choice $s\in\widetilde{E}_I^{n,\,0}\setminus\{0\}$ with $I=(\lambda_1,\,\lambda_2]$ and $\lambda_1>0$), we get: \begin{eqnarray}\label{eqn:2nd-choice_lambda_1_epsilon}\varepsilon_0'\leq\frac{CC_1}{\varepsilon^3\,k^{\frac{1}{b_2}}}\,\bigg(\frac{2k}{\lambda_1\,\varepsilon^3}\bigg)^{\frac{1}{2}} = \frac{\sqrt{2}\,CC_1}{\sqrt{\lambda_1}\,\varepsilon^{4+\frac{1}{2}}\,k^{\frac{1}{b_2} - \frac{1}{2}}}\end{eqnarray} for all $k$ and $\varepsilon$ satisfying the preliminary choice (\ref{eqn:1st-choice_lambda_1_epsilon}). 

Now, if we fix an arbitrary $\delta\in(0,\,2/b_2)$ and choose $$\lambda_1 = \lambda_1(k):=\frac{2(CC_1)^2}{k^{1+\delta}} \hspace{3ex}\mbox{and}\hspace{3ex} \varepsilon_0'=\varepsilon_0'(k):=\varepsilon_0''\,k$$ for a constant $\varepsilon_0''>0$ independent of $k$, (\ref{eqn:2nd-choice_lambda_1_epsilon}) translates to \begin{eqnarray}\label{eqn:epsilon''_k-fraction}0<\varepsilon_0''\leq\frac{1}{\varepsilon^{\frac{9}{2}}\,k^{\frac{1}{2}\,(\frac{2}{b_2} - \delta)}}.\end{eqnarray} 

This is impossible if we fix $k$ and choose $\varepsilon$ such that $\varepsilon^{\frac{9}{2}}>\frac{1}{\varepsilon_0''\,k^{\frac{1}{2}\,(\frac{2}{b_2} - \delta)}}$.

Alternatively, we can choose \begin{eqnarray}\label{eqn:choice_epsilon-k}\varepsilon=\frac{1}{k^a} \hspace{3ex}\mbox{with}\hspace{1ex} 0<a<\frac{1}{9}\,\bigg(\frac{2}{b_2} - \delta\bigg),\end{eqnarray} in which case (\ref{eqn:epsilon''_k-fraction}) further translates, for all $k$, to \begin{eqnarray*}0<\varepsilon_0''\leq\frac{1}{k^{\frac{1}{2}\,(\frac{2}{b_2} - \delta - 9a)}},\end{eqnarray*} which is impossible if $k\gg 1$ since the right-hand-side term converges to $0$ when $k\to +\infty$, while the left-hand-side term is a positive constant. 

In both cases, we get the sought-after contradiction proving that the map (\ref{eqn:map_spectra_injective}) is injective with the above choices when $k$ is large enough.

The conclusion of this discussion is summed up in the following

\begin{Prop}\label{Prop:spectra_n0-n1} Let $X$ be a complex manifold of dimension $n$ endowed with a {\bf complete} K\"ahler metric $\omega$ such that $b_2$ is {\bf finite}. Suppose $\alpha$ is a ${\cal C}^\infty$ real $(1,\,1)$-form on $X$ such that $d\alpha=0$ and $\alpha$ satisfies Assumption A for a choice of approximately holomorphic ${\cal C}^\infty$ complex Hermitian line bundles $(L_k,\,h_k,\,D_k)\to X$ approximating the classes $\{k\alpha\}$ with the properties $(\ref{eqn:approximation_1})-(\ref{eqn:approximation_3})$. 

Suppose that \begin{eqnarray}\label{eqn:pos-hyp_spectra_n0-n1}i\Theta(D_k)^{1,\,1}\geq k\,\omega + idw\wedge d\overline{w}\end{eqnarray} on $X$ for all $k\in\Sigma$ with $k\geq 2$. Fix constants $0<\delta<\frac{2}{b_2}$, $0<a<\frac{1}{9}\,\bigg(\frac{2}{b_2} - \delta\bigg)$, $A>e$, $\varepsilon_0\in(0,\,1)$ and $c_0>0$. Further suppose that 

\vspace{1ex}

(i)\, either $k\in\Sigma$ and $\varepsilon$ are such that \begin{eqnarray}\label{eqn:spectra_n0-n1_hyp1}\frac{c_0}{k^{\frac{1}{9}\,(\frac{2}{b_2}-\delta)}}<\varepsilon<\delta_0:=\sqrt{(A/e)-1};\end{eqnarray}

\vspace{1ex}

(ii)\, or $\varepsilon=k^{-a}$ and $k$ is large enough.

\vspace{1ex}

  Then, the spectra of the twisted Laplacians \begin{eqnarray*}\widetilde\Delta''_{k,\,(\eta_\varepsilon+\lambda_\varepsilon)} &= & \widetilde\Delta''_{k,\,(\varepsilon)}:=(\bar\partial_k^{\eta_\varepsilon+\lambda_\varepsilon})^{\star}(\bar\partial_k^{\eta_\varepsilon+\lambda_\varepsilon}):{\cal C}^\infty_{n,\,0}(X,\,L_k)\longrightarrow {\cal C}^\infty_{n,\,0}(X,\,L_k),\\
\Delta''_{k,\,(\eta_\varepsilon+\lambda_\varepsilon)} &= & \Delta''_{k,\,(\varepsilon)}:= (\bar\partial_k^{\eta_\varepsilon+\lambda_\varepsilon})(\bar\partial_k^{\eta_\varepsilon+\lambda_\varepsilon})^{\star} + (^{\eta_\varepsilon+\lambda_\varepsilon}\bar\partial_k)^{\star}(^{\eta_\varepsilon+\lambda_\varepsilon}\bar\partial_k):{\cal C}^\infty_{n,\,1}(X,\,L_k)\longrightarrow {\cal C}^\infty_{n,\,1}(X,\,L_k),\end{eqnarray*} defined using the sum of the auxiliary functions introduced in (\ref{eqn:eta-lambda-eps-def}), satisfy the following properties: \begin{eqnarray}\label{eqn:spectrum_n0}\mbox{Spec}^{n,\,0}(\widetilde\Delta''_{k,\,(\varepsilon)})\bigcap\bigg(\frac{c}{k^{1+\delta}},\,(1-\varepsilon_0)k\bigg] = \emptyset\end{eqnarray} and \begin{eqnarray}\label{eqn:spectrum_n1}\mbox{Spec}^{n,\,1}(\Delta''_{k,\,(\varepsilon)})\subset\bigg[k,\,+\infty\bigg),\end{eqnarray} where $c:=C\,\max\bigg\{\sup\limits_X|dw|^2_\omega,\,1\bigg\}$ for some constant $C>0$ independent of $k$, while $k_{\varepsilon_0}$ is a positive integer depending only on $\varepsilon_0$.

\end{Prop}

\noindent {\it Proof.} Only the proof of (\ref{eqn:spectrum_n0}) needs completing, as (\ref{eqn:spectrum_n1}) is an immediate consequence of (i) of Proposition \ref{Prop:BKN-consequence} in the case $q=1$. 

Now, in the notation of this $\S$\ref{subsection:spectra_n0-n1}, if we fix $\delta$ as in the statement and arbitrary $0<\varepsilon_0''<\varepsilon_0<1$ and then choose $$\lambda_2 = \lambda_2(k):=(1-\varepsilon_0)\,k,$$ we get $J= [0,\,\lambda_2 + \varepsilon_0'] =[0,\,\lambda_2 + \varepsilon_0''\,k] = [0,\,(1-\varepsilon_0+\varepsilon_0'')\,k]\subset[0,\,k)$. Thus, (\ref{eqn:spectrum_n1}) implies that $E_J^{n,\,1}=\{0\}$ for all $k\in\Sigma$ with $k\geq 2$. 

Hence, thanks to the injectivity of the map (\ref{eqn:map_spectra_injective}) proved above for $k$ large enough, we infer that $\widetilde{E}_I^{n,\,0}=\{0\}$ for $k$ large enough, where $$I=(\lambda_1,\,\lambda_2] = \bigg(\frac{2(CC_1)^2}{k^{1+\delta}},\, (1-\varepsilon_0)\,k \bigg].$$ This is equivalent to (\ref{eqn:spectrum_n0}) if we set $c:=2(CC_1)^2$, relabel $2C^2$ as $C$ and recall the definition of $C_1^2$ in the injectivity proof above the statement.  \hfill $\Box$

\subsection{The main a priori twisted $L^2$-estimate}\label{subsection:main-a-priori-estimate}  We are now in a position to derive the main ingredient for the proof of our non-integrable $L^2$ extension theorem. This is the non-integrable analogue of the a priori $L^2$ estimates for the integrable case proved in [Ohs95] (see also [Dem01, Proposition 12.4]). 

\begin{The}\label{The:main-a-priori-estimate} Let $(X,\,\omega)$ be a {\bf complete K\"ahler} manifold with $\mbox{dim}_{\C}X=n$ and $b_2(X)<+\infty$. Let $\alpha$ be a ${\cal C}^\infty$ real $(1,\,1)$-form on $X$ such that $d\alpha=0$ and $\alpha$ satisfies {\bf Assumption A}. Choose any sequence $(L_k,\,h_k,\,D_k=\partial_k+\bar\partial_k)\longrightarrow X$ ($k\in\Sigma\subset\N^{\star}$) of approximately holomorphic ${\cal C}^\infty$ Hermitian line bundles associated with $\alpha$ such that $$\frac{i}{2\pi}\,\Theta(D_k)=\alpha_k \hspace{2ex} \mbox{and} \hspace{2ex} ||\alpha_k-k\alpha||\leq_{{\cal C}^\infty}\frac{C}{k^{1/b_2}}.$$

 We make the following {\bf positivity assumption} on $\alpha$: \begin{eqnarray}\label{eqn:main-est_curvature-assumption}2\pi\alpha > \omega + idw\wedge d\overline{w} \hspace{3ex} \mbox{on} \hspace{1ex} X.\end{eqnarray}

Fix constants $0<\delta<\frac{2}{b_2}$, $0<a<\frac{1}{9}\,\bigg(\frac{2}{b_2} - \delta\bigg)$, $A>e$ and $\varepsilon_0\in(0,\,1)$. Consider the auxiliary functions $\eta_\varepsilon$ and $\lambda_\varepsilon$ introduced in (\ref{eqn:eta-lambda-eps-def}) and associated with a bounded holomorphic function $w:X\longrightarrow\C$ such that $dw(x)\neq 0$ whenever $w(x)=0$. Let $c:=C\,\max\bigg\{\sup\limits_X|dw|^2_\omega,\,1\bigg\}$ for some constant $C>0$ independent of $k$.

\vspace{1ex}

(i)\, Suppose that:

\vspace{1ex}

\hspace{6ex} $\bullet$ either $k\in\Sigma$ and $\varepsilon>0$ satisfy (\ref{eqn:spectra_n0-n1_hyp1}), 

\vspace{1ex}

\hspace{6ex} $\bullet$ or $\varepsilon:=k^{-a}$ and $k$ is large enough.

\vspace{1ex}

\noindent Then, for every $s\in L^2_{n,\,0}(X,\,L_k)$, letting $$s = s_{h,\,\varepsilon} + s_{nh,\,\varepsilon}$$ with $s_{h,\,\varepsilon}\in{\cal H}^{n,\,0}_{[0,\,c/k^{1+\delta}],\,\widetilde\Delta''_{k,\,(\varepsilon)}}(X,\,L_k) = \mbox{Im}\,\bigg(\chi_{[0,\,c/k^{1+\delta}]}(\widetilde\Delta''_{k,\,(\varepsilon)})\bigg):={\cal H}^{n,\,0}_{k,\,\varepsilon}(X,\,L_k)\subset L^2_{n,\,0}(X,\,L_k)$ and

\hspace{0.5ex} $s_{nh,\,\varepsilon}\in{\cal H}^{n,\,0}_{((1-\varepsilon_0)\,k,\,+\infty),\,\widetilde\Delta''_{k,\,(\varepsilon)}}(X,\,L_k) = \mbox{Im}\,\bigg(\chi_{((1-\varepsilon_0)\,k,\,+\infty)}(\widetilde\Delta''_{k,\,(\varepsilon)})\bigg):={\cal N}^{n,\,0}_{k,\,\varepsilon}(X,\,L_k)$, 

\noindent be the decomposition of $s$ w.r.t. the $L^2_\omega$-orthogonal splitting \begin{equation}\label{eqn:n0_2-space-decomp_epsilon-k}L^2_{n,\,0}(X,\,L_k) = {\cal H}^{n,\,0}_{k,\,\varepsilon}(X,\,L_k) \oplus {\cal N}^{n,\,0}_{k,\,\varepsilon}(X,\,L_k),\end{equation}

\noindent the following estimate holds: \begin{equation}\label{eqn:main-a-priori-estimate}||s_{nh,\,\varepsilon}||^2\leq \frac{1}{(1-\varepsilon_0)\,k}\,||\bar\partial_k^{\eta_\varepsilon + \lambda_\varepsilon}s||^2.\end{equation} 

\vspace{1ex}

(ii)\, Suppose that $k\in\Sigma$ and $\varepsilon>0$ are arbitrary. Then, for every $s\in L^2_{n,\,0}(X,\,L_k)$, letting $$s = s_{h,\,\varepsilon} + t_\varepsilon + s_{nh,\,\varepsilon}$$ with $t_\varepsilon\in{\cal H}^{n,\,0}_{(c/k^{1+\delta},\,(1-\varepsilon_0)\,k],\,\widetilde\Delta''_{k,\,(\varepsilon)}}(X,\,L_k) = \mbox{Im}\,\bigg(\chi_{(c/k^{1+\delta},\,(1-\varepsilon_0)\,k]}(\widetilde\Delta''_{k,\,(\varepsilon)})\bigg):={\cal M}^{n,\,0}_{k,\,\varepsilon}(X,\,L_k)\subset L^2_{n,\,0}(X,\,L_k)$ and, as in (i), $s_{h,\,\varepsilon}\in{\cal H}^{n,\,0}_{k,\,\varepsilon}(X,\,L_k)$ and $s_{nh,\,\varepsilon}\in{\cal N}^{n,\,0}_{k,\,\varepsilon}(X,\,L_k)$, be the decomposition of $s$ w.r.t. the $L^2_\omega$-orthogonal splitting \begin{equation}\label{eqn:n0_3-space-decomp_epsilon-k}L^2_{n,\,0}(X,\,L_k) = {\cal H}^{n,\,0}_{k,\,\varepsilon}(X,\,L_k) \oplus {\cal M}^{n,\,0}_{k,\,\varepsilon}(X,\,L_k) \oplus {\cal N}^{n,\,0}_{k,\,\varepsilon}(X,\,L_k),\end{equation} 

\noindent the following estimates hold: \begin{equation}\label{eqn:main-a-priori-estimate_rough}||t_\varepsilon||^2\leq \frac{k^{1+\delta}}{c}\,||\bar\partial_k^{\eta_\varepsilon + \lambda_\varepsilon}s||^2 \hspace{3ex} \mbox{and} \hspace{3ex} ||s_{nh,\,\varepsilon}||^2\leq \frac{1}{(1-\varepsilon_0)\,k}\,||\bar\partial_k^{\eta_\varepsilon + \lambda_\varepsilon}s||^2.\end{equation}

\end{The}

\noindent {\it Proof.} (i)\, After multiplication by $k$, the positivity assumption (\ref{eqn:main-est_curvature-assumption}) implies hypothesis (\ref{eqn:pos-hyp_spectra_n0-n1}), so Proposition \ref{Prop:spectra_n0-n1} applies and ensures the validity of (\ref{eqn:spectrum_n0}). In particular, the two-space decomposition (\ref{eqn:n0_2-space-decomp_epsilon-k}) holds.

Fix $s\in L^2_{n,\,0}(X,\,L_k)$. On the one hand, from $s_{nh,\,\varepsilon}\in{\cal N}^{n,\,0}_{k,\,\varepsilon}(X,\,L_k)$ we get the inequality below: \begin{eqnarray}\label{eqn:main-est_curvature-assumption_proof_1}||\bar\partial_k^{\eta_\varepsilon + \lambda_\varepsilon}s_{nh,\,\varepsilon}||^2 = \langle\langle\widetilde\Delta''_{k,\,(\varepsilon)}s_{nh,\,\varepsilon},\,s_{nh,\,\varepsilon}\rangle\rangle\geq (1-\varepsilon_0)k\,||s_{nh,\,\varepsilon}||^2,\end{eqnarray} where the equality follows from the definition $\widetilde\Delta''_{k,\,(\varepsilon)}:=(\bar\partial_k^{\eta_\varepsilon+\lambda_\varepsilon})^{\star}(\bar\partial_k^{\eta_\varepsilon+\lambda_\varepsilon})$ of the twisted Laplacian in bidegree $(n,\,0)$.  

On the other hand, we note that $\bar\partial_k^{\eta_\varepsilon + \lambda_\varepsilon}s_{h,\,\varepsilon}$ and $\bar\partial_k^{\eta_\varepsilon + \lambda_\varepsilon}s_{nh,\,\varepsilon}$ are orthogonal to each other since \begin{eqnarray*}\langle\langle\bar\partial_k^{\eta_\varepsilon + \lambda_\varepsilon}s_{h,\,\varepsilon},\,\bar\partial_k^{\eta_\varepsilon + \lambda_\varepsilon}s_{nh,\,\varepsilon}\rangle\rangle = \langle\langle\widetilde\Delta''_{k,\,(\varepsilon)}s_{h,\,\varepsilon},\,s_{nh,\,\varepsilon}\rangle\rangle = 0,\end{eqnarray*} where the last equality follows from $\widetilde\Delta''_{k,\,(\varepsilon)}s_{h,\,\varepsilon}\in{\cal H}^{n,\,0}_{k,\,\varepsilon}(X,\,L_k)$ (which is due to the fact that $s_{h,\,\varepsilon}\in{\cal H}^{n,\,0}_{k,\,\varepsilon}(X,\,L_k)$), from $s_{nh,\,\varepsilon}\in{\cal N}^{n,\,0}_{k,\,\varepsilon}(X,\,L_k)$ and from the subspaces ${\cal H}^{n,\,0}_{k,\,\varepsilon}(X,\,L_k)$ and ${\cal N}^{n,\,0}_{k,\,\varepsilon}(X,\,L_k)$ being $L^2_\omega$-orthogonal to each other. Therefore, we get: \begin{eqnarray}\label{eqn:main-est_curvature-assumption_proof_2}||\bar\partial_k^{\eta_\varepsilon + \lambda_\varepsilon}s||^2 = ||\bar\partial_k^{\eta_\varepsilon + \lambda_\varepsilon}s_{h,\,\varepsilon}||^2 + ||\bar\partial_k^{\eta_\varepsilon + \lambda_\varepsilon}s_{nh,\,\varepsilon}||^2\geq||\bar\partial_k^{\eta_\varepsilon + \lambda_\varepsilon}s_{nh,\,\varepsilon}||^2.\end{eqnarray}

Thus, putting (\ref{eqn:main-est_curvature-assumption_proof_1}) and (\ref{eqn:main-est_curvature-assumption_proof_2}) together, we get: \begin{eqnarray*}||s_{nh,\,\varepsilon}||^2\leq \frac{1}{(1-\varepsilon_0)\,k}\,||\bar\partial_k^{\eta_\varepsilon + \lambda_\varepsilon}s_{nh,\,\varepsilon}||^2\leq \frac{1}{(1-\varepsilon_0)\,k}\,||\bar\partial_k^{\eta_\varepsilon + \lambda_\varepsilon}s||^2,\end{eqnarray*} which proves (\ref{eqn:main-a-priori-estimate}).  

\vspace{2ex}

(ii)\, When no restriction is imposed on $k$ and $\varepsilon$, the spectral gap equality (\ref{eqn:spectrum_n0}) need not hold, so the middle term ${\cal M}^{n,\,0}_{k,\,\varepsilon}(X,\,L_k)$ of decomposition (\ref{eqn:n0_3-space-decomp_epsilon-k}) need not be non-trivial. 

Estimates (\ref{eqn:main-a-priori-estimate_rough}) are obtained in the same way as (\ref{eqn:main-a-priori-estimate}). For example, in the case of $t_\varepsilon$, we have: \begin{eqnarray*}||\bar\partial_k^{\eta_\varepsilon + \lambda_\varepsilon}s||^2 & = & ||\bar\partial_k^{\eta_\varepsilon + \lambda_\varepsilon}s_{h,\,\varepsilon}||^2 + ||\bar\partial_k^{\eta_\varepsilon + \lambda_\varepsilon}t_\varepsilon||^2 + ||\bar\partial_k^{\eta_\varepsilon + \lambda_\varepsilon}s_{nh,\,\varepsilon}||^2 \\
& \geq & ||\bar\partial_k^{\eta_\varepsilon + \lambda_\varepsilon}t_\varepsilon||^2 = \langle\langle\widetilde\Delta''_{k,\,(\varepsilon)}t_\varepsilon,\,t_\varepsilon\rangle\rangle\geq \frac{c}{k^{1+\delta}}\,||t_\varepsilon||^2,\end{eqnarray*} where the equality on the first line follows from the mutual $L^2_\omega$-orthogonality of $\bar\partial_k^{\eta_\varepsilon + \lambda_\varepsilon}s_{h,\,\varepsilon}$, $\bar\partial_k^{\eta_\varepsilon + \lambda_\varepsilon}t_\varepsilon$ and $\bar\partial_k^{\eta_\varepsilon + \lambda_\varepsilon}s_{nh,\,\varepsilon}$ (same argument as in (i)), while the last inequality follows from $t_\varepsilon$ lying in ${\cal M}^{n,\,0}_{k,\,\varepsilon}(X,\,L_k)$ and from the definition of this space.  \hfill $\Box$

\vspace{2ex}

\begin{Def}\label{Def:constrained} In the setting of Theorem \ref{The:main-a-priori-estimate}, if $k\in\Sigma$ and $\varepsilon>0$ satisfy any of the two hypotheses imposed on them under (i), we say that $\varepsilon$ is {\bf constrained by $k$}. Otherwise, we say that $k$ and $\varepsilon$ are {\bf mutually independent}. 

\end{Def}

\vspace{2ex}

As a complement to Theorem \ref{The:main-a-priori-estimate}, we notice an immediate consequence of the definition of the space ${\cal H}^{n,\,0}_{k,\,\varepsilon}(X,\,L_k)$, itself forced on us by the spectral gap formula (\ref{eqn:spectrum_n0}) in bidegree $(n,\,0)$.

\begin{Obs}\label{Obs:dbar_s_h-epsilon_estimate} Hypotheses and notation as in (ii) of Theorem \ref{The:main-a-priori-estimate}, then every $s_{h,\,\varepsilon}\in{\cal H}^{n,\,0}_{k,\,\varepsilon}(X,\,L_k)$ satisfies the estimate: \begin{eqnarray}\label{eqn:dbar_s_h-epsilon_estimate}||\bar\partial_k^{\eta_\varepsilon + \lambda_\varepsilon}s_{h,\,\varepsilon}||^2\leq\frac{c}{k^{1+\delta}}\,||s_{h,\,\varepsilon}||^2.\end{eqnarray}

\end{Obs}

\noindent {\it Proof.} We have: \begin{eqnarray*}||\bar\partial_k^{\eta_\varepsilon + \lambda_\varepsilon}s_{h,\,\varepsilon}||^2 = \langle\langle\widetilde\Delta''_{k,\,(\varepsilon)}s_{h,\,\varepsilon},\,s_{h,\,\varepsilon}\rangle\rangle\leq\frac{c}{k^{1+\delta}}\,||s_{h,\,\varepsilon}||^2,\end{eqnarray*} where the inequality follows from definition of ${\cal H}^{n,\,0}_{k,\,\varepsilon}(X,\,L_k)$ given in Theorem \ref{The:main-a-priori-estimate}. \hfill $\Box$

\vspace{3ex}

We will later need a variant of our above $L^2$-estimates in higher-order Sobolev spaces $W^l$ ($l\geq 1$), so we now briefly discuss a definition of Sobolev inner products that will suit our purposes.

\vspace{1ex}

$\bullet$ To begin with, suppose that $X$ is a compact $n$-dimensional complex manifold on which a Hermitian metric $\omega$ has been fixed. Let $\Delta'':C^\infty_{p,\,q}(X,\,\C)\longrightarrow C^\infty_{p,\,q}(X,\,\C)$ be the induced $\bar\partial$-Laplacian defined on the smooth $(p,\,q)$-forms on $X$ of any bidegree $(p,\,q)$ in the usual way as $\Delta'' = \bar\partial\bar\partial^\star + \bar\partial^\star\bar\partial$, where $\bar\partial^\star$ is the formal adjoint of $\bar\partial$ with respect to the $L^2$-inner product $\langle\langle\,\,\cdot\,\,,\,\,\cdot\,\,\rangle\rangle$ induced by $\omega$.

Recall the classical {\it G\r{a}rding inequality} satisfied by any elliptic non-negative self-adjoint differential operator, in particular by $\Delta^{''l}$ which is of order $2l$ on $X$: there exist constants $C_1(l),\,C_2(l)>0$ such that for every $(p,\,q)$-form $u$ on $X$, we have: \begin{eqnarray}\label{eqn:Garding_scalar}C_2(l)\,||u||^2_{W^l}\leq\langle\langle\Delta^{''l} u,\,u\rangle\rangle + C_1(l)\,||u||^2,\end{eqnarray} where $||u||_{W^l}^2$ is the standard squared $l$-th Sobolev norm of $u$ defined as the sum of the squared $L^2$-norms of all the derivatives of $u$ up to order $l$. 

For every positive integer $l$, we define the following inner product on $C^\infty_{p,\,q}(X,\,\C)$: \begin{eqnarray}\label{eqn:Sobolev-product_Laplacian}\langle\langle u,\,v\rangle\rangle_{\Delta^{''l}}:=\bigg\langle\bigg\langle\bigg(\Delta^{''l} + C_1(l)\,\mbox{Id}\bigg)u,\,v\bigg\rangle\bigg\rangle, \hspace{5ex} u,v\in C^\infty_{p,\,q}(X,\,\C),\end{eqnarray} where $\mbox{Id}$ is the identity operator.

\vspace{1ex}

The first observation is that the formal adjoint of $\bar\partial$ with respect to this inner product coincides with the standard formal adjoint.

\begin{Lem}\label{Lem:Sobolev-product_adjoints-equal} If $\bar\partial^\star_{\Delta^{''l}}:C^\infty_{p,\,q+1}(X,\,\C)\longrightarrow C^\infty_{p,\,q}(X,\,\C)$ is the formal adjoint of $\bar\partial:C^\infty_{p,\,q}(X,\,\C)\longrightarrow C^\infty_{p,\,q+1}(X,\,\C)$ with respect to the inner product $\langle\langle\,\,\cdot\,\,,\,\,\cdot\,\,\rangle\rangle_{\Delta^{''l}}$ and $\bar\partial^\star:C^\infty_{p,\,q+1}(X,\,\C)\longrightarrow C^\infty_{p,\,q}(X,\,\C)$ is the formal adjoint of $\bar\partial:C^\infty_{p,\,q}(X,\,\C)\longrightarrow C^\infty_{p,\,q+1}(X,\,\C)$ with respect to the $L^2$-inner product $\langle\langle\,\,\cdot\,\,,\,\,\cdot\,\,\rangle\rangle$, then \begin{eqnarray}\label{eqn:Sobolev-product_adjoints-equal}\bar\partial^\star_{\Delta^{''l}} = \bar\partial^\star.\end{eqnarray}

\end{Lem}

\noindent {\it Proof.} Let $u\in C^\infty_{p,\,q}(X,\,\C)$ and $v\in C^\infty_{p,\,q+1}(X,\,\C)$ be arbitrary.

On the one hand, by definition of $\bar\partial^\star_{\Delta^{''l}}$, we have: \begin{eqnarray}\label{eqn:Sobolev-product_adjoints-equal_proof_1}\langle\langle\bar\partial u,\,v\rangle\rangle_{\Delta^{''l}} = \langle\langle u,\,\bar\partial^\star_{\Delta^{''l}}v\rangle\rangle_{\Delta^{''l}}.\end{eqnarray}

On the other hand, by definition of the inner product introduced in (\ref{eqn:Sobolev-product_Laplacian}), we have: \begin{eqnarray}\label{eqn:Sobolev-product_adjoints-equal_proof_2}\nonumber\langle\langle\bar\partial u,\,v\rangle\rangle_{\Delta^{''l}} & = & \langle\langle(\Delta^{''l} + \mbox{Id})\,\bar\partial u,\,v\rangle\rangle = \langle\langle\bar\partial\Delta^{''l} u,\,v\rangle\rangle + \langle\langle\bar\partial u,\,v\rangle\rangle = \langle\langle(\Delta^{''l} + \mbox{Id})\,u,\,\bar\partial^\star v\rangle\rangle\\
  & = & \langle\langle u,\,\bar\partial^\star v\rangle\rangle_{\Delta^{''l}}.\end{eqnarray}

Putting (\ref{eqn:Sobolev-product_adjoints-equal_proof_1}) and (\ref{eqn:Sobolev-product_adjoints-equal_proof_2}) together, we get: \begin{eqnarray*}\langle\langle u,\,\bar\partial^\star_{\Delta^{''l}}v\rangle\rangle_{\Delta^{''l}} = \langle\langle u,\,\bar\partial^\star v\rangle\rangle_{\Delta^{''l}}.\end{eqnarray*} Since $u$ and $v$ are arbitrary, the conclusion follows. \hfill $\Box$

\vspace{1ex}

The next observation is that, due to the classical {\it G\r{a}rding inequality} (\ref{eqn:Garding_scalar}), the inner product defined in (\ref{eqn:Sobolev-product_Laplacian}) induces the same Sobolev space $W^l_{p,\,q}(X,\,\C)$ on the $(p,\,q)$-forms on $X$, equipped with the same topology, as the standard Sobolev inner product defined as the sum of the $L^2$-inner products of all the derivatives of $u$ and $v$ up to order $l$. Indeed, (\ref{eqn:Garding_scalar}) implies that the norms $||u||_{W^l}$ and $||u||_{\Delta^{''l}}$ are equivalent.

Together with Lemma \ref{Lem:Sobolev-product_adjoints-equal}, this further implies that the operators $\Delta''$ and $\bar\partial\bar\partial^\star_{\Delta^{''l}} + \bar\partial^\star_{\Delta^{''l}}\bar\partial$ coincide and that they have the same eigenvalues and the same eigenspaces as each other and independently of whether they are made to act on $L^2_{p,\,q}(X,\,\C)$ or on any other Sobolev space $W^l_{p,\,q}(X,\,\C)$ for any positive integer $l$. In particular, the eigenspaces (which are contained in $C^\infty_{p,\,q}(X,\,\C)$ by (hypo-)ellipticity) corresponding to any pair of distinct eigenvalues are mutually orthogonal for the $L^2$-inner product $\langle\langle\,\,\cdot\,\,,\,\,\cdot\,\,\rangle\rangle$ and for every inner product $\langle\langle\,\,\cdot\,\,,\,\,\cdot\,\,\rangle\rangle_{\Delta^{''l}}$ with $l$ a positive integer. 
  
\vspace{1ex}

$\bullet$ It is standard (see e.g. [Wel08, IV, $\S1$]) that the classical Sobolev inner products and the Sobolev spaces are defined on a manifold using local coordinates and local trivialisations followed by the use of a partition of unity to globalise the construction. In the same way, our modified Sobolev inner product introduced in (\ref{eqn:Sobolev-product_Laplacian}) for every $l\in\N^\star$ starting from the integrable operator $\bar\partial$ associated with the complex structure of the manifold $X$ and from the induced Laplacian $\Delta''$ can be generalised to the case where $\bar\partial$ is replaced by our (possibly non-integrable) operator $\bar\partial_k$ and $\Delta''$ is replaced by $\Delta''_k$ even when $X$ is non-compact, but only equipped with a complete metric $\omega$. Our results obtained above for the $L^2$-inner product remain valid for the modified Sobolev inner products induced by $\Delta^{''l}_k$ for every positive integer $l$. In particular, we have

\begin{Obs}\label{Obs:Sobolev_main-a-priori-estimate} Estimates (\ref{eqn:main-a-priori-estimate_rough}) also hold if the $L^2$-norm $||\,\,||$ is replaced by any modified Sobolev norm $||\,\,||_{\Delta^{''l}_k}$ with $l$ a non-negative integer. 

\end{Obs}

\section{Proof of the extension Theorem \ref{nonintOT}}\label{section:main-theorem} Let us fix an arbitrary $k\in\Sigma$. We will split the proof into several steps that broadly follow the ones in the classical Ohsawa-Takegoshi Theorem and will also point out the new ingredients needed in this non-integrable situation.

\vspace{1ex}

\noindent {\bf (A) Rough extension.} Let $(U_j)_{j\in J}$ be a covering of $Y$ by open coordinate patches $U_j\subset X$ such that $z_1^{(j)},\dots ,z_{n-1}^{(j)}, w$ are local holomorphic coordinates on $U_j$ and $Y\cap U_j =\{w=0\}$. Thus, $dz_1^{(j)}\wedge\dots\wedge dz_{n-1}^{(j)}$ is a non-vanishing holomorphic section of $K_Y$ on $Y\cap U_j$ and $dz_1^{(j)}\wedge\dots\wedge dz_{n-1}^{(j)}\wedge dw$ is a non-vanishing holomorphic section of $K_X$ on $U_j$. We may assume that the restriction of $L_k$ to each $U_j$ is the trivial $\C$-line bundle and let $e_k^{(j)}$ be a ${\cal C}^\infty$ frame of $L_{k|Y\cap U_j}$. We also consider a partition of unity $(\theta_j)_{j\in J}$ associated with $(U_j)_{j\in J}$. Thus, $\mbox{Supp}\,\theta_j\subset U_j$ and $$\sum\limits_{j\in J}\theta_j = 1, \hspace{3ex} \mbox{hence} \hspace{3ex} \sum\limits_{j\in J}\bar\partial\theta_j = 0,$$ on a neighbourhood of $Y$. Finally, for every $j\in J$, let $\pi_j:U_j\rightarrow Y\cap U_j$ be the holomorphic projection $(z_1^{(j)},\dots ,z_{n-1}^{(j)}, w)\mapsto(z_1^{(j)},\dots ,z_{n-1}^{(j)})$.

Let $f=f_k\in C^\infty_{0,\,0}(Y,\,K_Y\otimes L_{k|Y})$ such that $\int\limits_Y|f|^2_{h_k}\,dV_{Y,\,\omega}<+\infty$. In particular, its restriction to $Y\cap U_j$ is of the shape

$$f_{|Y\cap U_j} = f_j\,dz_1^{(j)}\wedge\dots\wedge dz_{n-1}^{(j)}\otimes e_k^{(j)},$$

\noindent for some ${\cal C}^\infty$ function $f_j$ on $Y\cap U_j$. We put

$$\widetilde{f}:=\sum\limits_{j\in J}\theta_j\,(f_j\circ\pi_j)\,dz_1^{(j)}\wedge\dots\wedge dz_{n-1}^{(j)}\wedge dw\otimes e_k^{(j)}\in {\cal C}^\infty_{n,\,0}(X,\,L_k).$$

\noindent Thus, each $\widetilde{f}_j:=f_j\circ\pi_j$ is a smooth function on $U_j$ such that $\widetilde{f}_{j|Y\cap U_j}=f_j$. Consequently, $\widetilde{f}_{|Y} = f\wedge dw$ on $Y$.

 We now compute $(\bar\partial_k\widetilde{f})_{|Y}$. We get

$$\bar\partial_k\widetilde{f} = \sum\limits_{j\in J}(\bar\partial\theta_j)\wedge(\widetilde{f}_j\,dz_1^{(j)}\wedge\dots\wedge dz_{n-1}^{(j)}\wedge dw\otimes e_k^{(j)}) + \sum\limits_{j\in J}\theta_j\,\bar\partial_k(\widetilde{f}_j\,dz_1^{(j)}\wedge\dots\wedge dz_{n-1}^{(j)}\wedge dw\otimes e_k^{(j)})$$

\noindent on $X$, hence by taking restrictions to $Y$, we get:

\begin{eqnarray}\nonumber (\bar\partial_k\widetilde{f})_{|Y} & = & \sum\limits_{j\in J}(\bar\partial\theta_j)_{|Y\cap U_j}\wedge(\widetilde{f}_j\,dz_1^{(j)}\wedge\dots\wedge dz_{n-1}^{(j)}\wedge dw\otimes e_k^{(j)})_{|Y\cap U_j} \\
\nonumber & + & \sum\limits_{j\in J}\theta_{j|Y}\,\bar\partial_k(\widetilde{f}_j\,dz_1^{(j)}\wedge\dots\wedge dz_{n-1}^{(j)}\wedge dw\otimes e_k^{(j)})_{|Y\cap U_j} \\
\nonumber & = & (\sum\limits_{j\in J}\bar\partial\theta_j)\wedge f + (\sum\limits_{j\in J}\theta_j)_{|Y}\,\bar\partial_k(f\wedge dw) = \bar\partial_k(f\wedge dw)  \end{eqnarray}

\noindent on $Y$.

\begin{Conc}\label{Conc:rough-extension} For any $f=f_k\in C^\infty_{0,\,0}(Y,\,K_Y\otimes L_{k|Y})$ such that $\int_Y|f|^2_{h_k}\,dV_{Y,\,\omega}<+\infty$, there exists a rough extension $\widetilde{f}=\widetilde{f}_k\in {\cal C}^\infty_{n,\,0}(X,\,L_k)$ satisfying the following identities on $Y$:

$$(i)\,\,\widetilde{f}_{|Y}=f\wedge dw  \hspace{3ex} \mbox{and} \hspace{3ex} (ii)\,\,  (\bar\partial_k\widetilde{f})_{|Y} = \bar\partial_k(f\wedge dw) = \bar\partial_kf\wedge dw.$$

\end{Conc}

\noindent {\bf (B) A cutoff function.} As is standard in such a context, we let $\rho:[0,\,+\infty)\longrightarrow[0,\,1]$ be a ${\cal C}^\infty$ function such that $$\rho(x)=1 \hspace{2ex} \mbox{for all} \hspace{2ex} x\in\bigg[0,\,\frac{1}{3}\bigg], \hspace{3ex} \rho(x)=0 \hspace{2ex} \mbox{for all} \hspace{2ex} x\in[1,\,+\infty)  \hspace{2ex} \mbox{and} \hspace{2ex} \sup\limits_{x\in[0,\,1]}|\rho'(x)|\leq 2.$$ 

\noindent Again as usual, for every constant $\varepsilon>0$, we put $$\rho_\varepsilon(w):=\rho\bigg(\frac{|w|^2}{\varepsilon^2}\bigg), \hspace{2ex} \mbox{hence} \hspace{2ex} \rho_\varepsilon(w) = 1 \hspace{1ex} \mbox{on}\hspace{1ex} \{|w|\leq\varepsilon/\sqrt{3}\} \hspace{2ex} \mbox{and} \hspace{2ex} \mbox{Supp}\,(\rho_\varepsilon)\subset\{|w|\leq\varepsilon\}.$$ 

\noindent The size of the tubular neighbourhood $\{|w|\leq\varepsilon\}$ of $Y$ will be shrunk by letting $\varepsilon$ converge to zero later on. Since the rough extension $\widetilde{f}$ is uncontrallable away from $Y$, we cut it off using the function $\rho_\varepsilon$, so we put

$$\widetilde{f}_\varepsilon=\widetilde{f}_{k,\,\varepsilon}:= \rho_\varepsilon(w)\,\widetilde{f}\in {\cal C}^\infty_{n,\,0}(X,\,L_k).$$

\begin{Lem}\label{Lem:f-tilda_epsilon} Besides the obvious properties $$\widetilde{f}_\varepsilon = \widetilde{f} \hspace{2ex} \mbox{on} \hspace{1ex} \{|w|\leq\varepsilon/\sqrt{3}\}  \hspace{3ex} \mbox{and} \hspace{3ex} \widetilde{f}_\varepsilon = 0 \hspace{2ex} \mbox{on} \hspace{1ex} X\setminus\{|w|\leq\varepsilon\},$$

\noindent the section $\widetilde{f}_\varepsilon$ also satisfies $$\bar\partial_k\widetilde{f}_\varepsilon=\rho'\bigg(\frac{|w|^2}{\varepsilon^2}\bigg)\,\frac{w}{\varepsilon^2}\,d\overline{w}\wedge\widetilde{f} + \rho_\varepsilon(w)\,\bar\partial_k\widetilde{f}\in {\cal C}^\infty_{n,\,1}(X,\,L_k).$$

\noindent In particular, $\bar\partial_k\widetilde{f}_\varepsilon = \bar\partial_k\widetilde{f}$ on $\{|w|\leq\varepsilon/\sqrt{3}\}$.

\end{Lem}

\noindent {\it Proof.} Immediate computation.  \hfill $\Box$

\begin{Lem}\label{Lem:v_epsilon} For every sufficiently small $\varepsilon>0$, the following section: $$v_\varepsilon:=\frac{1}{w-2\varepsilon}\,\bar\partial_k\widetilde{f}_\varepsilon = \frac{1}{\varepsilon^2}\,\frac{w}{w-2\varepsilon} \,\rho'\bigg(\frac{|w|^2}{\varepsilon^2}\bigg)\,d\overline{w}\wedge\widetilde{f} + \frac{\rho_\varepsilon(w)}{w-2\varepsilon}\,\bar\partial_k\widetilde{f}\in L^{1,\,loc}_{n,\,1}(X,\,L_k)$$

\noindent is ${\cal C}^\infty$ on $X\setminus Y_\varepsilon$ (where $Y_\varepsilon:=\{w=2\varepsilon\}$) and $L^{2-\delta}_{loc}$ on $X$ for all $0<\delta<1$, while having the following extra properties: $$\bar\partial_k v_\varepsilon = -\frac{2\pi i}{w-2\varepsilon}\,\alpha_k^{0,\,2}\wedge\widetilde{f}_\varepsilon \hspace{6ex} \mbox{on}\hspace{1ex} X\setminus Y_\varepsilon$$

\noindent and $$\bar\partial_k v_\varepsilon = -\frac{2\pi i}{w-2\varepsilon}\,\alpha_k^{0,\,2}\wedge\widetilde{f} + \bar\partial\bigg(\frac{1}{w-2\varepsilon}\bigg)\wedge\bar\partial_k\widetilde{f} \hspace{6ex} \mbox{on}\hspace{1ex} \{|w|\leq\varepsilon/\sqrt{3}\},$$

\noindent where $\bar\partial(1/(w-2\varepsilon))$ is a $(0,\,1)$-current supported on $Y_\varepsilon$ (since $w$ is holomorphic).

\end{Lem}

\noindent {\it Proof.} The function $1/(w-2\varepsilon)$ is ${\cal C}^\infty$ on $X\setminus Y_\varepsilon$ and $L^{2-\delta}_{loc}$ on $X$ for every $0<\delta<1$, so this induces the stated regularity for $v_\varepsilon$, while immediate computations yield the stated formulae for $\bar\partial_k v_\varepsilon$ since $\bar\partial_k^2=-2\pi i\alpha_k^{0,\,2}$. \hfill $\Box$

\vspace{2ex}

\noindent {\bf (C) Applying the main a priori twisted $L^2$-estimate.} We will now assume, until further notice, that $k\in\Sigma$ and $\varepsilon>0$ are mutually independent (cf. Definition \ref{Def:constrained}). We will see later how our results transform in the special case when $\varepsilon$ is constrained by $k$. 



\vspace{1ex}

$\bullet$ Using Theorem \ref{The:main-a-priori-estimate} and the $L^2_\omega$-orthogonal splitting (\ref{eqn:n0_3-space-decomp_epsilon-k}) featuring therein, we decompose \begin{equation}\label{eqn:1-w_f-tilde_epsilon_splitting}s:=\frac{1}{\sqrt{\eta_\varepsilon + \lambda_\varepsilon}}\,\frac{1}{w-2\varepsilon}\,\widetilde{f}_\varepsilon = s_{h,\,\varepsilon} + t_\varepsilon + s_{nh,\,\varepsilon},\end{equation} with $s_{h,\,\varepsilon}\in{\cal H}^{n,\,0}_{k,\,\varepsilon}(X,\,L_k)$, $t_\varepsilon\in{\cal M}^{n,\,0}_{k,\,\varepsilon}(X,\,L_k) $, $s_{nh,\,\varepsilon}\in{\cal N}^{n,\,0}_{k,\,\varepsilon}(X,\,L_k)$  and we infer from (\ref{eqn:main-a-priori-estimate_rough}) the following estimates: \begin{eqnarray}\label{eqn:s_nh-epsilon_estimate}\nonumber ||t_\varepsilon||^2 & \leq & \frac{k^{1+\delta}}{c}\,||\bar\partial_k^{\eta_\varepsilon + \lambda_\varepsilon}s||^2 = \frac{k^{1+\delta}}{c}\,\bigg|\bigg|\frac{1}{w-2\varepsilon}\,\bar\partial_k\widetilde{f}_\varepsilon\bigg|\bigg|^2\\
||s_{nh,\,\varepsilon}||^2 & \leq & \frac{1}{(1-\varepsilon_0)\,k}\,||\bar\partial_k^{\eta_\varepsilon + \lambda_\varepsilon}s||^2 = \frac{1}{(1-\varepsilon_0)\,k}\,\bigg|\bigg|\frac{1}{w-2\varepsilon}\,\bar\partial_k\widetilde{f}_\varepsilon\bigg|\bigg|^2.\end{eqnarray}


Note that estimates (\ref{eqn:s_nh-epsilon_estimate}) also hold if the $L^2$-norms $||\,\,||$ are replaced by any Sobolev norm $||\,\,||_{W^l(X)}$ with $l\geq 0$.

\begin{Def}\label{Def:F_epsilon_def} We set: \begin{eqnarray*}F_\varepsilon = F_{k,\,\varepsilon} & := & (w-2\varepsilon)\,\,\bigg(\frac{1}{w-2\varepsilon}\,\widetilde{f}_\varepsilon - \sqrt{\eta_\varepsilon + \lambda_\varepsilon}\,t_\varepsilon - \sqrt{\eta_\varepsilon + \lambda_\varepsilon}\,s_{nh,\,\varepsilon}\bigg) \\
    & = & \widetilde{f}_\varepsilon - (w-2\varepsilon)\, \sqrt{\eta_\varepsilon + \lambda_\varepsilon}\,(t_\varepsilon + s_{nh,\,\varepsilon}) \\
    & = & (w-2\varepsilon)\, \sqrt{\eta_\varepsilon + \lambda_\varepsilon}\,s_{h,\,\varepsilon}\in {\cal C}^\infty_{n,\,0}(X,\,L_k).\end{eqnarray*}

\end{Def}

\vspace{1ex}

$\bullet$ This choice of $F_\varepsilon$ can be justified by analogy with the integrable case. If we had $\bar\partial_k^2=0$, we would pick, for every $\varepsilon>0$, the minimal $L^2$-norm solution $u_\varepsilon$ of the equation \begin{eqnarray}\label{eqn:integrable-case_analogy_equation}\bar\partial_k^{\eta_\varepsilon+\lambda_\varepsilon}u_\varepsilon = \bar\partial_k^{\eta_\varepsilon+\lambda_\varepsilon}\bigg(\frac{1}{\sqrt{\eta_\varepsilon + \lambda_\varepsilon}}\,\frac{1}{w-2\varepsilon}\,\widetilde{f}_\varepsilon\bigg).\end{eqnarray} This equation is equivalent to each of the following two equations: \begin{eqnarray*}\bar\partial_k\bigg(\frac{1}{w-2\varepsilon}\,\widetilde{f}_\varepsilon - \sqrt{\eta_\varepsilon + \lambda_\varepsilon}\,u_\varepsilon\bigg) = 0 \iff \bar\partial_k\bigg(\widetilde{f}_\varepsilon - (w-2\varepsilon)\,\sqrt{\eta_\varepsilon + \lambda_\varepsilon}\,u_\varepsilon\bigg) = 0,\end{eqnarray*} where the last equivalence is a consequence of the function $w\longmapsto w-2\varepsilon$ being holomorphic on $X$.

We would then set $F_\varepsilon:=\widetilde{f}_\varepsilon - (w-2\varepsilon)\,\sqrt{\eta_\varepsilon + \lambda_\varepsilon}\,u_\varepsilon$. Each $L_k$-valued $(n,\,0)$-form $F_\varepsilon$ would be holomorphic, by construction. The sought-after extension $F$ would be the limit (whose existence in a certain topology we would deduce from the $L^2$-estimate we had previously obtained) of a subsequence $F_{\varepsilon_j}$ with $\varepsilon_j\downarrow 0$.

In our (possibly) non-integrable context, the part of $u_\varepsilon$ is played by $t_\varepsilon + s_{nh,\,\varepsilon}$. Both $u_\varepsilon$ and $t_\varepsilon + s_{nh,\,\varepsilon}$ represent, each in its own case, the minimal $L^2$-norm correction of $\widetilde{f}_\varepsilon$ to a holomorphic, resp. approximately holomorphic, $L_k$-valued $(n,\,0)$-form.

This non-integrable analogue of the classical integrable $\bar\partial$-equation resolution was introduced in [Pop13] for non-twisted Laplacians. The new feature of the present case is the twisting of $\bar\partial_k$ to $\bar\partial_k^{\eta_\varepsilon+\lambda_\varepsilon}$ and the use of the associated twisted Laplacians $\widetilde\Delta''_{k,\,(\varepsilon)}$ and $\Delta''_{k,\,(\varepsilon)}$ in bidegree $(n,\,0)$, resp. $(n,\,1)$.

\vspace{1ex}

$\bullet$ To estimate the $L^2$-norm of $F_\varepsilon$, we need to estimate the r.h.s. of (\ref{eqn:s_nh-epsilon_estimate}). Using the expression for $v_\varepsilon= \bar\partial_k((1/(w-2\varepsilon))\,\widetilde{f}_\varepsilon)$ on $X\setminus Y_\varepsilon$ given in Lemma \ref{Lem:v_epsilon}, the squared pointwise norm splits as: \begin{eqnarray}\label{eqn:B-epsilon_inverse1}\nonumber|v_\varepsilon|^2 = \bigg|\frac{1}{w-2\varepsilon}\,\bar\partial_k\widetilde{f}_\varepsilon\bigg|^2 & = & \frac{1}{\varepsilon^4}\,\bigg|\frac{w}{w-2\varepsilon}\bigg|^2 \,\bigg|\rho'\bigg(\frac{|w|^2}{\varepsilon^2}\bigg)\bigg|^2\,|d\overline{w}\wedge\widetilde{f}|^2 + \rho'\bigg(\frac{|w|^2}{\varepsilon^2}\bigg)\,\frac{w\,\rho_\varepsilon(w)}{\varepsilon^2\,|w-2\varepsilon|^2}\,\langle d\overline{w}\wedge\widetilde{f},\,\bar\partial_k\widetilde{f}\rangle\\
& + & \rho'\bigg(\frac{|w|^2}{\varepsilon^2}\bigg)\,\frac{\overline{w}\,\rho_\varepsilon(w)}{\varepsilon^2\,|w-2\varepsilon|^2}\,\langle \bar\partial_k\widetilde{f},\,d\overline{w}\wedge\widetilde{f}\rangle +  \frac{|\rho_\varepsilon(w)|^2}{|w-2\varepsilon|^2}\,|\bar\partial_k\widetilde{f}|^2.\end{eqnarray}







Now, $\rho'(|w|^2/\varepsilon^2)$ is supported on $\{\varepsilon/\sqrt{3}\leq|w|\leq\varepsilon\}\subset X$ and $|\rho'(|w|^2/\varepsilon^2)|^2\leq 4$, while $\rho_\varepsilon$ is supported on $\{|w|\leq\varepsilon\}\subset X$ and $0\leq\rho_\varepsilon(w)\leq 1$. Thus, putting together (\ref{eqn:s_nh-epsilon_estimate}), (\ref{eqn:B-epsilon_inverse1}) and the general inequality $|u\wedge v|\leq|u|\,|v|$ for any forms $u,\,v$, we get

\begin{Lem}\label{Lem:s_nh_epsilon_estimate} The sections $s_{nh,\,\varepsilon}$ and $t_\varepsilon$ defined in (\ref{eqn:1-w_f-tilde_epsilon_splitting}) satisfy the following estimates: \begin{eqnarray}\nonumber ||s_{nh,\,\varepsilon}||^2 & \leq & \frac{4}{(1-\varepsilon_0)\,k}\,\int\limits_{\frac{\varepsilon}{\sqrt{3}}\leq|w|\leq\varepsilon}\frac{1}{\varepsilon^4}\,\bigg|\frac{w}{w-2\varepsilon}\bigg|^2\,|dw|^2_\omega\,|\widetilde{f}|^2_{h_k,\,\omega}\,dV_\omega \\
\nonumber & + & \frac{4}{(1-\varepsilon_0)\,k}\,\int\limits_{\frac{\varepsilon}{\sqrt{3}}\leq|w|\leq\varepsilon}\frac{1}{\varepsilon^2}\,\frac{|w|\,|dw|_\omega}{|w-2\varepsilon|^2}\,|\widetilde{f}|_{h_k,\,\omega}\,|\bar\partial_k\widetilde{f}|_{h_k,\,\omega}\,dV_\omega + \frac{1}{(1-\varepsilon_0)\,k}\,\int\limits_{|w|\leq\varepsilon}\frac{1}{|w-2\varepsilon|^2}\,|\bar\partial_k\widetilde{f}|^2\,dV_\omega  \\
\nonumber & =: &  \frac{4}{(1-\varepsilon_0)\,k}\,I_1(\varepsilon) + \frac{4}{(1-\varepsilon_0)\,k}\,I_2(\varepsilon) + \frac{1}{(1-\varepsilon_0)\,k}\,I_3(\varepsilon)\end{eqnarray}

\noindent and

\begin{eqnarray*}||t_\varepsilon||^2\leq\frac{4}{c}\,k^{1+\delta}\,\bigg(I_1(\varepsilon) + I_2(\varepsilon) + \frac{1}{4}\,I_3(\varepsilon)\bigg)\end{eqnarray*} for all $k\in\Sigma$ and $\varepsilon>0$.

\end{Lem}

We will now proceed to estimate separately each of the integrals $I_1(\varepsilon)$, $I_2(\varepsilon)$, $I_3(\varepsilon)$.

\vspace{2ex}

\noindent $\bullet$ {\bf Preliminaries.} 

\vspace{1ex}

 With respect to a fixed ${\cal C}^\infty$ frame $e_k$ of $L_{k|U}\simeq U\times\C$, where $U\subset X$ is a trivialising open subset for $L_k$ about a point of $Y$, we can write $$\widetilde{f}=\widetilde{u}\otimes e_k \hspace{1ex} \mbox{on} \hspace{1ex} U \hspace{3ex} \mbox{and} \hspace{3ex} f=u\otimes e_{k|Y} \hspace{1ex} \mbox{on} \hspace{1ex} U\cap Y,$$

 \noindent where $\widetilde{u}$ is a $\C$-valued $(n,\,0)$-form on $U$ and $u$ is a $\C$-valued $(n-1,\,0)$-form on $U\cap Y$. Since $\widetilde{f}_{|Y} = f\wedge dw$ on $Y$ (in the sense that $\widetilde{f}(z,\,0) = f(z)\wedge dw$ for every $z\in Y$), we get $\widetilde{u}_{|Y} = u\wedge dw$ on $U\cap Y$, i.e. $\widetilde{u}(z,\,0) = u(z)\wedge dw$ for every $z\in U\cap Y$.

 Moreover, if we denote by $z=(z_1,\dots , z_{n-1})$ a system of holomorphic coordinates on $U\cap Y$, we have $u(z) = f_U(z)\,dz_1\wedge\dots\wedge dz_{n-1}$ on $U\cap Y$. If $\pi:U\rightarrow U\cap Y$ denotes the holomorphic projection $(z_1,\dots , z_{n-1},\,w)\mapsto(z_1,\dots , z_{n-1})$, the way the {\it rough extension} $\widetilde{f}$ of $f$ was constructed in (A) above means that \begin{eqnarray}\label{eqn:u-tilde_u_dw}\widetilde{u}(z,\,w) = (f_U\circ\pi)(z,\,w)\,dz_1\wedge\dots\wedge dz_{n-1}\wedge dw = f_U(z)\,dz_1\wedge\dots\wedge dz_{n-1}\wedge dw = u(z)\wedge dw\end{eqnarray} for all $(z,\,w)\in U$ (including those for which $w\neq 0$).

\vspace{2ex}

On the other hand, recall that any $(0,\,n)$-form is primitive and for any {\it primitive} $(p,\,q)$-form $v$ with $p+q=k\leq n$, the Hodge star operator is given by the standard formula (see e.g. [Voi02, Proposition 6.29, p.  150]): \begin{eqnarray}\label{eqn:stardard-formula_primitive}\star v = (-1)^{k(k+1)/2}i^{p-q}\,v\wedge\omega^{n-k}/(n-k)!\end{eqnarray} In our case, we have: $|\widetilde{u}|^2_\omega\,dV_\omega = \widetilde{u}\wedge\star\overline{\widetilde{u}}$ (a consequence of the definition of the Hodge star operator $\star$) and $\star\overline{\widetilde{u}}=i^{n^2}\,\overline{\widetilde{u}}$ (a consequence of the standard formula (\ref{eqn:stardard-formula_primitive}) applied to the $(0,\,n)$-form $v=\overline{\widetilde{u}}$). Thus, \begin{eqnarray}\label{eqn:u-tilde_squared-norm}i^{n^2}\,\widetilde{u}\wedge\overline{\widetilde{u}} = |\widetilde{u}|^2_\omega\,dV_\omega,  \hspace{5ex} \mbox{on}\hspace{1ex} U\subset X.\end{eqnarray} The analogue of this on the $(n-1)$-dimensional $Y$ is \begin{eqnarray}\label{eqn:u_squared-norm}i^{(n-1)^2}\,u\wedge\bar{u} = |u|_\omega^2\,dV_{Y,\,\omega}  \hspace{5ex} \mbox{on}\hspace{1ex} U\cap Y\subset Y,\end{eqnarray} where $dV_{Y,\,\omega}:=(\omega_{|Y})^{n-1}/(n-1)!$ is the volume form induced on $Y$ by the restriction of $\omega$.

\vspace{2ex}

From the shape of $\widetilde{f}$ and from (\ref{eqn:u-tilde_squared-norm}) we get: \begin{eqnarray}\label{eqn:f-tilde_local_norm_U}|\widetilde{f}|^2_{h_k,\,\omega}\,dV_\omega = |e_k|_{h_k}^2\,|\widetilde{u}|^2_\omega\,dV_\omega = |e_k|_{h_k}^2\,i^{n^2}\,\widetilde{u}\wedge\overline{\widetilde{u}} \hspace{5ex} \mbox{on}\hspace{1ex} U.\end{eqnarray}

\noindent Then, from $\widetilde{u}(z,\,w) = u(z)\wedge dw$ (see (\ref{eqn:u-tilde_u_dw})) we get: \begin{eqnarray*}i^{n^2}\,\widetilde{u}(z,\,w)\wedge\overline{\widetilde{u}}(z,\,w) = \bigg(i^{(n-1)^2}\,u(z)\wedge\bar{u}(z)\bigg)\wedge(idw\wedge d\overline{w}), \hspace{5ex} (z,\,w)\in U,\end{eqnarray*} so (\ref{eqn:f-tilde_local_norm_U}) reads: \begin{eqnarray}\label{eqn:f-tilde_local_norm_U-Y}\nonumber|\widetilde{f}(z,\,w)|^2_{h_k,\,\omega}\,dV_\omega(z,\,w) & = & |e_k(z,\,w)|_{h_k}^2\,\bigg(i^{(n-1)^2}\,u(z)\wedge\bar{u}(z)\bigg)\wedge(idw\wedge d\overline{w})\\
  & = & |e_k(z,\,w)|_{h_k}^2\,\bigg(|u(z)|_\omega^2\,dV_{Y,\,\omega}(z)\bigg)\wedge(idw\wedge d\overline{w}), \hspace{5ex}  (z,\,w)\in U,\end{eqnarray} where (\ref{eqn:u_squared-norm}) was used to get the last equality.

\vspace{2ex}

\noindent $\bullet$ {\bf Estimating $I_1(\varepsilon)$ (same method as in the classical integrable case).} 

\vspace{1ex}

Using (\ref{eqn:f-tilde_local_norm_U-Y}), we get the second equality below: \begin{eqnarray*}I_1(\varepsilon) & = & \int\limits_{\frac{\varepsilon}{\sqrt{3}}\leq|w|\leq\varepsilon}\frac{1}{\varepsilon^4}\,\,\bigg|\frac{w}{w-2\varepsilon}\bigg|^2\,|dw|^2_\omega\,|\widetilde{f}(z,\,w)|^2_{h_k,\,\omega}\,dV_\omega(z,\,w) \\
  & = & \int\limits_{\frac{\varepsilon}{\sqrt{3}}\leq|w|\leq\varepsilon}\frac{1}{\varepsilon^4}\,\,\bigg|\frac{w}{w-2\varepsilon}\bigg|^2\,|dw|^2_\omega\,|e_k(z,\,w)|_{h_k}^2\,\bigg(|u(z)|_\omega^2\,dV_{Y,\,\omega}(z)\bigg)\wedge(idw\wedge d\overline{w}) \\
  & = & \int\limits_{\frac{1}{\sqrt{3}}\leq|\zeta|\leq 1}\frac{|\zeta|^2\,|d\zeta|^2_\omega}{|\zeta-2|^2}\,|e_k(z,\,\varepsilon\,\zeta)|_{h_k}^2\,\bigg(|u(z)|_\omega^2\,dV_{Y,\,\omega}(z)\bigg)\wedge(id\zeta\wedge d\overline\zeta),\end{eqnarray*} where for the last equality we used the change of variable $w=\varepsilon\,\zeta$.

Now, letting $\varepsilon\downarrow 0$ and applying the Fubini Theorem, we get: \begin{eqnarray}\label{eqn:lim_I_1_epsilon}\nonumber\lim\limits_{\varepsilon\rightarrow 0}I_1(\varepsilon) & = & \bigg(\int\limits_{\frac{1}{\sqrt{3}}\leq|\zeta|\leq 1}\frac{|\zeta|^2\,|d\zeta|^2_\omega}{|\zeta-2|^2}\,id\zeta\wedge d\overline\zeta \bigg)\,\bigg(\int\limits_Y|u(z)|_\omega^2\,|e_k(z,\,0)|_{h_k}^2\,dV_{Y,\,\omega}(z)\bigg)\\
 & = & c_1\,\int\limits_Y|f|^2_{h_k,\,\omega}\,dV_{Y,\,\omega},\end{eqnarray} where the constant $c_1>0$ is defined by the integral over $\{\zeta\in\C\,\mid\,1/\sqrt{3}\leq|\zeta|\leq 1\}\subset\C$ of the expression depending on $\zeta$ on the previous line and the last equality is a consequence of $|f(z)|_{h_k,\,\omega} = |u(z)|_\omega\,|e_k(z,\,0)|_{h_k}$ for all $z\in Y$.

\vspace{2ex}

\noindent $\bullet$ {\bf Estimating $I_2(\varepsilon)$.} We have: \begin{eqnarray}\label{eqn:I_2_varepsilon_rel1}I_2(\varepsilon) & = & \int\limits_{\frac{\varepsilon}{\sqrt{3}}\leq|w|\leq\varepsilon}\frac{1}{\varepsilon^2}\,\frac{|w|\,|dw|_\omega}{|w-2\varepsilon|^2}\,|\widetilde{f}|_{h_k,\,\omega}\,|\bar\partial_k\widetilde{f}|_{h_k,\,\omega}\,dV_\omega.\end{eqnarray}

Now, for all $(z,\,w)\in U$, we have: \begin{eqnarray*}\widetilde{f}(z,\,w) = \widetilde{u}(z,\,w)\otimes e_k(z,\,w) = \bigg(u(z)\otimes e_k(z,\,w)\bigg)\wedge dw,\end{eqnarray*} where we used (\ref{eqn:u-tilde_u_dw}) to get the latter equality. Taking $\bar\partial_k$, we get: \begin{eqnarray}\label{eqn:d-bar-_k_f-tilde_local_norm_U-Y}(\bar\partial_k\widetilde{f})(z,\,w) = \bigg((\bar\partial u)(z)\otimes e_k(z,\,w) + (-1)^{n-1}\,u(z)\otimes(\bar\partial_k e_k)(z,\,w)\bigg)\wedge dw\end{eqnarray} for $(z,\,w)\in U$.

  On the other hand, in the complex $1$-dimensional direction of $w$, the standard formula (\ref{eqn:stardard-formula_primitive}) applied to the (necessarily primitive) $(0,\,1)$-form $d\overline{w}$ in this $1$-dimensional ambient space (so, we take $n=1$ in the standard formula) spells: $\star(d\overline{w}) = idw$. Denoting by $dV_{Y^\perp,\,\omega}$ the volume form (a $(1,\,1)$-form) induced by $\omega$ in the direction transversal to $Y$ ($=$ the direction of $w$), we have: \begin{eqnarray*}idw\wedge d\overline{w} = dw\wedge\star(d\overline{w}) = |dw|^2_\omega\,dV_{Y^\perp,\,\omega}.\end{eqnarray*} This leads to \begin{eqnarray}\label{eqn:dV_omega_local}dV_\omega(z,\,w) = dV_{Y,\,\omega}(z)\wedge dV_{Y^\perp,\,\omega}(w) = \frac{1}{|dw|^2_\omega}\,dV_{Y,\,\omega}(z)\wedge(idw\wedge d\overline{w})\end{eqnarray} for $(z,\,w)\in U$. 

  Thus, using (\ref{eqn:f-tilde_local_norm_U-Y}), (\ref{eqn:d-bar-_k_f-tilde_local_norm_U-Y}) and the general inequality $|u\wedge v|\leq|u|\,|v|$ for forms $u,v$, from (\ref{eqn:I_2_varepsilon_rel1}) we get: \begin{eqnarray*}I_2(\varepsilon) \leq \widetilde{I}_2(\varepsilon):=\int\limits_{\frac{\varepsilon}{\sqrt{3}}\leq|w|\leq\varepsilon}\frac{1}{\varepsilon^2}\,\frac{|w|\,|dw|_\omega}{|w-2\varepsilon|^2}\,\bigg|(\bar\partial u)(z)\otimes e_k(z,\,w) & + & (-1)^{n-1}\,u(z)\otimes(\bar\partial_k e_k)(z,\,w)\bigg|_{h_k,\,\omega}\cdot \\
    &  & |dw|_\omega\,|e_k(z,\,w)|_{h_k}\,|u(z)|_\omega\,|dw|_\omega\,dV_\omega(z,\,w).\end{eqnarray*}  

 Now, two of the factors $|dw|_\omega$ in the above integral cancel out the factor $1/|dw|^2_\omega$ in the expression (\ref{eqn:dV_omega_local}) for $dV_\omega(z,\,w)$. Thus, after the change of variable $w=\varepsilon\zeta$, we get: \begin{eqnarray*}I_2(\varepsilon) \leq \int\limits_{\frac{1}{\sqrt{3}}\leq|\zeta|\leq 1}\frac{|\zeta|\,|d\zeta|_\omega}{|\zeta-2|^2} \bigg|(\bar\partial u)(z)\otimes e_k(z,\,\varepsilon\zeta) & + & (-1)^{n-1}\,u(z)\otimes(\bar\partial_k e_k)(z,\,\varepsilon\zeta)\bigg|_{h_k,\,\omega}\cdot \\
    & & |e_k(z,\,\varepsilon\zeta)|_{h_k}\,|u(z)|_\omega\,dV_{Y,\,\omega}(z)\wedge(id\zeta\wedge d\overline{\zeta}).\end{eqnarray*}

  Now, letting $\varepsilon\downarrow 0$ and applying the Fubini Theorem, we get: \begin{eqnarray*}\nonumber\lim\limits_{\varepsilon\rightarrow 0}I_2(\varepsilon) & \leq & \bigg(\int\limits_{\frac{1}{\sqrt{3}}\leq|\zeta|\leq 1}\frac{|\zeta|\,|d\zeta|_\omega}{|\zeta-2|^2}\,id\zeta\wedge d\overline{\zeta}\bigg)\cdot\\
    & & \bigg(\int\limits_Y\bigg|(\bar\partial u)(z)\otimes e_k(z,\,0) + (-1)^{n-1}\,u(z)\otimes(\bar\partial_k e_k)(z,\,0)\bigg|_{h_k,\,\omega}\,|u(z)|_\omega\,|e_k(z,\,0)|_{h_k}\,dV_{Y,\,\omega}(z)\bigg).\end{eqnarray*} Since $(\bar\partial_kf)(z) = (\bar\partial u)(z)\otimes e_k(z,\,0) + (-1)^{n-1}\,u(z)\otimes(\bar\partial_k e_k)(z,\,0)$ and $|f(z)|_{h_k,\,\omega} = |u(z)|_\omega\,|e_k(z,\,0)|_{h_k}$ for all $z\in Y$, this means that \begin{eqnarray}\label{eqn:lim_I_2_epsilon}\lim\limits_{\varepsilon\rightarrow 0}I_2(\varepsilon) \leq c_2\,\int\limits_Y|f|_{h_k,\,\omega}\,|\bar\partial_kf|_{h_k,\,\omega}\,dV_{Y,\,\omega},\end{eqnarray} where the constant $c_2>0$ is defined by the integral over $\{\zeta\in\C\,\mid\,1/\sqrt{3}\leq|\zeta|\leq 1\}\subset\C$ of the expression depending on $\zeta$ in the previous equality for $\lim_{\varepsilon\rightarrow 0}I_2(\varepsilon)$.

\vspace{2ex}

\noindent $\bullet$ {\bf Estimating $I_3(\varepsilon)$.} We have: \begin{eqnarray*}\nonumber I_3(\varepsilon) = \int\limits_{|w|\leq\varepsilon}\frac{1}{|w-2\varepsilon|^2}\,|\bar\partial_k\widetilde{f}|_{h_k,\,\omega}^2\,dV_\omega.\end{eqnarray*} Now, using the expression (\ref{eqn:d-bar-_k_f-tilde_local_norm_U-Y}) for $(\bar\partial_k\widetilde{f})(z,\,w)$, the general inequality $|u\wedge v|\leq|u|\,|v|$ for any forms $u,v$ and formula (\ref{eqn:dV_omega_local}) for $dV_\omega(z,\,w)$, we get: \begin{eqnarray*}I_3(\varepsilon)\leq\widetilde{I}_3(\varepsilon):=\int\limits_{|w|\leq\varepsilon}\frac{1}{|w-2\varepsilon|^2}\,\bigg|(\bar\partial u)(z)\otimes e_k(z,\,w) & + & (-1)^{n-1}\,u(z)\otimes(\bar\partial_k e_k)(z,\,w)\bigg|_{h_k,\,\omega}^2\,|dw|^2_\omega\cdot\\
  & & \frac{1}{|dw|^2_\omega}\,dV_{Y,\,\omega}(z)\wedge(idw\wedge d\bar{w}).\end{eqnarray*}

Thus, after the change of variable $w=\varepsilon\zeta$, we get: \begin{eqnarray*}I_3(\varepsilon)\leq\int\limits_{|\zeta|\leq 1}\frac{1}{|\zeta-2|^2}\,\bigg|(\bar\partial u)(z)\otimes e_k(z,\,\varepsilon\zeta) + (-1)^{n-1}\,u(z)\otimes(\bar\partial_k e_k)(z,\,\varepsilon\zeta)\bigg|_{h_k,\,\omega}^2 \,dV_{Y,\,\omega}(z)\wedge(id\zeta\wedge d\bar{\zeta}).\end{eqnarray*}

Now, letting $\varepsilon\downarrow 0$ and applying the Fubini Theorem, we get:   \begin{eqnarray*}\label{eqn:I_3_varepsilon_rel1}\lim\limits_{\varepsilon\rightarrow 0}I_3(\varepsilon)\leq \bigg(\int\limits_{|\zeta|\leq 1}\frac{1}{|\zeta-2|^2}\,id\zeta\wedge d\bar{\zeta}\bigg)\, \int\limits_Y\bigg|(\bar\partial u)(z)\otimes e_k(z,\,0) + (-1)^{n-1}\,u(z)\otimes(\bar\partial_k e_k)(z,\,0)\bigg|_{h_k,\,\omega}^2\,dV_{Y,\,\omega}(z).\end{eqnarray*} Since $(\bar\partial_kf)(z) = (\bar\partial u)(z)\otimes e_k(z,\,0) + (-1)^{n-1}\,u(z)\otimes(\bar\partial_k e_k)(z,\,0)$ for all $z\in Y$, this means that  \begin{eqnarray}\label{eqn:I_3_varepsilon_rel1}\lim\limits_{\varepsilon\rightarrow 0}I_3(\varepsilon)\leq c_3\int\limits_Y|\bar\partial_kf|^2_{h_k,\,\omega}\,dV_{Y,\,\omega},\end{eqnarray} where the constant $c_3>0$ is defined by the integral over $\{\zeta\in\C\,\mid\,|\zeta|\leq 1\}\subset\C$ of the expression depending on $\zeta$ in the previous upper estimate of $\lim_{\varepsilon\rightarrow 0}I_3(\varepsilon)$.

This is where we see why $w-2\varepsilon$ has been used in place of the $w$ used in the integrable case (see e.g. [Siu96] for this latter case): the integral defining $c_3$ would have been the divergent $\int_{|\zeta|\leq 1}(1/|\zeta|^2)\,id\zeta\wedge d\bar{\zeta}$. The integrals $I_2(\varepsilon)$ and $I_3(\varepsilon)$ do not feature in the integrable case since $f$ is holomorphic then.

\vspace{3ex}

Summing up, Lemma \ref{Lem:s_nh_epsilon_estimate} and the above considerations yield the following

\begin{Conc}\label{Conc:s_nh_epsilon_estimate_limit} For all $k\in\Sigma$ and $\varepsilon>0$, the sections $s_{nh,\,\varepsilon}$ and $t_\varepsilon$ defined in (\ref{eqn:1-w_f-tilde_epsilon_splitting}) satisfy the following estimates: \begin{eqnarray*}||s_{nh,\,\varepsilon}||^2 \leq \frac{4}{(1-\varepsilon_0)\,k}\,\bigg(I_1(\varepsilon) + \widetilde{I}_2(\varepsilon) + \frac{1}{4}\,\widetilde{I}_3(\varepsilon)\bigg),\end{eqnarray*} and \begin{eqnarray*}||t_\varepsilon||^2 \leq \frac{4}{c}\,k^{1+\delta}\,\bigg(I_1(\varepsilon) + \widetilde{I}_2(\varepsilon) + \frac{1}{4}\,\widetilde{I}_3(\varepsilon)\bigg),\end{eqnarray*} with $\varepsilon$-dependent quantities $I_1(\varepsilon), \widetilde{I}_2(\varepsilon), \widetilde{I}_3(\varepsilon)>0$ such that \begin{eqnarray*}\lim\limits_{\varepsilon\rightarrow 0}I_1(\varepsilon) & = & c_1\,\int\limits_Y|f|^2_{h_k,\,\omega}\,dV_{Y,\,\omega} \\
\lim\limits_{\varepsilon\rightarrow 0}\widetilde{I}_2(\varepsilon) & = & c_2\,\int\limits_Y|f|_{h_k,\,\omega}\,|\bar\partial_kf|_{h_k,\,\omega}\,dV_{Y,\,\omega} \\
 \lim\limits_{\varepsilon\rightarrow 0}\widetilde{I}_3(\varepsilon) & = & c_3\,\int\limits_Y|\bar\partial_kf|^2_{h_k,\,\omega}\,dV_{Y,\,\omega},\end{eqnarray*} for $k$-independent constants $c_1, c_2, c_3>0$.

\end{Conc}

\vspace{3ex}

Based on this conclusion and the boundedness of $w$, we can now estimate the $L^2$-norms of the $L_k$-valued $(n,\,0)$-forms $F_\varepsilon=\widetilde{f}_\varepsilon - (w-2\varepsilon)\,\sqrt{\eta_\varepsilon + \lambda_\varepsilon}\,(t_\varepsilon + s_{nh,\,\varepsilon})$ introduced in Definition \ref{Def:F_epsilon_def}. To stress the dependence of $F_\varepsilon$ and $\widetilde{f}_\varepsilon$ on $k$, we will henceforth denote: \begin{eqnarray*}F_{k,\,\varepsilon}:=F_\varepsilon \hspace{3ex}\mbox{and}\hspace{3ex} \widetilde{f}_{k,\,\varepsilon}:=\widetilde{f}_\varepsilon.\end{eqnarray*} Although $t_\varepsilon$ and $s_{nh,\,\varepsilon}$ depend on $k$ too, we will omit the $k$-subscript in their case for the sake of notational simplicity. We will distinguish between the two cases described in Definition \ref{Def:constrained} according to whether $k$ and $\varepsilon$ are mutually independent or not.

\begin{Cor}\label{Cor:F_epsilon_L2-norm_estimate} Suppose the assumptions are as in Theorem \ref{The:main-a-priori-estimate}.

  \vspace{1ex}

  (i)\, If $k\in\Sigma$ and $\varepsilon>0$ are {\bf mutually independent}, the following inequality holds: \begin{eqnarray}\label{eqn:F_epsilon_L2-norm_estimate_mutually-ind}||F_{k,\,\varepsilon}||^2 \leq 3||\widetilde{f}_{k,\,\varepsilon}||^2 + 3\cdot 2^6\bigg(\frac{k^{1+\delta}}{c} + \frac{1}{(1-\varepsilon_0)\,k}\bigg)\,\bigg(I_1(\varepsilon) + \widetilde{I}_2(\varepsilon) + \frac{1}{4}\,\widetilde{I}_3(\varepsilon)\bigg).\end{eqnarray}

  \vspace{1ex}

 In particular, for every $k\in\Sigma$, there exists a decreasing sequence of positive reals $\varepsilon_j\to 0$ such that $F_{k,\,\varepsilon_j}$ converges, as $j\to +\infty$, in the weak topology of the $L^2$-space to some $L_k$-valued $L^2$ $(n,\,0)$-form $F_k$ on $X$ whose squared $L^2$-norm satisfies the estimate: \begin{eqnarray}\label{eqn:F_L2-norm_estimate_mutually-ind}\nonumber||F_k||^2 \leq 3\cdot 2^6\bigg(\frac{k^{1+\delta}}{c} + \frac{1}{(1-\varepsilon_0)\,k}\bigg)\,\bigg(c_1\,\int\limits_Y|f|^2_{h_k,\,\omega}\,dV_{Y,\,\omega} & + & c_2\,\int\limits_Y|f|_{h_k,\,\omega}\,|\bar\partial_kf|_{h_k,\,\omega}\,dV_{Y,\,\omega} \\
    & + & \frac{c_3}{4}\,\int\limits_Y|\bar\partial_kf|^2_{h_k,\,\omega}\,dV_{Y,\,\omega}\bigg),\end{eqnarray} where $c_1, c_2, c_3>0$ are universal constants.

\vspace{1ex}

 (ii)\, If $\varepsilon$ is {\bf constrained by $k$}, the following inequality holds: \begin{eqnarray}\label{eqn:F_epsilon_L2-norm_estimate}||F_{k,\,\varepsilon}||^2 \leq 2||\widetilde{f}_{k,\,\varepsilon}||^2 + \frac{2^7}{(1-\varepsilon_0)\,k}\,\bigg(I_1(\varepsilon) + \widetilde{I}_2(\varepsilon) + \frac{1}{4}\,\widetilde{I}_3(\varepsilon)\bigg).\end{eqnarray}

\end{Cor}

\noindent {\it Proof.} (i)\, We have: \begin{eqnarray*}||F_{k,\,\varepsilon}||^2 & = & ||\widetilde{f}_{k,\,\varepsilon} - (w-2\varepsilon)\,\sqrt{\eta_\varepsilon + \lambda_\varepsilon}\,(t_\varepsilon + s_{nh,\,\varepsilon})||^2 \\
  & \leq & 3\,\bigg[||\widetilde{f}_{k,\,\varepsilon}||^2 + \sup\limits_X\bigg(|w-2\varepsilon|^2\,(\eta_\varepsilon + \lambda_\varepsilon)\bigg)\,(||t_\varepsilon||^2 + ||s_{nh,\,\varepsilon}||^2)\bigg].\end{eqnarray*} Since $\eta_\varepsilon\leq\lambda_\varepsilon$ for $\varepsilon$ small enough, we get: $$|w-2\varepsilon|^2\,(\eta_\varepsilon + \lambda_\varepsilon)\leq 2\,\frac{|w-2\varepsilon|^2}{|w|^2 + \varepsilon^2}\leq 4\,\frac{|w|^2 +4\varepsilon^2}{|w|^2 + \varepsilon^2}\leq 16.$$ We get (\ref{eqn:F_epsilon_L2-norm_estimate_mutually-ind}) from these inequalities combined with the upper bound obtained in Conclusion \ref{Conc:s_nh_epsilon_estimate_limit} on $||s_{nh,\,\varepsilon}||^2$. 

Now, $\widetilde{f}$ is $L^2$, so by dominated convergence, $||\widetilde{f}_{k,\,\varepsilon}||^2$ converges to $0$ as $\varepsilon\to 0$. Indeed, we have: \begin{eqnarray*}||\widetilde{f}_{k,\,\varepsilon}||^2 & = & \int\limits_X|\rho_\varepsilon(w)|^2\,|\widetilde{f}|^2\,dV_{X,\,\omega}\\
& \stackrel{(a)}{\leq} & \int\limits_{|w|\leq\varepsilon}|\widetilde{f}|^2\,dV_{X,\,\omega} = \int\limits_X|\widetilde{f}|^2\,\chi_{\{|w|\leq\varepsilon\}}\,dV_{X,\,\omega} \stackrel{\varepsilon\to 0}{\longrightarrow} \int\limits_X|\widetilde{f}|^2\,\chi_Y\,dV_{X,\,\omega} = \int\limits_Y|\widetilde{f}|^2\,dV_{X,\,\omega} = 0,\end{eqnarray*} where $\chi_E$ stands for the characteristic function of a given set $E$, inequality (a) follows from $\rho_\varepsilon$ being supported on $\{|w|\leq\varepsilon\}$ (which shrinks to $Y$ as $\varepsilon\to 0$) and from $|\rho_\varepsilon|\leq 1$, while the last equality follows from $Y$ being Lebesgue negligible.

Thus, (\ref{eqn:F_epsilon_L2-norm_estimate_mutually-ind}) and the last part of Conclusion \ref{Conc:s_nh_epsilon_estimate_limit} imply the uniform boundedness in the $L^2$-norm of the family $(F_{k,\,\varepsilon})_{\varepsilon>0}$. Therefore, by Alaoglu's theorem, we can extract a subsequence $\varepsilon_j\downarrow 0$ such that $(F_{k,\,\varepsilon_j})_j$ converges in the weak topology of the $L^2$-space to some $L_k$-valued $L^2$ $(n,\,0)$-form $F_k$ on $X$. The upper bound (\ref{eqn:F_L2-norm_estimate_mutually-ind}) on the squared $L^2$-norm of $F_k$ follows by taking $\varepsilon = \varepsilon_j$ in (\ref{eqn:F_epsilon_L2-norm_estimate_mutually-ind}) and then passing to the limit as $j\to +\infty$  and by elementary arguments (including the Cauchy-Schwarz inequality $\langle\langle F_{k,\,\varepsilon_j},\,F_k\rangle\rangle\leq||F_{k,\,\varepsilon_j}||\,||F_k||$) starting from the convergence $\langle\langle F_{k,\,\varepsilon_j},\,F_k\rangle\rangle\longrightarrow||F_k||^2$, itself a consequence of the weak convergence $F_{k,\,\varepsilon_j}\longrightarrow F_k$, as $j\to +\infty$.

\vspace{1ex}

(ii)\, When $\varepsilon$ is constrained by $k$, $t_\varepsilon=0$ thanks to (i) of Theorem \ref{The:main-a-priori-estimate}. Estimate (\ref{eqn:F_epsilon_L2-norm_estimate}) follows as (\ref{eqn:F_epsilon_L2-norm_estimate_mutually-ind}) did above, starting from \begin{eqnarray*}||F_{k,\,\varepsilon}||^2 = ||\widetilde{f}_{k,\,\varepsilon} - (w-2\varepsilon)\,\sqrt{\eta_\varepsilon + \lambda_\varepsilon}\,s_{nh,\,\varepsilon}||^2\leq 2\,\bigg(||\widetilde{f}_{k,\,\varepsilon}||^2 + \sup\limits_X\bigg(|w-2\varepsilon|^2\,(\eta_\varepsilon + \lambda_\varepsilon)\bigg)\,||s_{nh,\,\varepsilon}||^2\bigg).\end{eqnarray*}

\hfill $\Box$

\vspace{3ex}

Note that we cannot fix $k$ and let $\varepsilon\downarrow 0$ in (ii) of the above Corollary \ref{Cor:F_epsilon_L2-norm_estimate} (unlike what we did in (i)) since $\varepsilon$ is bounded below by a positive quantity depending on $k$ when it is constrained by $k$.

\vspace{3ex}

\noindent {\it End of proof of Theorem \ref{nonintOT}.} (a) The $L_k$-valued $(n,\,0)$-form $F_k$ on $X$, constructed from the original $f_k\wedge dw$ on $Y$, was obtained in (i) of Corollary \ref{Cor:F_epsilon_L2-norm_estimate} where its $L^2$-norm estimate (\ref{eqn:nonintOT_estimate}) (= its property (iii)) was obtained as (\ref{eqn:F_L2-norm_estimate_mutually-ind}). It remains to prove that $F_k$ has properties (i) and (ii).

\vspace{1ex}

$\bullet$ {\it Proof of property (i) for $F_k$.} Recall that \begin{eqnarray*}F_{k,\,\varepsilon}=\widetilde{f}_{k,\,\varepsilon} - (w-2\varepsilon)\,\sqrt{\eta_\varepsilon + \lambda_\varepsilon}\,(t_\varepsilon + s_{nh,\,\varepsilon}) = (w-2\varepsilon)\,\sqrt{\eta_\varepsilon + \lambda_\varepsilon}\,s_{h,\,\varepsilon},\end{eqnarray*} with $s_{h,\,\varepsilon}\in{\cal H}^{n,\,0}_{k,\,\varepsilon}(X,\,L_k)$. Thus, by Observation \ref{Obs:dbar_s_h-epsilon_estimate}, $s_{h,\,\varepsilon}$ satisfies estimate (\ref{eqn:dbar_s_h-epsilon_estimate}), leading to: \begin{eqnarray*}\frac{c}{k^{1+\delta}}\,||s_{h,\,\varepsilon}||^2 & \geq & ||\bar\partial_k^{\eta_\varepsilon + \lambda_\varepsilon}s_{h,\,\varepsilon}||^2 = ||\bar\partial_k(\sqrt{\eta_\varepsilon + \lambda_\varepsilon}\,s_{h,\,\varepsilon})||^2 = \bigg|\bigg|\bar\partial_k\bigg(\frac{1}{w-2\varepsilon}\,F_{k,\,\varepsilon}\bigg)\bigg|\bigg|^2 \\
  & = & \bigg|\bigg|\frac{1}{w-2\varepsilon}\,\bar\partial_kF_{k,\,\varepsilon}\bigg|\bigg|^2\geq\frac{1}{(1+2\varepsilon)^2}\,||\bar\partial_kF_{k,\,\varepsilon}||^2,\end{eqnarray*} the last inequality following from the hypothesis $\sup\limits_X|w|\leq 1$ yielding $|w-2\varepsilon|\leq|w| + 2\varepsilon\leq 1 +2\varepsilon$.

Since $F_{k,\,\varepsilon} = (w-2\varepsilon)\sqrt{\eta_\varepsilon + \lambda_\varepsilon}\,s_{h,\,\varepsilon} = \widetilde{f}_{k,\,\varepsilon} - (w-2\varepsilon)\sqrt{\eta_\varepsilon + \lambda_\varepsilon}\,(t_\varepsilon + s_{nh,\,\varepsilon})$, we infer: \begin{eqnarray}\label{eqn:nonintOT_proof_1}\nonumber||\bar\partial_kF_{k,\,\varepsilon}||^2 & \leq & \frac{c\,(1+2\varepsilon)^2}{k^{1+\delta}}\,||s_{h,\,\varepsilon}||^2 = \frac{c\,(1+2\varepsilon)^2}{k^{1+\delta}}\,\bigg|\bigg|\frac{1}{(w-2\varepsilon)\sqrt{\eta_\varepsilon + \lambda_\varepsilon}}\,F_{k,\,\varepsilon}\bigg|\bigg|^2 \\
\nonumber & \leq & \frac{4c}{k^{1+\delta}}\,\bigg|\bigg|\frac{1}{(w-2\varepsilon)\sqrt{\eta_\varepsilon + \lambda_\varepsilon}}\,\widetilde{f}_{k,\,\varepsilon} - (t_\varepsilon + s_{nh,\,\varepsilon})\bigg|\bigg|^2 \\
 & \leq & \frac{12c}{k^{1+\delta}}\,\bigg(\bigg|\bigg|\frac{1}{(w-2\varepsilon)\sqrt{\eta_\varepsilon + \lambda_\varepsilon}}\,\widetilde{f}_{k,\,\varepsilon}\bigg|\bigg|^2 + ||t_\varepsilon||^2 + ||s_{nh,\,\varepsilon}||^2\bigg) \\
\nonumber & \leq & \frac{12c}{k^{1+\delta}}\,\bigg[\bigg|\bigg|\frac{1}{(w-2\varepsilon)\sqrt{\eta_\varepsilon + \lambda_\varepsilon}}\,\widetilde{f}_{k,\,\varepsilon}\bigg|\bigg|^2 + \bigg(\frac{4}{c}\,k^{1+\delta} + \frac{4}{(1-\varepsilon_0)\,k}\bigg)\,\bigg(I_1(\varepsilon) + I_2(\varepsilon) + \frac{1}{4}\,I_3(\varepsilon)\bigg)\bigg],\end{eqnarray} whenever $\varepsilon\in(0,\,1/2)$. We used Lemma \ref {Lem:s_nh_epsilon_estimate} to get the last inequality.

As for the first squared $L^2$-norm on the right of the last inequality, we get: \begin{eqnarray}\label{eqn:nonintOT_proof_2}\nonumber\bigg|\bigg|\frac{1}{(w-2\varepsilon)\sqrt{\eta_\varepsilon + \lambda_\varepsilon}}\,\widetilde{f}_{k,\,\varepsilon}\bigg|\bigg|^2 & = & \int\limits_X\frac{1}{|w-2\varepsilon|^2\,(\eta_\varepsilon + \lambda_\varepsilon)}\,|\rho_\varepsilon(w)|^2\,|\widetilde{f}_k|^2\,dV_{X,\,\omega} \\
& \leq & \int\limits_{|w|\leq\varepsilon}\frac{1}{|w-2\varepsilon|^2\,(\eta_\varepsilon + \lambda_\varepsilon)}\,|\widetilde{f}_k|^2\,dV_{X,\,\omega},\end{eqnarray} since $|\rho_\varepsilon(w)|\leq 1$ and $\mbox{Supp}\,\rho_\varepsilon\subset\{|w|\leq\varepsilon\}$. Now, whenever $|w|\leq\varepsilon$, we have: $|w-2\varepsilon|\geq 2\varepsilon - |w|\geq 2\varepsilon - \varepsilon = \varepsilon$, hence \begin{eqnarray*}\frac{1}{|w-2\varepsilon|^2\,(\eta_\varepsilon + \lambda_\varepsilon)}\leq\frac{|w|^2 + \varepsilon^2}{|w-2\varepsilon|^2} \leq \frac{2\,\varepsilon^2}{\varepsilon^2} = 2,\end{eqnarray*} having used the relations $\eta_\varepsilon + \lambda_\varepsilon\geq\lambda_\varepsilon = 1/(|w|^2 + \varepsilon^2)$ to get the first inequality above. Using this and (\ref{eqn:nonintOT_proof_2}), we get: \begin{eqnarray}\label{eqn:nonintOT_proof_3}\bigg|\bigg|\frac{1}{(w-2\varepsilon)\sqrt{\eta_\varepsilon + \lambda_\varepsilon}}\,\widetilde{f}_{k,\,\varepsilon}\bigg|\bigg|^2 \leq 2\,\int\limits_{|w|\leq\varepsilon}|\widetilde{f}_k|^2\,dV_{X,\,\omega} = 2\,\int\limits_X|\widetilde{f}_k|^2\,\chi_{\{|w|\leq\varepsilon\}}\,dV_{X,\,\omega}\stackrel{\varepsilon\to 0}{\longrightarrow} 0,\end{eqnarray} where the last convergence was already seen in the proof of (i) of Corollary \ref{Cor:F_epsilon_L2-norm_estimate}.

We infer from (\ref{eqn:nonintOT_proof_3}), from the previously seen inequalities $I_j(\varepsilon)\leq\widetilde{I}_j(\varepsilon)$ ($j=1,2$) and from Conclusion \ref{Conc:s_nh_epsilon_estimate_limit}, that the right-hand side of (\ref{eqn:nonintOT_proof_1}) has an upper bound that converges to \begin{eqnarray}\label{eqn:nonintOT_proof_4}\nonumber L:=\frac{12c}{k^{1+\delta}}\,\bigg(\frac{4}{c}\,k^{1+\delta} + \frac{4}{(1-\varepsilon_0)\,k}\bigg)\,\bigg(c_1\,\int\limits_Y|f|^2_{h_k,\,\omega}\,dV_{Y,\,\omega} & + & c_2\,\int\limits_Y|f|_{h_k,\,\omega}\,|\bar\partial_kf|_{h_k,\,\omega}\,dV_{Y,\,\omega} \\
& + & \frac{c_3}{4}\,\int\limits_Y|\bar\partial_kf|^2_{h_k,\,\omega}\,dV_{Y,\,\omega}\bigg)\end{eqnarray} as $\varepsilon\downarrow 0$.

 In particular, from this and from (\ref{eqn:nonintOT_proof_1}), we infer that the family $(\bar\partial_kF_{k,\,\varepsilon})_{\varepsilon>0}$ is uniformly bounded in the $L^2$-norm. Thus, Alaoglu's theorem implies the existence of a subsequence $\varepsilon_j\downarrow 0$ such that $(\bar\partial_kF_{k,\,\varepsilon_j})_j$ converges in the weak topology of the $L^2$-space to some $L_k$-valued $L^2$ $(n,\,1)$-form $G_k$ on $X$. In particular, $\bar\partial_kF_{k,\,\varepsilon_j}$ also converges to $G_k$ in the weak topology of currents (since the latter topology is weaker than the former) as $j\to +\infty$. 

On the other hand, we know from the analogous argument we applied to the family $(F_{k,\,\varepsilon})_{\varepsilon>0}$ in the proof of (i) of Corollary \ref{Cor:F_epsilon_L2-norm_estimate} that, after possibly extracting a subsequence, $F_{k,\,\varepsilon_j}$ converges to $F_k$ in the weak  topology of the $L^2$-space, hence also in the weak topology of currents. Since $\bar\partial_k$ is continuous for the latter topology, we infer that $\bar\partial_kF_{k,\,\varepsilon_j}$ converges in the weak topology of currents to $\bar\partial_kF_k$.

By uniqueness of the limit, we infer that $G_k=\bar\partial_kF_k$, hence $\bar\partial_kF_{k,\,\varepsilon_j}$ converges to $\bar\partial_kF_k$  in the weak  topology of the $L^2$-space. From this we deduce, as in the last part of the proof of (i) of Corollary \ref{Cor:F_epsilon_L2-norm_estimate}, after taking $\varepsilon=\varepsilon_j$ in (\ref{eqn:nonintOT_proof_1}) and letting $j\to +\infty$, that \begin{eqnarray*}||\bar\partial_kF_k||^2\leq L,\end{eqnarray*} where $L>0$ was defined in (\ref{eqn:nonintOT_proof_4}). This proves (i) of part (a) of Theorem \ref{nonintOT}.

\vspace{1ex}

$\bullet$ {\it Proof of property (ii) for $F_k$.} Fix an arbitrary non-negative integer $l$ and consider the $l$-th Sobolev inner product $\langle\langle\,\,,\,\,\rangle\rangle_{W^l}$ and the $l$-th Sobolev norm $||\,\,||_{W^l}$ for $L_k$-valued forms on $X$ in the modified form introduced in (\ref{eqn:Sobolev-product_Laplacian}) with $\Delta''_k$ now replacing $\Delta''$ and further modified through a local procedure as explained just before Observation \ref{Obs:Sobolev_main-a-priori-estimate}. This lighter notation replaces the notation $\langle\langle\,\,,\,\,\rangle\rangle_{\Delta_k^{''l}}$ and $||\,\,||_{\Delta_k^{''l}}$ one would have expected in the light of (\ref{eqn:Sobolev-product_Laplacian}). The analogue of Lemma \ref{Lem:Sobolev-product_adjoints-equal} ensures that $\bar\partial_k$ has the same formal adjoint with respect to $\langle\langle\,\,,\,\,\rangle\rangle_{W^l}$, for every positive integer $l$, as with respect to the $L^2$-inner product $\langle\langle\,\,,\,\,\rangle\rangle$.

On the one hand, we have: \begin{eqnarray*}\bigg\langle\bigg\langle\Delta''_k\bigg((w-2\varepsilon)\sqrt{\eta_\varepsilon+\lambda_\varepsilon}\,s_{h,\,\varepsilon}\bigg),\,(w-2\varepsilon)\sqrt{\eta_\varepsilon+\lambda_\varepsilon}\,s_{h,\,\varepsilon}\bigg\rangle\bigg\rangle_{W^l} = \bigg|\bigg|\bar\partial_k\bigg((w-2\varepsilon)\sqrt{\eta_\varepsilon+\lambda_\varepsilon}\,s_{h,\,\varepsilon}\bigg)\bigg|\bigg|_{W^l}^2 \end{eqnarray*}
\begin{eqnarray}\label{eqn:Sobolev_1} & \stackrel{(a)}{=} & \nonumber \bigg|\bigg|(w-2\varepsilon)\,\bar\partial_k\bigg(\sqrt{\eta_\varepsilon+\lambda_\varepsilon}\,s_{h,\,\varepsilon}\bigg)\bigg|\bigg|_{W^l}^2   \stackrel{(b)}{\leq} 2(1+4\varepsilon^2)\,\bigg|\bigg|\bar\partial_k\bigg(\sqrt{\eta_\varepsilon+\lambda_\varepsilon}\,s_{h,\,\varepsilon}\bigg)\bigg|\bigg|_{W^l}^2 \\
  \nonumber & = & 2(1+4\varepsilon^2)\, \bigg\langle\bigg\langle\Delta''_k\bigg(\sqrt{\eta_\varepsilon+\lambda_\varepsilon}\,s_{h,\,\varepsilon}\bigg),\,\sqrt{\eta_\varepsilon+\lambda_\varepsilon}\,s_{h,\,\varepsilon}\bigg\rangle\bigg\rangle_{W^l} \\
  \nonumber & \stackrel{(c)}{=} & 2(1+4\varepsilon^2)\, \bigg\langle\bigg\langle\widetilde\Delta''_{k,\,(\varepsilon)}\,s_{h,\,\varepsilon},\,s_{h,\,\varepsilon}\bigg\rangle\bigg\rangle_{W^l} \stackrel{(d)}{\leq} 2(1+4\varepsilon^2)\,\delta_k\,||s_{h,\,\varepsilon}||_{W^l}^2 \\
  & \leq & 10\delta_k\,||s_{h,\,\varepsilon}||^2_{W^l},\end{eqnarray} for all $\varepsilon\in(0,\,1)$, where we put $\delta_k=c/k^{1+\delta}$ with the constant $c:=C\,\max\bigg\{\sup\limits_X|dw|^2_\omega,\,1\bigg\}$ of Theorem \ref{The:main-a-priori-estimate} for some constant $C>0$ independent of $k$. Equality (a) follows from the function $w\mapsto w-2\varepsilon$ being holomorphic; (b) follows from $|w|\leq 1$, which implies $|w-2\varepsilon|^2\leq 2(1+4\varepsilon^2)$, and from the vanishing of all the derivatives of $w-2\varepsilon$ except for $(\partial/\partial w)(w-2\varepsilon) = 1$; (c) follows from equality (\ref{eqn:map_Laplacians_ker_proof-extra}); (d) follows from $$s_{h,\,\varepsilon}\in{\cal H}^{n,\,0}_{[0,\,\delta_k],\,\widetilde\Delta''_{k,\,(\varepsilon)}}(X,\,L_k)={\cal H}^{n,\,0}_{k,\,\varepsilon}(X,\,L_k)= \mbox{Im}\,\bigg(\chi_{[0,\,\delta_k]}(\widetilde\Delta''_{k,\,(\varepsilon)})\bigg).$$

\vspace{2ex}

On the other hand, we will apply the classical {\it G\r{a}rding inequality} to the elliptic non-negative self-adjoint differential operator $\Delta''_k$.  

If $X$ is compact, we consider the constant function $\phi=1$ on $X$. If $X$ is not compact, let $K_1\subset K_2\subset X$ be arbitrary compacts and let $\phi$ be a $C^\infty$ function on $X$ such that \begin{eqnarray*}\phi = 1 \hspace{2ex}\mbox{in a neighbourhood of}\hspace{1ex} K_1, \hspace{5ex} \mbox{Supp}\,\phi\subset\mathring{K}_2, \hspace{5ex} 0\leq\phi\leq 1.\end{eqnarray*} 

For every non-negative integer $l$, the {\it G\r{a}rding inequality} ensures the existence of constants $\delta_1(l),\,\delta_2(l)>0$ depending on $\Delta''_k$, hence on $k$, but not on $\varepsilon$, such that for every $L_k$-valued $(n,\,0)$-form $u$ on $X$, we have: \begin{eqnarray*}\delta_2(l)\,\bigg|\bigg|\phi u\bigg|\bigg|^2_{W^{l+1}}\leq\bigg\langle\bigg\langle\Delta_k''(\phi u),\,\phi u\bigg\rangle\bigg\rangle_{W^l} + \delta_1(l)\,||\phi u||^2.\end{eqnarray*}

Now, choosing $u$ to be the $L_k$-valued $(n,\,0)$-form $F_{k,\,\varepsilon} = (w-2\varepsilon)\,\sqrt{\eta_\varepsilon + \lambda_\varepsilon}\,s_{h,\,\varepsilon}$ in G\r{a}rding's inequality and using (\ref{eqn:Sobolev_1}), we get: \begin{eqnarray}\label{eqn:Garding_consequence}\delta_2(l)\,\bigg|\bigg|F_{k,\,\varepsilon}\bigg|\bigg|^2_{W^{l+1}(K_1)} \leq C'\,\bigg(10\delta_k\,||s_{h,\,\varepsilon}||^2_{W^l} + \delta_1(l)\,||F_{k,\,\varepsilon}||^2\bigg)\end{eqnarray} for some constant $C'>0$ depending on the derivatives of $\phi$ (all of which are bounded). By $||\,\,\,||_{W^l(K_1)}$ we mean the $l$-th Sobolev norm on $K_1$.

Moreover, we can handle the $l$-th Sobolev norm of $s_{h,\,\varepsilon}$ as we did its $L^2$-norm in (\ref{eqn:nonintOT_proof_1}). We get: \begin{eqnarray*}||s_{h,\,\varepsilon}||^2_{W^l} & = & \bigg|\bigg|\frac{1}{(w-2\varepsilon)\sqrt{\eta_\varepsilon + \lambda_\varepsilon}}\,\widetilde{f}_{k,\,\varepsilon} - (t_\varepsilon + s_{nh,\,\varepsilon})\bigg|\bigg|^2_{W^l} \\
  & \leq & 3\,\bigg(\bigg|\bigg|\frac{1}{(w-2\varepsilon)\sqrt{\eta_\varepsilon + \lambda_\varepsilon}}\,\widetilde{f}_{k,\,\varepsilon}\bigg|\bigg|^2_{W^l} + ||t_\varepsilon||^2_{W^l} + ||s_{nh,\,\varepsilon}||^2_{W^l}\bigg)\end{eqnarray*} and the last expression is seen to be bounded independently of $\varepsilon$ when using for the Sobolev norm the arguments that led to (\ref{eqn:nonintOT_proof_3}) and to Conclusion \ref{Conc:s_nh_epsilon_estimate_limit} for the $L^2$-norm. (See also Observation \ref{Obs:Sobolev_main-a-priori-estimate}.)

We infer from this and from (\ref{eqn:Garding_consequence}) that the family $(F_{k,\,\varepsilon})_{\varepsilon>0}$ of $L_k$-valued $(n,\,0)$-forms is uniformly bounded in the $(l+1)$-st Sobolev norm on $K_1$. This implies, thanks to the classical Sobolev Lemma, that for any non-negative integer $r$, if we choose $l$ such that $l>r+n$ (where $n=\mbox{dim}_\C X$), the family $(F_{k,\,\varepsilon})_{\varepsilon>0}$ is uniformly bounded on $K_1$ in the $C^r$-norm. Thus, choosing arbitrarily large $l$'s and compacts $K_1\subset X$, we infer that the family $(F_{k,\,\varepsilon})_{\varepsilon>0}$ and the families of all its derivatives are uniformly bounded on the compacts on $X$. We can, therefore, extract a sequence $\varepsilon_j\downarrow 0$ such that the sequence $(F_{k,\,\varepsilon_j})_j$ converges in the $C^\infty$-topology of $L_k$-valued $(n,\,0)$-forms on $X$.

In particular, the sequence $(F_{k,\,\varepsilon_j})_j$ converges uniformly on the compacts of $X$. Its limit must be the previously found $L_k$-valued $(n,\,0)$-form $F_k$. Since $$\widetilde{f}_{k,\,\varepsilon|Y} = f_k\wedge dw,  \hspace{5ex} \varepsilon>0,$$ (see Lemma \ref{Lem:f-tilda_epsilon} and Conclusion \ref{Conc:rough-extension}), the definition of $F_{k,\,\varepsilon}$ and the uniform convergence on compacts $F_{k,\,\varepsilon_j}\longrightarrow F_k$ as $j\to\infty$ imply that $F_{k|Y} = f_k\wedge dw$. This proves property (ii) 
for $F_k$.

\vspace{2ex}

(b)\, If $f\in{\cal H}^{n-1,\,0}_{[0,\,\delta_k],\,\Delta''_k}(Y,\,L_{k|Y})$, we get the inequality below (see e.g. (i) of Observation \ref{Obs:inclusion_Qs}): \begin{eqnarray*}||\bar\partial_kf||^2_{h_k,\,\omega} = \langle\langle\Delta_k''f,\,f\rangle\rangle_{h_k,\,\omega}\leq\delta_k\,||f||^2_{h_k,\,\omega},\end{eqnarray*} where the equality follows from the fact (due to $f$ being of pure type $(n-1,\,0)$) that $\Delta_k''f = \bar\partial_k^\star\bar\partial_kf$.

This inequality, combined with the Cauchy-Schwarz inequality \begin{eqnarray*}\int\limits_Y|f|_{h_k,\,\omega}\,|\bar\partial_kf|_{h_k,\,\omega}\,dV_{Y,\,\omega}\leq||f||_{h_k,\,\omega}\,||\bar\partial_kf||_{h_k,\,\omega}\end{eqnarray*} leads to the deduction of estimate (\ref{eqn:nonintOT_estimate_bis}) from the estimate (\ref{eqn:nonintOT_estimate}). \hfill $\Box$

\section{Appendix}\label{section:appendix1}

In this section, we give the proofs of some computational results that were used in the text. 

\subsection{Commutation relations}\label{subsection:appendix1}

The first group of them were probably known earlier, but the details were written down in [Pop03, $\S.1.0.2$] and did not appear elsewhere. They can be viewed as commutation relations for zeroth-order operators. 

\begin{Lem}\label{Lem:com} Let $(X,\,\omega)$ be a complex Hermitian manifold. Consider the operators $L:= \omega\wedge \cdot$ and $\Lambda:= L^{\star}$, where $^\star$ denotes (here and throughout this statement) the adjoint w.r.t. the {\bf pointwise} inner product induced on differential forms by the metric $\omega$. Fix any {\bf real-valued} ${\cal C}^\infty$ function $\eta$ on $X$. 

The following identities hold pointwise for arbitrary differential forms of any degree on $X$. \\

\hspace{2ex} (a)\,\,$[\partial\eta\wedge\cdot, \, \Lambda]=i\, (\bar\partial\eta\wedge\cdot)^{\star}$, \hspace{4ex}  $[\bar\partial\eta\wedge\cdot, \, \Lambda] = -i\, (\partial\eta\wedge\cdot)^{\star}$ \\

\hspace{6ex} $[L, \, (\partial\eta\wedge\cdot)^{\star}] = -i\, \bar\partial\eta\wedge\cdot$,  \hspace{4ex} $[(\bar\partial\eta\wedge\cdot)^{\star}, \, L] = -i\, \partial\eta\wedge\cdot$. 

\noindent \begin{eqnarray}\nonumber (b)\,\, [i\, \partial\eta\wedge\bar\partial\eta\wedge\cdot, \, \Lambda] & = & (\partial\eta\wedge\cdot) \, (\partial\eta\wedge\cdot)^{\star} - (\bar\partial\eta\wedge\cdot)^{\star}\, (\bar\partial\eta\wedge\cdot) \\
\nonumber & = & (\bar\partial\eta\wedge\cdot)\, (\bar\partial\eta\wedge\cdot)^{\star} - (\partial\eta\wedge\cdot)^{\star}  \, (\partial\eta\wedge\cdot).\end{eqnarray}

\end{Lem}

\noindent {\it Proof.} (a) It suffices to prove the fourth identity. The remaining three ones follow from it by taking conjugates or adjoints. Thus, we have to prove that for every $(p,q)$-form $v$, we have: $$[(\bar\partial\eta\wedge\cdot)^{\star}, L]v = -i\, \partial\eta \wedge v \iff (\bar\partial\eta\wedge\cdot)^{\star}(\omega \wedge v)-\omega \wedge ((\bar\partial\eta\wedge\cdot)^{\star} \, v) = -i\, \partial\eta \wedge v.$$

 Let $x\in X$ be an arbitrary point and let $z_1, \dots , z_n$ be local holomorphic coordinates about $x$. With respect to this coordinate system, we have: $$(\bar\partial\eta\wedge\cdot)^{\star}v=\sum\limits_j \frac{\partial \eta}{\partial z_j}\,\bigg(\frac{\partial}{\partial \bar{z}_j} \lrcorner v\bigg).$$

\noindent Indeed, this follows from the equalities: \begin{eqnarray*}\langle\!\langle u, (\bar\partial\eta\wedge\cdot)^{\star}v \rangle\!\rangle & = & \langle\!\langle \bar\partial\eta \wedge u, v\rangle\!\rangle = \int \sum\limits_j \langle \frac{\partial \eta}{\partial \bar{z}_j} \, d\bar{z}_j \wedge u, v\rangle \, dV \\
 &= & \int \sum\limits_j \langle d\bar{z}_j \wedge u, \frac{\partial \eta}{\partial z_j}  v\rangle dV = \int \sum\limits_j \langle u, \frac{\partial}{\partial \bar{z}_j} \lrcorner (\frac{\partial \eta}{\partial z_j} \, v) \rangle  dV = \langle\!\langle u, \sum\limits_j \frac{\partial \eta}{\partial z_j} \, (\frac{\partial}{\partial \bar{z}_j} \lrcorner v) \rangle\!\rangle,\end{eqnarray*} that hold for all compactly supported forms $u, v$. \\

 Thus, we get: \begin{eqnarray*}[(\bar\partial \eta\wedge\cdot)^{\star}, L]v & = & \sum\limits_j \frac{\partial \eta}{\partial z_j} \, \frac{\partial}{\partial \bar{z}_j} \lrcorner (\omega \wedge v) - \omega\wedge\bigg(\sum\limits_j\frac{\partial \eta}{\partial z_j}\,\bigg(\frac{\partial}{\partial \bar{z}_j} \lrcorner v\bigg)\bigg)\\
 & = & \sum\limits_j\frac{\partial \eta}{\partial z_j}\,\bigg(\frac{\partial}{\partial \bar{z}_j} \lrcorner \omega\bigg)\wedge v + \sum\limits_j \frac{\partial \eta}{\partial z_j} \, \omega \wedge \bigg(\frac{\partial}{\partial \bar{z}_j} \lrcorner v\bigg) - \sum\limits_j \frac{\partial \eta}{\partial z_j} \, \omega \wedge \bigg(\frac{\partial}{\partial \bar{z}_j} \lrcorner v\bigg) \\
 & = & \sum\limits_j \frac{\partial \eta}{\partial z_j} \, \bigg(\frac{\partial}{\partial \bar{z}_j} \lrcorner \omega\bigg) \wedge v.\end{eqnarray*}

 Since we are trying to prove a pointwise identity for zeroth-order operators, we may assume that $\omega = i \, \sum\limits_k dz_k \wedge d\bar{z}_k$ at the fixed point $x$. Then, we get: $$\frac{\partial}{\partial \bar{z}_j} \lrcorner \omega = -i \, \sum\limits_k dz_k \wedge \bigg(\frac{\partial}{\partial \bar{z}_j} \lrcorner d\bar{z}_k\bigg)= -i\, dz_j.$$

\noindent This proves that $[(\bar\partial\eta)^{\star}, L]v = -i\, \sum\limits_j \frac{\partial \eta}{\partial z_j} \, dz_j \wedge v = -i\, \partial\eta \wedge v$ for all $v$.

\noindent (b) We have: $[i\partial\eta \wedge \bar\partial\eta\wedge\cdot, \Lambda] = i \partial\eta \wedge \bar\partial\eta\wedge\Lambda - i \Lambda \, \partial\eta \wedge \bar\partial\eta\wedge\cdot = $ \\

$ = i\, \partial\eta \wedge (\bar\partial\eta\wedge \Lambda - \Lambda \, \bar\partial\eta\wedge\cdot) + i\, \partial\eta \wedge \Lambda \, \bar\partial\eta\wedge\cdot - i(\Lambda \, \partial\eta -\partial\eta \wedge \Lambda) \, \bar\partial\eta\wedge\cdot -i\, \partial\eta \wedge \Lambda\,(\bar\partial\eta\wedge\cdot)$ \\

$= i\, \partial\eta \wedge [\bar\partial\eta\wedge\cdot, \Lambda] -i\, [\Lambda, \partial\eta\wedge\cdot] \,(\bar\partial\eta\wedge\cdot) = (\partial\eta\wedge\cdot)\,(\partial\eta\wedge\cdot)^{\star} + (\bar\partial\eta\wedge\cdot) \,(\bar\partial\eta\wedge\cdot)^{\star}$. \\

\noindent We have used the identities obtained under (a). \\

\noindent The second identity is equivalent to  \\

$[(\partial\eta, (\partial\eta)^{\star}] = [\bar\partial\eta, (\bar\partial\eta)^{\star}]$. In local coordinates $z_1, \dots ,z_n$, we have\\

\noindent  $[d\bar{z}_j \wedge \cdot ,\,  \frac{\partial}{\partial \bar{z}_k}\lrcorner \cdot] = \delta_{jk}$ and $ [\bar\partial\eta, (\bar\partial\eta)^{\star}] = \sum\limits_j \frac{\partial \eta}{\partial z_j} \, \frac{\partial \eta}{\partial \bar{z}_j}$. By taking conjugates, we get the same expression for $[\partial\eta , (\partial\eta)^{\star}]$, hence the desired identity.      \hfill $\Box$ \\

\vspace{2ex}

The next computation bears again on zeroth-order operators similar to those above.

\begin{Lem}\label{Lem:d-bar_eta_star_no-star} Let $(X,\,\omega)$ be a complex Hermitian manifold with $\mbox{dim}_\C X=n$ and let $\eta$ be a complex-valued ${\cal C}^\infty$ function on $X$. 

For every $(p,\,q)$-form $u$, the following identity holds everywhere on $X$.  

$$(\bar\partial\eta\wedge\cdot)^\star(\bar\partial\eta\wedge u) = |\bar\partial\eta|^2_\omega\,u - \bar\partial\eta\wedge(\bar\partial\eta\wedge\cdot)^\star u  $$

\end{Lem}

\noindent {\it Proof.} We fix an arbitrray point $x\in X$ and choose local holomorphic coordinates $z_1,\dots , z_n $ about $x$ such that $\omega(x) = \sum_{j=1}^n idz_j\wedge d\bar{z}_j$. Then, $(d\bar{z}_j\wedge\cdot)^\star = (\partial/\partial\bar{z}_j)\lrcorner\cdot$ at $x$, so

\vspace{1ex}

\hspace{6ex} $\displaystyle(\bar\partial\eta\wedge\cdot)^\star = \bigg(\sum\limits_{j=1}^n\frac{\partial\eta}{\partial\bar{z}_j}\,d\bar{z}_j\wedge\cdot\bigg)^\star =  \sum\limits_{j=1}^n\frac{\partial\bar\eta}{\partial z_j}\,\frac{\partial}{\partial\bar{z}_j}\lrcorner\cdot  \hspace{3ex} \mbox{at} \hspace{1ex} x.$

\vspace{1ex} 

\noindent Let $u = \sum\limits_{|I|=p,\,|J|=q}u_{I\bar{J}}\,dz_I\wedge d\bar{z}_J$ be the local expression in the chosen coordinates of an arbitrary $(p,\,q)$-form on $X$. We get at $x$: \begin{eqnarray}\nonumber(\bar\partial\eta\wedge\cdot)^\star(\bar\partial\eta\wedge u) & = & \sum\limits_{j,\,k=1}^n\sum\limits_{|I|=p,\,|J|=q}\frac{\partial\bar\eta}{\partial z_j}\,\frac{\partial\eta}{\partial\bar{z}_k}\,u_{I\bar{J}}\,\frac{\partial}{\partial\bar{z}_j}\lrcorner(d\bar{z}_k\wedge dz_I\wedge d\bar{z}_J).\end{eqnarray}

\noindent Since $\frac{\partial}{\partial\bar{z}_j}\lrcorner(d\bar{z}_k\wedge dz_I\wedge d\bar{z}_J) = \delta_{jk}\,dz_I\wedge d\bar{z}_J - d\bar{z}_k\wedge\bigg[\frac{\partial}{\partial\bar{z}_j}\lrcorner\bigg(dz_I\wedge d\bar{z}_J\bigg)\bigg]$, we further get: \begin{eqnarray}\nonumber(\bar\partial\eta\wedge\cdot)^\star(\bar\partial\eta\wedge u) & = & \bigg(\sum\limits_{j=1}^n\bigg|\frac{\partial\bar\eta}{\partial z_j}\bigg|^2\bigg)\,u \\
\nonumber & - & \bigg(\sum\limits_{k=1}^n\frac{\partial\eta}{\partial\bar{z}_k}\,d\bar{z}_k\bigg)\wedge\bigg[\bigg(\sum\limits_{j=1}^n\frac{\partial\bar\eta}{\partial z_j}\,\frac{\partial}{\partial\bar{z}_j}\bigg)\lrcorner\bigg(\sum\limits_{|I|=p,\,|J|=q}u_{I\bar{J}}\,dz_I\wedge d\bar{z}_J\bigg)\bigg].\end{eqnarray}

\noindent Since $\omega$ is given by the identity matrix at $x$, $\sum\limits_{j=1}^n|\frac{\partial\eta}{\partial z_j}|^2 = |\partial\eta|_\omega^2$ at $x$, so the above identity proves the claim.  \hfill $\Box$

\vspace{3ex}

We end this appendix by recalling the standard {\bf Hermitian commutation relations} (see e.g. [Dem97, VII, $\S.1$]).

\begin{Prop}\label{Prop:Hermitian_commutation-rel} Let $(X,\,\omega)$ be a complex Hermitian manifold. Then: \begin{eqnarray}\label{eqn:standard-comm-rel}\nonumber &  & (i)\,\,(\partial + \tau)^{\star} = i\,[\Lambda,\,\bar\partial];  \hspace{3ex} (ii)\,\,(\bar\partial + \bar\tau)^{\star} = - i\,[\Lambda,\,\partial]; \\
&  & (iii)\,\, \partial + \tau = -i\,[\bar\partial^{\star},\,L]; \hspace{3ex} (iv)\,\,
\bar\partial + \bar\tau = i\,[\partial^{\star},\,L],\end{eqnarray}

\noindent where the upper symbol $\star$ stands for the formal adjoint w.r.t. the $L^2$ inner product induced by $\omega$, $L=L_{\omega}:=\omega\wedge\cdot$ is the Lefschetz operator of multiplication by $\omega$, $\Lambda=\Lambda_{\omega}:=L^{\star}$ and $\tau:=[\Lambda,\,\partial\omega\wedge\cdot]$ is the torsion operator (of order zero and type $(1,\,0)$) associated with the metric $\omega$.

\end{Prop}

If $\omega$ is {\it K\"ahler}, $\tau=0$ and relations (\ref{eqn:standard-comm-rel}) reduce to the standard {\it K\"ahler commutation relations}.

\vspace{3ex}

\noindent {\bf References.}\\

\noindent [Dem 97]\, J.-P. Demailly --- {\it Complex Analytic and Algebraic Geometry}---http://www-fourier.ujf-grenoble.fr/~demailly/books.html

\vspace{1ex}

\noindent [Dem01]\, J.-P. Demailly --- {\it Multiplier Ideal Sheaves and Analytic Methods in Algebraic Geometry} --- in {\it Vanishing Theorems and Effective Results in Algebraic Geometry} ed. J.-P. Demailly, L. G\"ottsche, R. Lazarsfeld, ICTP lecture Notes, 6, Abdus Salam International Centre for Theoretical Physics, Trieste, 2001, pp. 1-148.

\vspace{1ex}

\noindent [Fin21]\, S. Finski --- {\it Semiclassical Ohsawa-Takegoshi Extension Theorem and Asymptotics of the Orthogonal Bergman Kernel} --- arXiv e-print DG 2109.06851

\vspace{1ex}

\noindent [GH78]\,P. Griffiths, J. Harris --- {\it Principles of Algebraic Geometry} --- John Wiley \& Sons, Inc., New York, 1994.

\vspace{1ex}

\noindent [Lae02]\, L. Laeng --- {\it Estimations spectrales asymptotiques en g\'eom\'etrie hermitienne} --- Th\`ese de doctorat de l'Universit\'e Joseph Fourier, Grenoble (octobre 2002), http://tel.archives-ouvertes.fr/tel-00002098/en/.

\vspace{1ex}

\noindent [Man93]\, L. Manivel --- {\it Un th\'eor\`eme de prolongement $L^2$ de sections holomorphes d'un fibré hermitien} --- Math. Z. {\bf 212} (1993), 107-122.

\vspace{1ex}

\noindent [OT87] T. Ohsawa, K. Takegoshi --- {\it On The Extension of $L^2$ Holomorphic Functions} --- Math. Z. {\bf 195} (1987) 197-204.

\vspace{1ex}

\noindent [Ohs88] T. Ohsawa --- {\it On The Extension of $L^2$ Holomorphic Functions, II} --- Publ. RIMS, Kyoto Univ., {\bf 24} (1988), 265-275.

\vspace{1ex}

\noindent [Ohs94] T. Ohsawa --- {\it On The Extension of $L^2$ Holomorphic Functions, IV: A New Density Concept} --- Geometry and Analysis on Complex Manifolds (T. Mabuchi et al., eds.), Festschrift for Professor S. Kobayashi’s 60th birthday, Singapore: World Scientific (1994), pp. 157–170.

\vspace{1ex}

\noindent [Ohs95] T. Ohsawa --- {\it On The Extension of $L^2$ Holomorphic Functions, III: Negligible Weights} --- Math. Z. {\bf 219} (1995), 215-225.

\vspace{1ex}

\noindent [Ohs95]\, T. Ohsawa --- {\it On the Extension of $L^2$ Holomorphic Functions, III. Negligible Weights} --- Math. Z. {\bf 219} (1995), 215-225.

\vspace{1ex}

\noindent [Pop03]\, D. Popovici --- {\it Quelques applications des m\'ethodes effectives en g\'eom\'etrie analytique} --- PhD thesis, University Joseph Fourier (Grenoble 1), http://tel.ccsd.cnrs.fr/documents/ 

\vspace{1ex}

\noindent [Pop05]\, D. Popovici --- {\it $L^2$ Extension for Jets of Holomorphic Sections of a Hermitian Line Bundle} --- Nagoya Math. J. {\bf 180} (2005), 1-34. 

\vspace{1ex}

\noindent [Pop13]\, D. Popovici --- {\it Transcendental K\"ahler Cohomology Classes} -- Publ. RIMS Kyoto Univ. {\bf 49} (2013), 313-360.

\vspace{1ex}

\noindent [RS80]\, M. Reed, B. Simon --- {\it Methods of Modern Mathematical Physics; Volume I: Functional Analysis} --- Academic Press, Inc (1980).

\vspace{1ex}

\noindent [RT20]\, S. Rao, I-H. Tsai --- {\it Invariance of Plurigenera and Chow-Type Lemma} --- arXiv 2011.03306.

\vspace{1ex}

\noindent [Siu96]\, Y.-T. Siu --- {\it The Fujita Conjecture and the Extension Theorem of Ohsawa-Takegoshi} --- in {\it Geometric Complex Analysis} ed. Junjiro Noguchi {\it et al}, World Scientific Publishing Co. 1996, pp. 577 - 592.

\vspace{1ex}

\noindent [Siu98]\, Y.-T. Siu --- {\it Invariance of Plurigenera} --- Invent. Math. {\bf 134} (1998), 661-673. 

\vspace{1ex}

\noindent [Siu02]\, Y.-T. Siu --- {\it Extension of Twisted Pluricanonical Sections with Plurisubharmonic Weight and Invariance of Semipositively Twisted Plurigenera for Manifolds Not Necessarily of General Type} --- in {\it Complex Geometry} (G\"ottingen, 2000), 223-277, Springer, Berlin, 2002.

\vspace{1ex}

\noindent [Voi02]\, C. Voisin --- {\it Hodge Theory and Complex Algebraic Geometry. I.} --- Cambridge Studies in Advanced Mathematics, {\bf 76}, Cambridge University Press, Cambridge, 2002.

\vspace{1ex}

\noindent [Wel08]\, R. O. Wells -- {\it Differential Analysis on Complex Manifolds} -- Graduate Texts in Mathematics, {\bf 65}, Springer, 2008.

\vspace{6ex}

\noindent Universit\'e Paul Sabatier, Institut de Math\'ematiques de Toulouse

\noindent 118, route de Narbonne, 31062, Toulouse Cedex 9, France

\noindent Email: popovici@math.univ-toulouse.fr

\end{document}